\documentclass{amsart}

\usepackage{amscd,amsmath,xypic,amssymb,combelow,tikz-cd,etoolbox,calligra,mathrsfs,enumitem,mathtools,epigraph,hyperref,float}

\makeatletter
\patchcmd{\@settitle}{\uppercasenonmath\@title}{}{}{}
\makeatother

\xyoption{all}
\CompileMatrices

\makeatletter
\@addtoreset{equation}{section}
\makeatother

\newtheorem{theorem}[subsection]{Theorem}

\newtheorem{proposition}[subsection]{Proposition}
\newtheorem{lemma}[subsection]{Lemma}

\newtheorem{definition}[subsection]{Definition}

\newtheorem{claim}[subsection]{Claim}

\newtheorem{remark}[subsection]{Remark}

\def\loccit{\emph{loc. cit. }}

\def\fg{{\mathfrak{g}}}

\def\fsl{{\mathfrak{sl}}}

\def\fZ{{\mathfrak{Z}}}

\def\BA{{\mathbb{A}}}
\def\BC{{\mathbb{C}}}
\def\BK{{\mathbb{K}}}

\def\BN{{\mathbb{N}}}

\def\BQ{{\mathbb{Q}}}
\def\BZ{{\mathbb{Z}}}

\def\woo{\widehat{\otimes}}

\def\CL{{\mathcal{L}}}

\def\CS{{\mathcal{S}}}

\def\CV{{\mathcal{V}}}
\def\CW{{\mathcal{W}}}

\def\coop{{\textrm{coop}}}

\def\Hom{\textrm{Hom}}

\def\e{\varepsilon}

\def\vs{\varsigma}

\def\pt{\textrm{pt}}
\def\and{\textrm{ }\&\textrm{ }}

\def\sym{\textrm{sym}}
\def\Sym{\textrm{Sym}}
\def\esym{\emph{sym}}
\def\eSym{\emph{Sym}}

\def\oCS{\mathring{\CS}}

\def\barx{{\bar{x}}}
\def\bary{{\bar{y}}}
\def\barz{{\bar{z}}}

\def\nn{{{\BN}}^I}
\def\zz{{{\BZ}}^I}

\def\fg{\mathfrak{g}}

\def\uu{U_q(\fg)}

\def\uupm{U_q^\pm(\fg)}
\def\uug{U_q^{\geq}(\fg)}
\def\uul{U_q^{\leq}(\fg)}

\def\tuup{\widetilde{U}_q^+(\fg)}
\def\tuum{\widetilde{U}_q^-(\fg)}
\def\tuupm{\widetilde{U}_q^\pm(\fg)}
\def\tuug{\widetilde{U}_q^{\geq}(\fg)}
\def\tuul{\widetilde{U}_q^{\leq}(\fg)}

\def\UU{U_q(L\fg)}
\def\UUp{U_q^+(L\fg)}
\def\UUm{U_q^-(L\fg)}
\def\UUpm{U_q^\pm(L\fg)}
\def\UUmp{U_q^\mp(L\fg)}
\def\UUg{U_q^{\geq}(L\fg)}
\def\UUl{U_q^{\leq}(L\fg)}

\def\tUUp{\widetilde{U}_q^+(L\fg)}
\def\tUUm{\widetilde{U}_q^-(L\fg)}
\def\tUUpm{\widetilde{U}_q^\pm(L\fg)}
\def\tUUmp{\widetilde{U}_q^\mp(L\fg)}
\def\tUUg{\widetilde{U}_q^{\geq}(L\fg)}
\def\tUUl{\widetilde{U}_q^{\leq}(L\fg)}

\def\tUUpT{\widetilde{U}_q^+(L\fg)^T}

\def\bs{{\boldsymbol{\vs}}}

\def\ow{{\overline{w}}}

\def\bn{\boldsymbol{n}}

\def\tU{\widetilde{U}}

\def\op{\text{op}}

\def\oij{\overrightarrow{ij}}

\def\loc{\text{loc}}
\def\eloc{\emph{loc}}

\def\tUpsilon{\widetilde{\Upsilon}}

\begin{document}

\title[Quantum loop groups for symmetric Cartan matrices]{\Large{\textbf{Quantum loop groups for symmetric Cartan matrices}}}

\author[Andrei Negu\cb t]{Andrei Negu\cb t}
\address{MIT, Department of Mathematics, Cambridge, MA, USA}
\address{Simion Stoilow Institute of Mathematics, Bucharest, Romania}
\email{andrei.negut@gmail.com}

\maketitle

\begin{abstract} We introduce a quantum loop group associated to a general symmetric Cartan matrix, by imposing just enough relations between the usual generators $\{e_{i,k}, f_{i,k}\}_{i \in I, k \in \BZ}$ in order for the natural Hopf pairing between the positive and negative halves of the quantum loop group to be perfect. As an application, we describe the localized $K$-theoretic Hall algebra of any quiver without loops, endowed with a particularly important $\mathbb{C}^*$ action.

\end{abstract}

$$$$

\epigraph{\emph{To Ivo and Elio, my favorite perfect pairing}}

\section{Introduction}

\medskip

\subsection{} 
\label{sub:quantum group intro} 

Let us fix any symmetric Cartan matrix, i.e. $C = \{d_{ij}\}_{i,j \in I}$ with
\begin{equation}
\label{eqn:cartan}
\begin{cases} d_{ii} = 2& \text{if } i = j \\ 
d_{ij} = d_{ji} \in \{0,-1,-2,\dots\} &\text{if }i \neq j \end{cases}
\end{equation}
for some finite set $I$. Let $\fg$ denote the Kac-Moody Lie algebra associated to the Cartan matrix $C$. A well-known and very important Hopf algebra deformation of $U(\fg)$ is the Drinfeld-Jimbo quantum group, which may be defined by the following procedure. We will work over a ground field $\BK$ of characteristic 0 (usually taken to be $\BC(q)$), endowed with an element $q \in \BK^\times$ which is not a root of unity. 

\medskip

\begin{enumerate}[leftmargin=*]

\item[\emph{(1)}] Start from the algebra
$$
\tuup = \BK  \langle e_i  \rangle_{i \in I} 
$$

\medskip

\item[\emph{(2)}] To make $\tuup$ into a bialgebra, we first enlarge it
$$
\tuug = \tuup [h_i^{\pm 1}]_{i \in I} \Big/ \Big( h_i e_j = q^{d_{ij}} e_j h_i \Big)
$$
and define the coproduct
\begin{align*}
&\Delta(e_i) = h_i \otimes e_i + e_i \otimes 1 \\
&\Delta(h_i) = h_i \otimes h_i
\end{align*}

\medskip

\item[\emph{(3)}] Define $\tuul = \tuug^{\coop}$, with generators denoted by $f_i, h_i'$ instead of $e_i,h_i$. Then there is a bialgebra pairing (see \eqref{eqn:bialg 1}--\eqref{eqn:bialg 2})
\begin{equation}
\label{eqn:pairing intro 1}
\tuug \otimes \tuul \xrightarrow{\langle \cdot , \cdot \rangle} \BK
\end{equation}
completely determined by
$$
\Big \langle e_i,f_i \Big \rangle = \frac 1{q^{-1}-q}, \qquad \Big \langle h_i, h'_j \Big \rangle = q^{d_{ij}}
$$
and all other pairings between generators being 0. 

\medskip

\item[\emph{(4)}] Consider the radical of the pairing \eqref{eqn:pairing intro 1}, namely
$$
I^+_{\circ} \subset \tuup, \qquad x \in I^+_{\circ} \Leftrightarrow \Big \langle x, \tuum \Big \rangle = 0
$$
and define $I^-_{\circ} \subset \tuum$ analogously. Since $I^\pm_{\circ}$ are ideals, the quotients
$$
\uupm = \tuupm \Big/ I^\pm_{\circ}
$$
are algebras.

\medskip

\item[\emph{(5)}] We may define the bialgebras $\uug$ and $\uul$ by removing all the tildes in items \emph{(2)}--\emph{(3)}. Then the pairing \eqref{eqn:pairing intro 1} descends to a bialgebra pairing
\begin{equation}
\label{eqn:pairing intro 2}
\uug \otimes \uul \xrightarrow{\langle \cdot , \cdot \rangle} \BK
\end{equation}

\medskip

\item[\emph{(6)}] Define the quantum group as the vector space
\begin{equation}
\label{eqn:quantum group intro}
\uu = \uug \otimes \uul \Big/ \Big( h_i \otimes h'_i = 1 \otimes 1 \Big)
\end{equation}
made into an algebra by imposing the Drinfeld double relation \eqref{eqn:dd} between the subalgebras $\uug \otimes 1$ and $1 \otimes \uul$ of $\uu$. Note that the Drinfeld double relation only takes as input the bialgebra structures and the pairing defined in item \emph{(5)} \footnote{All the bialgebras studied in the present paper are also Hopf algebras, but we will not need the antipode map (and thus will not explicitly describe it, although it is straightforward to do so).}, and it gives rise to the well-known relation
$$
[e_i,f_j] = \delta_{ij} \cdot \frac {h_i-h_i'}{q-q^{-1}}
$$
in $\uu$, for all $i,j \in I$ (we identify $h_i$ with $h_i \otimes 1$ and $h_i'$ with $1 \otimes h_i'$ in \eqref{eqn:quantum group intro}).

\medskip

\end{enumerate}

\noindent Although somewhat dry, the procedure above yields beautiful formulas, when one asks to describe the radicals $I^\pm_{\circ}$ explicitly. For example, it was shown in \cite[Theorem 33.1.3]{Lu} that $I^+_{\circ}$ is generated as a two-sided ideal by the LHS of the relation 
\begin{equation}
\label{eqn:serre}
\sum_{k = 0}^{1-d_{ij}} (-1)^k {1-d_{ij} \choose k}_{q} e_i^k e_j e_i^{1-d_{ij}-k}  = 0
\end{equation}
\footnote{Above, we use the notation ${n \choose k}_q = \frac {[n]!_q}{[k]!_q [n-k]!_q}$ where $[n]!_q = [1]_q \dots [n]_q$ and $[n]_q = \frac {q^n - q^{-n}}{q-q^{-1}}$.} as $(i,j)$ runs over all pairs of distinct elements of $I$. The analogous result holds in $I^-_{\circ}$ if one replaces all the $e$'s by $f$'s. Relation \eqref{eqn:serre} is not too surprising, in light of the fact that as $q \rightarrow 1$, it converges to the usual Serre relation that holds in $\fg$
\begin{equation}
\label{eqn:serre limit}
[e_i,[e_i,[\dots,[e_i,e_j]\dots]]] = 0
\end{equation}
where the number of brackets is $1-d_{ij}$. The fact that relations \eqref{eqn:serre limit} generate the ($q\rightarrow 1$ limit of the) ideal $I_\circ^+$ was proved by Serre for Lie algebras $\fg$ of finite type, and by Gabber-Kac for general Kac-Moody Lie algebras $\fg$.

\medskip

\subsection{} 
\label{sub:quantum loop group intro}

The main purpose of the present paper is to carry out the procedure above for quantum loop groups, i.e. to appropriately define a deformation of the universal enveloping algebra of $L\fg = \fg[t^{\pm 1}]$ with Lie bracket given by
\begin{equation}
\label{eqn:lie loop}
[xt^k, yt^l] = [x,y] t^{k+l}
\end{equation}
for all $x,y \in \fg$ and all $k,l \in \BZ$. Quantum loop groups are important players in the theory of quantum integrable systems and quantum field theory. Moreover, together with their Yangian degenerations and their elliptic versions, quantum loop groups have played big roles in geometric representation theory in recent decades.

\medskip

\noindent The natural analogue of the six-step procedure in Subsection \ref{sub:quantum group intro} is the following. The construction of quantum loop groups that we are about to recall is due to \cite{Dr}, who used it \footnote{Here, ``it" more appropriately refers to a central extension of $\UU$.} to produce an alternate construction of quantum affine groups (this corresponds to the case when $C$ is a Cartan matrix of finite type).

\medskip

\begin{enumerate}[leftmargin=*]

\item[\emph{(1)}] Start from the algebra
\begin{equation}
\label{eqn:def quad intro}
\tUUp = \BK \langle e_{i,k} \rangle_{i \in I, k \in \BZ} \Big / \text{relation \eqref{eqn:quad intro}}
\end{equation}
where we consider the formal series $e_i(z) = \sum_{k \in \BZ} \frac {e_{i,k}}{z^k}$ and impose the relations
\begin{equation}
\label{eqn:quad intro}
e_i(z)e_j(w) (z-wq^{d_{ij}}) = e_j(w)e_i(z) (zq^{d_{ij}}-w)
\end{equation}
for all $i,j \in I$. To motivate relation \eqref{eqn:quad intro}, note that as $q \rightarrow 1$ it converges to
$$
[e_i(z), e_j(w)](z-w) = 0
$$
Unpacking this relation says that the commutator $[e_{i,k}, e_{j,l}]$ only depends on $k+l$, which is to be expected from the Lie bracket equality \eqref{eqn:lie loop}. 

\medskip

\item[\emph{(2)}] Enlarge $\tUUp$ as follows
$$
\tUUg = \frac {\tUUp [h_{i,0}^{\pm 1}, h_{i,1}, h_{i,2},\dots]_{i \in I} }{\left( h_i(z) e_j(w) = e_j(w) h_i(z) \frac {zq^{d_{ij}} - w}{z - w q^{d_{ij}}} \right)}
$$
where $h_i(z) = \sum_{k=0}^{\infty} \frac {h_{i,k}}{z^k}$, and define the topological coproduct
\begin{align*}
&\Delta(e_i(z)) = h_i(z) \otimes e_i(z) + e_i(z) \otimes 1 \\
&\Delta(h_i(z)) = h_i(z) \otimes h_i(z)
\end{align*}

\medskip

\item[\emph{(3)}] Define $\tUUl = \tUUg^{\coop}$, with generators denoted by $f_{i,-k}, h_{i,-k}'$ instead of $e_{i,k}, h_{i,k}$. Then there is a topological bialgebra pairing
\begin{equation}
\label{eqn:pairing intro loop 1}
\tUUg \otimes \tUUl \xrightarrow{\langle \cdot, \cdot \rangle} \BK
\end{equation}
completely determined by
\begin{equation}
\label{eqn:pairing from the intro}
\Big \langle e_{i,k}, f_{i,-k} \Big \rangle = \frac 1{q^{-1}- q}, \qquad \Big \langle h_i(z), h'_j(w)\Big\rangle = \frac {zq^{d_{ij}}-w}{z-wq^{d_{ij}}}
\end{equation}
and all other pairings between generators being 0.

\medskip

\item[\emph{(4)}] Consider the radical
$$
I^\pm \subset \widetilde{U}_q^\pm(L\fg)
$$
of the pairing \eqref{eqn:pairing intro loop 1}, and define the algebras $U_q^\pm(L\fg) = \widetilde{U}_q^\pm(L\fg) \Big/ I^\pm$. 

\medskip

\item[\emph{(5)}] We may define the bialgebras $\UUg$ and $\UUl$ by removing all the tildes in items \emph{(2)}--\emph{(3)}. Then the pairing \eqref{eqn:pairing intro loop 1} descends to a bialgebra pairing
\begin{equation}
\label{eqn:pairing intro loop 2}
\UUg \otimes \UUl \xrightarrow{\langle \cdot , \cdot \rangle} \BK
\end{equation}

\medskip

\item[\emph{(6)}] Define the quantum group as
\begin{equation}
\label{eqn:quantum loop group intro}
\UU = \UUg \otimes \UUl \Big/ \Big( h_{i,0} \otimes h'_{i,0} = 1 \otimes 1 \Big)
\end{equation}
with the multiplication governed by \eqref{eqn:dd}. It implies the well-known relation
\begin{equation}
\label{eqn:rel 6 full}
[e_i(z), f_j(w)] =\delta_{ij} \delta \left(\frac zw\right) \cdot \frac { h_i(z) - h_i'(w)}{q-q^{-1}}
\end{equation}
where $\delta(x) = \sum_{k \in \BZ} x^k$ is a formal series and $\delta_{ij}$ is the Kronecker delta function.

\medskip

\end{enumerate}

\subsection{} 
\label{sub:zig-zag intro}

In the present paper, we will complete the definition of the quantum loop group $\UU$ by explicitly describing the radical $I^+$ (the negative part $I^-$ will be analogous, with $e_{i,k}$'s replaced by $f_{i,-k}$'s); specifically, we develop generators for the ideal $I^+$ that provide quantum loop group versions of relations \eqref{eqn:serre}. 

\medskip

\noindent The first combinatorial tool that we need for this is the notion of \textbf{distinguished zig-zags} $Z$, developed in Subsection \ref{sub:zig-zag}. In brief, for any $i \neq j \in I$ and any arithmetic progressions $s,\dots,t$ and $s',\dots, t'$ with step 2, such that $s+t=s'+t'$, we draw the following oriented graph

\begin{figure}[h]
\centering
\includegraphics[scale=0.25]{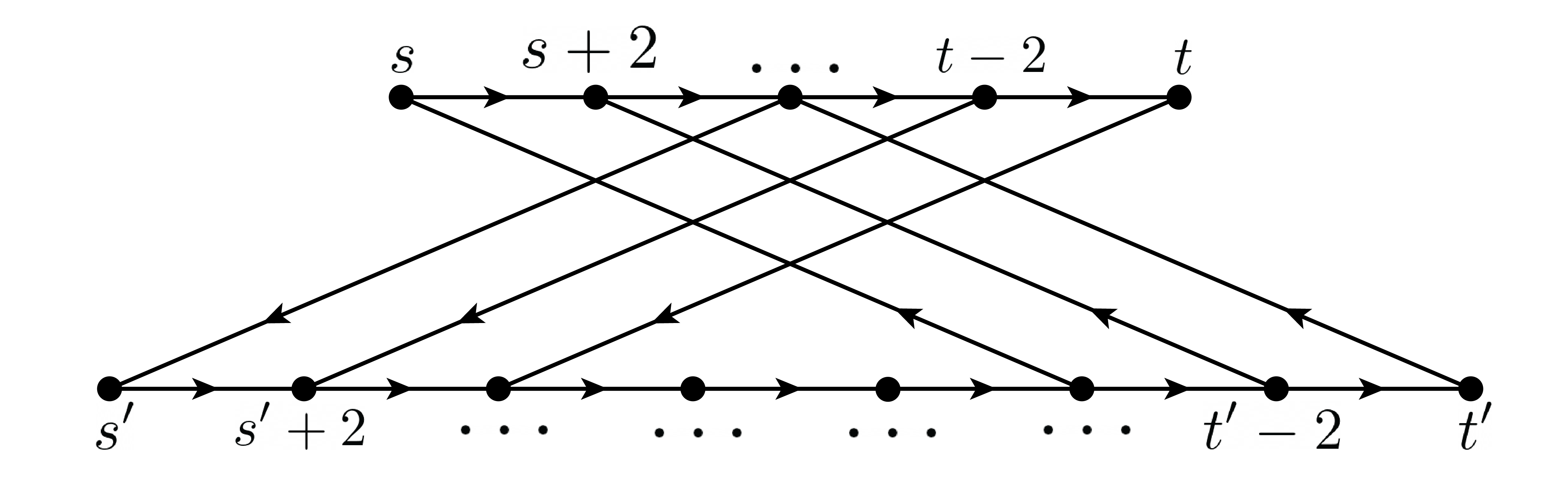} 
\caption{A distinguished zig-zag $Z$}
\end{figure}

\noindent The vertices of the graph are placed on two horizontal lines at coordinates prescribed by the chosen arithmetic progressions, and the diagonal edges correspond to all $(a,b) \in \{s,\dots,t\} \times \{s', \dots, t'\}$ such that $a = b\pm d_{ij}$. We may think of the top (respectively bottom) vertices of the graph $Z$ as being colored by $i$ (respectively $j$). 

\medskip

\noindent The second combinatorial tool that we need is that of a \textbf{refined selection}, developed in Subsection \ref{sub:refined selection}, which is a particular type of multiset $\mathscr{S}$ of edges in the oriented graph of Figure 1, whose removal leaves the graph without any oriented cycles. For any $i \neq j$ and any distinguished zig-zag $Z$, consider the formal series

\begin{equation}
\label{eqn:rho intro}
\rho_Z(x_s,\dots,x_t,y_{s'},\dots,y_{t'}) = \sum_{\mathscr{S} \text{ refined selection}} (-1)^{\sigma(\mathscr{S})} \cdot 
\end{equation}
$$
\frac {\prod_{\text{vertices } c < c'} \left[ (-1)^{\delta_{\iota(c)i} \delta_{\iota(c')j}+ \delta_{\iota(c)\iota(c')} \nu_{c'-c}} (z_c - z_{c'} q^{-d_{\iota(c)\iota(c')}}) \right]}{\prod_{\text{edges } c' \rightarrow c} \left(\barz_{c'} - \barz_c \right)^{1-\mu_{\mathscr{S}}(c'\rightarrow c)}}  \mathop{\prod_{\text{vertices } c \text{ in}}}_{\text{descending order}} e_{\iota(c)}(z_c) 
$$
in $\tUUp [\![ x^{\pm 1}_s,\dots,x^{\pm 1}_t, y^{\pm 1}_{s'}, \dots, y^{\pm 1}_{t'} ]\!]$. Above,

\medskip

\begin{itemize}[leftmargin=*]

\item $c$ and $c'$ go over the vertex set of the graph in Figure 1, i.e. $\{s,\dots,t\} \sqcup \{s',\dots,t'\}$;

\medskip

\item $x_s,\dots,x_t,y_{s'},\dots,y_{t'}$ are formal variables, and all other notations in the formula above are defined explicitly in \eqref{eqn:notation 1}--\eqref{eqn:notation 4}; the various $z_c, \barz_c$ are equal to various $x_i,y_j$ times powers of $q$;

\medskip

\item the sign $(-1)^{\sigma(\mathscr{S})}$ is defined in \eqref{eqn:sign}, and 
$$
\mu_{\mathscr{S}}(c' \rightarrow c) \in \{0,1,\dots\}
$$
denotes the multiplicity of the edge $c' \rightarrow c$ in the multiset $\mathscr{S}$;

\medskip

\item $c<c'$ denotes any total order of the vertex set, such that there is no edge in the complement of $\mathscr{S}$ that points from a smaller vertex to a larger vertex \footnote{This choice is made possible by the fact that removing $\mathscr{S}$ leaves no oriented cycles in the graph; the element on the second line of \eqref{eqn:rho intro} does not depend on the choice of total order associated to a refined selection $\mathscr{S}$, as shown in Proposition \ref{prop:doesn't matter}.};

\item $\delta_{kl}$ is the Kronecker delta symbol for any $k,l \in I$, and $\nu_{c'-c}$ is 1 (respectively 0) if $c'$ is greater (respectively smaller) than $c$ with respect to the usual ordering of the integers, which may differ from the total order $<$ in the previous bullet.

\end{itemize}

\medskip

\noindent Note that the second line of \eqref{eqn:rho intro} is a Laurent polynomial in the variables $x_s,\dots,x_t$, $y_{s'},\dots,y_{t'}$ times various series $e_{i\text{ or }j}(x_a \text{ or }y_b)$, despite the apparent denominators.

\medskip

\begin{theorem}
\label{thm:main} 

The ideal $I^+$ is generated by the coefficients of the series $\rho_Z$ of \eqref{eqn:rho intro}, as $Z$ goes over all distinguished zig-zags, for any $i \neq j$ in $I$. Thus, the positive half of the quantum loop group may be defined as
\begin{equation}
\label{eqn:main}
\UUp = \BK \langle e_{i,k} \rangle_{i \in I, k \in \BZ} \Big /  \text{relations \eqref{eqn:quad intro} and } \{\rho_Z = 0\}_{Z \text{ distinguished zig-zag}}
\end{equation}
The negative half $\UUm$ can be defined analogously (with $f_{i,-k}$'s instead of $e_{i,k}$'s). 

\end{theorem}

\medskip

\noindent Since $\rho_Z$ are formal series in variables $x_s,\dots,x_t, y_{s'},\dots, y_{t'}$, one might be concerned that \eqref{eqn:main} requires setting all the coefficients of said formal series to 0, thus yielding ``too many" relations. In fact, \eqref{eqn:main} still holds if we only ask for the vanishing of a single generic coefficient of $\rho_Z$ of any given homogeneous degree in the variables $x_s,\dots,x_t, y_{s'},\dots, y_{t'}$, for any distinguished zig-zag $Z$ as in Figure 1. The exact meaning of the word ``generic" in the previous sentence is given in Remark \ref{rem:replace}. 

\medskip

\noindent Distinguished zig-zags with a single northwest-pointing (equivalently, southwest-pointing) diagonal arrow are called \textbf{minimal}, and for fixed $i \neq j$, they are in one-to-one correspondence with pairs of non-negative integers $k,l$ such that
\begin{equation}
\label{eqn:integers intro}
k+l = -d_{ij} 
\end{equation}
where $k+1$ (respectively $l+1$) is the number of vertices on the top (respectively bottom) row of the zig-zag. For example, when $k=2$ and $l=5$, the minimal zig-zag takes the form in Figure 2 (note that $t = s+2k$ and $t'=s'+2l$ therein). 

\medskip

\begin{figure}[h]
\centering
\includegraphics[scale=0.25]{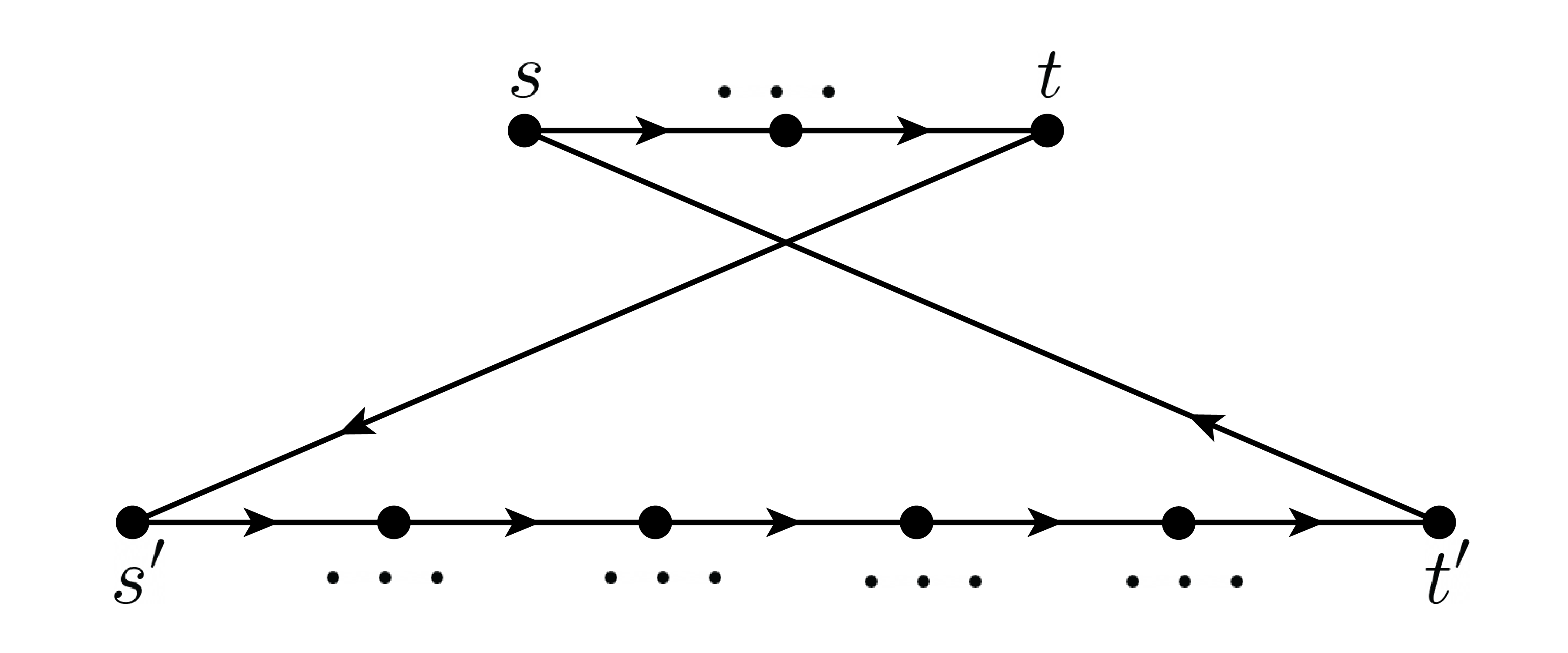} 
\caption{A minimal zig-zag $Z$}
\end{figure}

\noindent When $Z$ is the minimal zig-zag associated to $k,l$ as in \eqref{eqn:integers intro}, then the refined selections $\mathscr{S}$ that appear in \eqref{eqn:rho intro} are simply all one-edge subsets of the oriented graph in Figure 2. Thus, there are as many refined selections as there are total orders of the vertices $\{s,\dots,t\} \sqcup \{s',\dots,t'\}$ that respect the orientation of the zig-zag with one edge removed. We list all these possible orders below, as they dictate the terms that appear on the second line of \eqref{eqn:rho intro} (meanwhile, the signs $(-1)^{\sigma(\mathscr{S})}$ that appear on the first line \eqref{eqn:rho intro} are all 1 in the case at hand)
\begin{multline*}
x_{s+2a} > x_{s+2(a+1)} > \dots > x_{t-2} > x_{t} > \\ > y_{s'} > y_{s'+2} > \dots > y_{t'-2} > y_{t'} > \\ > x_s > x_{s+2} >  \dots > x_{s+2(a-2)} > x_{s+2(a-1)}
\end{multline*}
for all $a \in \{1,\dots,k\}$,
\begin{multline*}
y_{s'+2b} > y_{s'+2(b+1)} > \dots > y_{t'-2} > y_{t'} > \\ > x_{s} > x_{s+2} > \dots > x_{t-2} > x_{t} > \\ > y_{s'} > y_{s'+2} >  \dots > y_{s'+2(b-2)} > y_{s'+2(b-1)}
\end{multline*}
for all $b \in \{1,\dots,l\}$, as well as
$$
x_s > x_{s+2} >  \dots > x_{t-2} > x_t > y_{s'} > y_{s'+2} > \dots > y_{t'-2} > y_{t'}
$$
and
$$
y_{s'} > y_{s'+2} > \dots > y_{t'-2} > y_{t'} > x_s > x_{s+2} >  \dots > x_{t-2} > x_t
$$
We believe that the relations $\rho_Z = 0$ were not known before for general $Z$. However, in the particular case of the minimal zig-zag corresponding to $k=-d_{ij}$, $l=0$,

\begin{figure}[h]
\centering
\includegraphics[scale=0.25]{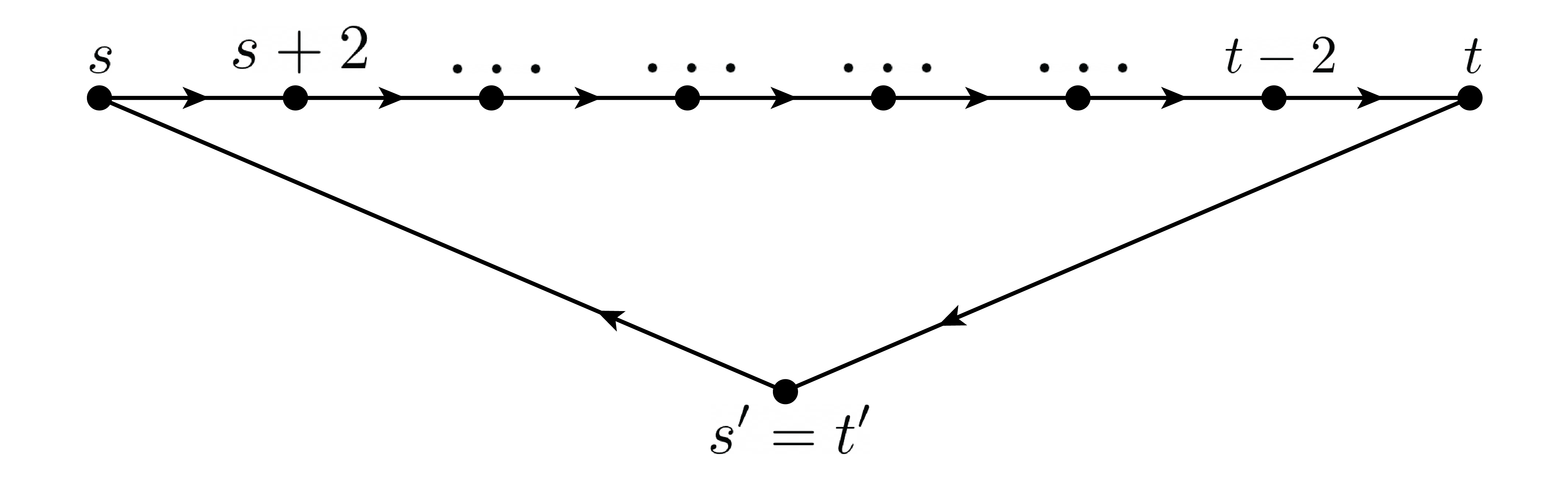} 
\caption{A particular minimal zig-zag $Z$}
\end{figure}

\noindent we will show in Proposition \ref{prop:equivalent} that the relations $\rho_Z = 0$ are equivalent to the following well-known loop versions of relations \eqref{eqn:serre} 
\begin{equation}
\label{eqn:loop serre}
\Sym \left[ \sum_{k = 0}^{1-d_{ij}} (-1)^k {1-d_{ij} \choose k}_{q} e_i(z_1) \dots e_i(z_k) e_j(w) e_i(z_{k+1}) \dots e_i(z_{1-d_{ij}}) \right] = 0
\end{equation}
where $\Sym$ denotes symmetrization in the variables $z_1,\dots,z_{1-d_{ij}}$. This phenomenon underscores the fact that our formulas for $\rho_Z$ are not unique; one could add arbitrary multiples of relation \eqref{eqn:quad intro} to $\rho_Z$ and obtain equivalent formulas. This also reveals the novelty behind our more general relations \eqref{eqn:rho intro}: while \eqref{eqn:loop serre} are all the necessary higher order relations when $d_{ij} \in \{0,-1\}$, they are already not sufficient when $d_{ij}=-2$ since they do not account for the minimal zig-zag with $t-s = t'-s'=2$.

\medskip

\subsection{}
\label{sub:shuffle intro}

Our main technical tool for proving Theorem \ref{thm:main} is the following trigonometric degeneration of the Feigin-Odesskii (\cite{FO}) shuffle algebra
\begin{equation}
\label{eqn:big shuf intro}
\CV^+  = \CV^- = \bigoplus_{\bn \in \nn} \frac {\BK [z_{i1}^{\pm 1}, \dots, z_{in_i}^{\pm 1}]_{i \in I}^{\sym}}{\prod_{i < j \in I} \prod_{a \leq n_i, b \leq n_j} (z_{ia} - z_{jb})}
\end{equation}
with the multiplication in $\CV^\pm$ defined in Subsection \ref{sub:def shuf} (the superscript ``sym" in the numerator of \eqref{eqn:big shuf intro} denotes symmetric Laurent polynomials in $z_{i1},z_{i2},\dots$ for all $i \in I$ separately; meanwhile, the symbol $<$ in the denominator of \eqref{eqn:big shuf intro} refers to an arbitrary fixed total order on $I$). The following facts were noted in \cite{E}:

\medskip

\begin{itemize}[leftmargin=*]

\item There are algebra homomorphisms $\widetilde{\Upsilon}^\pm : \tUUpm \rightarrow \CV^\pm$ given by $e_{i,k}, f_{i,k} \mapsto z_{i1}^k$

\medskip

\item There are bialgebra pairings \footnote{Strictly speaking, to claim that the pairings below have the bialgebra property requires one to replace $\tUUpm$, $\CV^\pm$ by $\tUUg$, $\tUUl$, $\CV^{\geq}$, $\CV^{\leq}$, where the latter two bialgebras are defined by adding generators $h_{i,k}, h'_{i,k}$ to $\CV^\pm$, akin to items \emph{(2)}--\emph{(3)} of Subsection \ref{sub:quantum loop group intro}; see Subsection \ref{sub:ext}.}
\begin{align}
&\tUUp \otimes \CV^{-} \xrightarrow{\langle \cdot, \cdot \rangle_{\tU\CV}} \BK \label{eqn:pairing shuffle intro 1}
\\
&\CV^{+} \otimes \tUUm \xrightarrow{\langle \cdot, \cdot \rangle_{\CV\tU}} \BK \label{eqn:pairing shuffle intro 2}
\end{align}
which will be shown in Proposition \ref{prop:non-deg} to be non-degenerate in both arguments. 

\end{itemize}

\medskip 

\noindent The pairings \eqref{eqn:pairing intro loop 1} and \eqref{eqn:pairing shuffle intro 1}--\eqref{eqn:pairing shuffle intro 2} are compatible, in the sense that
\begin{align*}
&\Big \langle x^+, y^- \Big \rangle = (q^{-1}-q)^{-|\deg x^+|} \Big \langle x^+, \widetilde{\Upsilon}^-(y^-) \Big \rangle_{\tU\CV} \\
&\Big \langle y^+, x^- \Big \rangle = (q^{-1}-q)^{|\deg x^-|} \Big \langle \widetilde{\Upsilon}^+(y^+),x^- \Big \rangle_{\CV\tU}
\end{align*}
for all $x^\pm, y^\pm \in \tUUpm$ (see \eqref{eqn:grading u}--\eqref{eqn:total degree} for the notation $|\deg x^\pm|$).

\medskip

\subsection{}
\label{sub:wheel intro}

In Definition \ref{def:wheel}, we will construct subalgebras
\begin{equation}
\label{eqn:small shuffle intro}
\CS^\pm \subset \CV^\pm
\end{equation}
We refer the reader to Subsection \ref{sub:wheel} for the full details and notation, but in brief, $\CS^\pm$ consists of those rational functions
$$
\frac {r(z_{i1},\dots,z_{in_i})_{i \in I}}{\prod_{i < j \in I} \prod_{a \leq n_i, b \leq n_j} (z_{ia} - z_{jb})} \in \CV^\pm
$$
such that for any distinguished zig-zag as in Figure 1, the specialization of the Laurent polynomial $r$ at the variables
\begin{align}
&z_{i1} = xq^s, \ \quad z_{i2} = xq^{s+2}, \ \ \quad z_{i,\frac {t-s}2} = xq^{t-2}, \ \ \quad z_{i, \frac {t-s}2+1} = xq^t \label{eqn:spec 1 intro} \\
&z_{j1} = yq^{s'}, \quad z_{j2} = yq^{s'+2}, \quad z_{j,\frac {t'-s'}2} = yq^{t'-2}, \quad z_{j, \frac {t'-s'}2+1} = yq^{t'} \label{eqn:spec 2 intro}
\end{align}
(with $\{s,\dots,t\} \sqcup \{s',\dots,t'\}$ as in Subsection \ref{sub:zig-zag intro}) has the property that
\begin{equation}
\label{eqn:wheel intro}
r\Big|_{\text{\eqref{eqn:spec 1 intro}--\eqref{eqn:spec 2 intro}}} \text{ is divisible by }(x-y)^{\text{number of southwest pointing arrows in Figure 1}}
\end{equation}
When the zig-zag is the specific one displayed in Figure 3, the aforementioned divisibility condition was known to \cite{E, FO}, but we believe the divisibility conditions provided by general distinguished zig-zags are new.

\medskip

\begin{theorem}
\label{thm:shuf}

(Proposition \ref{prop:coincide}) We have $\CS^\pm = \emph{Im }\widetilde{\Upsilon}^\pm$.

\end{theorem}

\medskip

\noindent The connection between the subalgebras $\CS^\pm$ and Theorem \ref{thm:main} is the following. 

\medskip

\begin{theorem} 
\label{thm:connect}

The subalgebras $\CS^\pm \subset \CV^\pm$ and the quotients $\UUmp \twoheadleftarrow \tUUmp$ are mutually dual, i.e. the pairings \eqref{eqn:pairing shuffle intro 1}--\eqref{eqn:pairing shuffle intro 2} descend to pairings given by the diagonal arrows in the diagrams below 
\begin{equation}
\label{eqn:desired diagram intro}
\begin{tikzcd}[column sep=0.05em]
\tUUp \arrow[d,twoheadrightarrow] & \otimes & \CV^- \arrow{rrrrrrrrrrrrrrrrrrrrrrrrrrrrrrrrrrrrrrrrrrrrrrrr}{\langle \cdot , \cdot \rangle_{\tU\CV}} & & & & & & & & & & & & & & & & & & & & & & & & & & & & & & & & & & & & & & & & & & & & & & & & \BK\\
\UUp  & \otimes & \CS^- \arrow[u, hook] \arrow{rrrrrrrrrrrrrrrrrrrrrrrrrrrrrrrrrrrrrrrrrrrrrrrru}[swap]{\langle \cdot , \cdot \rangle_{U\CS}}
\end{tikzcd} \qquad \qquad \begin{tikzcd}[column sep=0.05em]
\CV^+  & \otimes & \tUUm \arrow{rrrrrrrrrrrrrrrrrrrrrrrrrrrrrrrrrrrrrrrrrrrrrrrr}{\langle \cdot, \cdot \rangle_{\CV\tU}} \arrow[d,twoheadrightarrow] & & & & & & & & & & & & & & & & & & & & & & & & & & & & & & & & & & & & & & & & & & & & & & & & \BK\\
\CS^+ \arrow[u, hook] & \otimes & \UUm  \arrow{rrrrrrrrrrrrrrrrrrrrrrrrrrrrrrrrrrrrrrrrrrrrrrrru}[swap]{\langle \cdot, \cdot \rangle_{\CS U}}
\end{tikzcd}
\end{equation}
The pairings $\langle \cdot, \cdot \rangle_{U \CS}$ and $\langle \cdot, \cdot \rangle_{\CS U}$ are non-degenerate in both arguments. 

\end{theorem}

\medskip

\noindent Another way to rephrase Theorem \ref{thm:connect} is that the generators of the ideals $I^\pm$ can be constructed as linear functionals on $\CV^\pm$ which realize the conditions \eqref{eqn:wheel intro} that define the subalgebra $\CS^\pm$.

\medskip

\begin{theorem} 
\label{thm:iso}

The induced homomorphisms
\begin{equation}
\label{eqn:iso intro}
\Upsilon^\pm : \UUpm \rightarrow \CS^\pm
\end{equation}
are isomorphisms. The pairings $\langle \cdot, \cdot \rangle_{U \CS}$ and $\langle \cdot, \cdot \rangle_{\CS U}$  both coincide with a pairing
\begin{equation}
\label{eqn:pairing intro final}
\CS^{+} \otimes \CS^{-} \xrightarrow{\langle \cdot, \cdot \rangle_{\CS\CS}} \BK
\end{equation}
which is non-degenerate in both arguments.

\end{theorem}

\medskip

\subsection{} The present paper follows the philosophy of \cite{Quiver 1, Quiver 3}, with the following important observation: various specializations of the parameters that appear in shuffle algebras lead to wildly different generators-and-relations presentations of the associated quantum loop groups. In more detail, consider a quiver with vertex set $I$ and $-d_{ij}$ total arrows between any distinct vertices $i \neq j$. The shuffle algebra considered in relation to this quiver in \cite{Quiver 1, Quiver 3} involved ``generic" parameters $t_e$ associated to the aforementioned arrows; this led to the corresponding shuffle algebra being determined by 3-variable wheel conditions (\cite[equation (2.14)]{Quiver 1}), and the corresponding quantum loop group being determined by cubic relations (\cite[equation (1.6)]{Quiver 3}). The algebra studied in our Definition \ref{def:geom} corresponds to the specialization 
\begin{equation}
\label{eqn:weights intro}
q^{d_{ij}+2}, \dots, q^{-d_{ij}-2}, q^{-d_{ij}}
\end{equation}
of the parameters associated to the $-d_{ij}$ arrows between $i$ and $j$. This choice yields the more complicated conditions \eqref{eqn:wheel intro} and the more complicated relations $\rho_Z = 0$. We will study the situation of more general specializations of parameters in \cite{Quiver 2}.

\medskip

\subsection{} 
\label{sub:k-ha}

An important feature of quantum loop groups in geometric representation theory is their relation to $K$-theoretic Hall algebras of quivers (we follow the presentation of \cite{VV}, wherein the interested reader may also find an overview of the theory and its history). Specifically, to any quiver $Q$ without loops, we may associate the symmetric Cartan matrix $C$ given by
\begin{equation}
\label{eqn:number of edges}
-d_{ij} = \Big | \Big\{ \text{arrows between } i \text{ and }j \Big\} \Big |
\end{equation}
In Section \ref{sec:k-ha}, we will recall the construction of the $K$-theoretic Hall algebra \footnote{The torus $\BC^*$ in \eqref{eqn:k-ha intro} acts on the linear maps corresponding to the $-d_{ij}$ arrows from $i$ to $j$ in the doubled quiver via the characters \eqref{eqn:weights intro}; these are known as ``normal weights" in \cite{VV}.}
\begin{equation}
\label{eqn:k-ha intro}
K^{\text{nilp}}_{\BC^*, \loc} = \bigoplus_{\bn \in \nn} K_{\BC^*} (T^*\fZ_{\bn})_{\Lambda_{\bn}} \bigotimes_{\BZ[q^{\pm 1}]} \BQ(q)
\end{equation}
One \footnote{The word ``one" here conceals decades of foundational work in geometric representation theory; the constructions herein have appeared in many contexts in the work of many mathematicians, but the most relevant to our purposes are \cite{G,SV Hilb,VV}.} defines a convolution product on $K^{\text{nilp}}_{\BC^*, \loc}$, and a linear map
\begin{equation}
\label{eqn:k-ha hom}
K^{\text{nilp}}_{\BC^*, \loc} \longrightarrow \CV^{+,\text{geom}} := \bigoplus_{\bn \in \nn} \BQ(q)[z_{i1}^{\pm 1}, \dots, z_{in_i}^{\pm 1}]^{\sym}_{i \in I}
\end{equation}
which is an algebra morphism, once $\CV^{+,\text{geom}}$ is made into an algebra using the shuffle product in Subsection \ref{sub:zeta geom}. We have the following commutative diagram of algebra homomorphisms
$$
\begin{tikzcd}
\CV^+ \arrow[r,"\Omega"] & \CV^{+,\text{geom}} \\
\CS^+ \arrow[u, hook] \arrow[r,"\sim"] & \CS^{+,\text{geom}} \arrow[u, hook]
\end{tikzcd}
$$
with the notation as in \eqref{eqn:v geom}--\eqref{eqn:s geom}. Denote the image of the map \eqref{eqn:k-ha hom} by
$$
\overline{K}^{\text{nilp}}_{\BC^*, \loc} \subset \CV^{+,\text{geom}}
$$
In simply-laced types (i.e. $d_{ij} \in \{0,-1\}$ for all $i\neq j$), we have $\overline{K}^{\text{nilp}}_{\BC^*, \loc} = K^{\text{nilp}}_{\BC^*, \loc}$; this is not true in general due to the failure of the map \eqref{eqn:k-ha hom} to be injective.

\medskip

\begin{theorem}
\label{thm:k-ha}

For any quiver $Q$ without loops, we have 
\begin{equation}
\label{eqn:k-ha}
\overline{K}^{\emph{nilp}}_{\BC^*, \emph{loc}} =  \CS^{+,\emph{geom}} \stackrel{\Omega}{\cong} \CS^+ \stackrel{\eqref{eqn:iso intro}}{\cong} \UUp
\end{equation}

\end{theorem}

\medskip

\noindent The isomorphism $K^{\text{nilp}}_{\BC^*, \text{loc}} \cong \UUp
$ was previously known for $\fg$ of finite and affine type other than $A_1^{(1)}$ (see \cite[Theorem A]{VV}). 

\medskip

\subsection{} The plan of the paper is the following.

\medskip

\begin{itemize}[leftmargin=*]

\item In Section \ref{sec:shuffle}, we define the algebras $\oCS^\pm \subseteq \CS^\pm \subset \CV^\pm$. The fact that $\oCS^\pm = \CS^\pm$ is stated in Proposition \ref{prop:coincide}, and will be proved in Section \ref{sec:words}. 

\medskip

\item In Section \ref{sec:quantum}, we realize the quotient $\tUUpm \twoheadrightarrow \UUpm$ as dual to $\CS^\mp \subset \CV^\mp$, and prove Theorems \ref{thm:main}, \ref{thm:connect} and \ref{thm:iso} (modulo a technical result, Proposition \ref{prop:final non-deg}, that will be proved in Section \ref{sec:words}).

\medskip

\item In Section \ref{sec:words}, we use the combinatorics of words to prove two outstanding technical results, namely Propositions \ref{prop:coincide} and \ref{prop:final non-deg}.

\medskip

\item In Section \ref{sec:k-ha}, we discuss $K$-theoretic Hall algebras in relation to shuffle algebras, and prove Theorem \ref{thm:k-ha}.

\end{itemize}

\medskip

\noindent I would like to thank Alexander Tsymbaliuk for many years of stimulating and productive conversations on shuffle algebras and quantum groups. I would also like to thank Pavel Etingof, Andrei Okounkov, Olivier Schiffmann, Michela Varagnolo and Éric Vasserot for all the geometric representation theory they taught me. I gratefully acknowledge NSF grant DMS-$1845034$, as well as support from the Alfred P.\ Sloan Foundation and the MIT Research Support Committee.

\medskip

\subsection{}
\label{sec:notation}

Let us summarize the main notations used in the Introduction above, as they will be used repeatedly throughout the paper. The main algebras involved are
\begin{equation}
\label{eqn:main algebras}
\widetilde{\Upsilon}^\pm : \tUUpm \twoheadrightarrow \CS^\pm \subset \CV^\pm
\end{equation}
A priori, the image of $\widetilde{\Upsilon}^\pm$ will be denoted by $\oCS^\pm \subset \CV^\pm$; one of our most important technical results is the fact that $\oCS^\pm = \CS^\pm$, to be given in Proposition \ref{prop:coincide}. One of our main results is the fact that \eqref{eqn:main algebras} descends to an isomorphism
\begin{equation}
\label{eqn:main iso}
\Upsilon^\pm : \UUpm \xrightarrow{\sim} \CS^\pm
\end{equation}
where $\UUpm$ is the explicit quotient in \eqref{eqn:main}. Throughout the present paper, we will encounter pairings
\begin{align}
\tUUp \otimes \tUUm &\xrightarrow{\langle \cdot , \cdot \rangle_{\widetilde{U}\widetilde{U}}} \BK \label{eqn:first pair} \\
\tUUp \otimes \CV^- &\xrightarrow{\langle \cdot , \cdot \rangle_{\widetilde{U}\CV}} \BK \\
\CV^+ \otimes \tUUm &\xrightarrow{\langle \cdot , \cdot \rangle_{\CV\widetilde{U}}} \BK \\
\tUUp \otimes \oCS^- &\xrightarrow{\langle \cdot , \cdot \rangle_{\widetilde{U}\oCS}} \BK \\
\oCS^+ \otimes \tUUm &\xrightarrow{\langle \cdot , \cdot \rangle_{\oCS\widetilde{U}}} \BK \\
\tUUp \otimes \CS^- &\xrightarrow{\langle \cdot , \cdot \rangle_{\widetilde{U}\CS}} \BK \\
\CS^+ \otimes \tUUm &\xrightarrow{\langle \cdot , \cdot \rangle_{\CS\widetilde{U}}} \BK \\
\oCS^+ \otimes \CS^- &\xrightarrow{\langle \cdot , \cdot \rangle_{\oCS\CS}} \BK \\
\CS^+ \otimes \oCS^- &\xrightarrow{\langle \cdot , \cdot \rangle_{\CS\oCS}} \BK \\
\CS^+ \otimes \CS^- &\xrightarrow{\langle \cdot , \cdot \rangle_{\CS\CS}} \BK \label{eqn:last pair}
\end{align}
In all cases, the indices under $\langle \cdot, \cdot \rangle$ are meant to reflect the domain of the pairing in question (which will be important when we discuss issues like non-degeneracy). All the pairings above are compatible with each other under the natural inclusion maps, quotient maps, as well as the homomorphisms \eqref{eqn:main algebras}--\eqref{eqn:main iso}. 

\medskip

\noindent One can extend all the algebras above by adding generators $h_{i,k}$, $h_{i,k}'$ in a compatible way. The resulting objects are topological bialgebras, and they will be denoted by
\begin{equation}
\label{eqn:extended}
\tUUg, \tUUl, \UUg, \UUl, \oCS^{\geq}, \oCS^{\leq}, \CS^{\geq}, \CS^{\leq}, \CV^{\geq}, \CV^{\leq}.
\end{equation}
All the pairings \eqref{eqn:first pair}--\eqref{eqn:last pair} extend to bialgebra pairings between the bialgebras \eqref{eqn:extended}, in the sense that relations \eqref{eqn:bialg 1}--\eqref{eqn:bialg 2} hold. 

\medskip

\subsection{}
\label{sec:hopf algebras}

Let us now introduce some general notation and terminology pertaining to bialgebras, that will be used throughout the paper. We will work over a base field $\BK$, endowed with an element $q \in \BK^\times$ which is not a root of unity. All our algebras will be unital, associative algebras over $\BK$. All our coalgebras will be counital and coassociative over $\BK$, and all our bialgebras $A$ will be such that the coproduct
$$
\Delta : A \rightarrow A \woo A
$$
is an algebra homomorphism. The hat on top of the $\otimes$ sign in the formula above denotes completion with respect to a certain topology; in all cases studied in the present paper, it will be obvious which completion one needs to take in order to obtain a well-defined coproduct. Thus, we will henceforth drop the term ``topological" in front of notions such as ``coproduct", ``bialgebra" etc. 

\medskip

\subsection{} 
\label{sub:bialg pairing}

Given bialgebras $A^{\geq}$ and $A^{\leq}$, we call
\begin{equation}
\label{eqn:bialg pairing}
A^{\geq} \otimes A^{\leq} \xrightarrow{\langle \cdot, \cdot \rangle} \BK
\end{equation}
a bialgebra pairing if it satisfies the properties
\begin{align}
&\Big \langle x, yy' \Big \rangle = \Big \langle \Delta(x), y \otimes y' \Big \rangle \label{eqn:bialg 1} \\
&\Big \langle xx', y \Big \rangle = \Big \langle x \otimes x', \Delta^{\op}(y) \Big \rangle \label{eqn:bialg 2}
\end{align}
for any $x,x' \in A^{\geq}$, $y,y' \in A^{\leq}$, where $\Delta^{\op}$ denotes the opposite coproduct. The radical(s) of the pairing \eqref{eqn:bialg pairing} are the subsets
\begin{align}
&A^{\geq} \supset I^{\geq} = \left\{x \in A^{\geq} \text{ s.t. } \Big \langle x, A^{\leq} \Big \rangle = 0 \right\} \\
&A^{\leq} \supset I^{\leq} = \left\{y \in A^{\leq} \text{ s.t. } \Big \langle A^{\geq} , y \Big \rangle = 0 \right\}
\end{align}
Because of \eqref{eqn:bialg 1} and \eqref{eqn:bialg 2}, it is easy to see that $I^{\geq}, I^{\leq}$ are two-sided ideals. 

\medskip

\begin{definition}
\label{def:non-deg}

If $I^{\geq} = 0$ (respectively $I^{\leq} = 0$), then we call the pairing \eqref{eqn:bialg pairing} non-degenerate in the first (respectively second) argument.

\end{definition}

\medskip

\subsection{} Given a bialgebra pairing such as \eqref{eqn:bialg pairing}, one defines the Drinfeld double as
\begin{equation}
\label{eqn:dd def}
A = A^{\geq} \otimes A^{\leq}
\end{equation}
where the multiplication between the tensor factors is governed by the relation
\begin{equation}
\label{eqn:dd}
\sum_i \sum_j x_{1,i} y_{1,j} \Big \langle x_{2,i}, y_{2,j} \Big \rangle = \sum_i \sum_j \Big \langle x_{1,i}, y_{1,j} \Big \rangle y_{2,j}, x_{2,i}
\end{equation}
for any $x \in A^{\geq}$ and $y \in A^{\leq}$ with $\Delta(x) = \sum_i x_{1,i} \otimes x_{2,i}$ and  $\Delta(y) = \sum_j y_{1,j} \otimes y_{2,j}$. 

\medskip

\begin{remark} Relation \eqref{eqn:dd} is a convenient reformulation (the author first encountered it in \cite[Section 3.1]{BS}) of the well-known Drinfeld double commutation relation. Let us explain how formula \eqref{eqn:dd} defines an algebra structure on the vector space $A$ of \eqref{eqn:dd def}. First of all, one identifies 
\begin{align*}
&x \in A^{\geq} \quad \text{with} \quad x \otimes 1 \in A \\
&y \in A^{\leq} \quad \text{with} \quad 1 \otimes y \in A
\end{align*}
for all $x \in A^{\geq}$ and $y \in A^{\leq}$. Then one defines for all $x,x' \in A^{\geq}$ and $y,y' \in A^{\leq}$
$$
(x' \otimes y) \cdot (x \otimes y') = \sum_i \sum_j \Big \langle S^{-1}(x_{1,i}), y_{1,j} \Big \rangle x' x_{2,i} \otimes y_{2,j}y' \Big \langle x_{3,i}, y_{3,j}\Big \rangle
$$
Indeed, if $\Delta^{(2)}(x) = \sum_i x_{1,i} \otimes x_{2,i} \otimes x_{3,i}$, $\Delta^{(2)}(y) = \sum_j y_{1,j} \otimes y_{2,j} \otimes y_{3,j}$, the relation
$$
\sum_i \sum_j \Big \langle S^{-1}(x_{1,i}), y_{1,j} \Big \rangle x_{2,i} y_{2,j} \Big \langle x_{3,i}, y_{3,j} \Big\rangle = yx
$$
is equivalent to \eqref{eqn:dd} by general Hopf algebra properties. Above, $S: A^{\geq} \rightarrow A^{\geq}$ denotes the antipode map, which exists and satisfies all the required properties for all bialgebras considered in the present paper (we will not write it down explicitly). 

\end{remark}

\medskip 

\noindent Once \eqref{eqn:dd def} is made into an algebra as above, it is easily seen to be a Hopf algebra by requiring that the natural inclusions $A^{\geq}, A^{\leq} \rightarrow A$ be Hopf algebra homomorphisms.

\medskip

\section{Shuffle algebras} 
\label{sec:shuffle}

\medskip

\noindent In the present Section, we will define and study the shuffle algebras $\CV^\pm$ (and their important subalgebras $\CS^\pm$), which will provide models for quantum loop groups. 

\medskip

\subsection{} 
\label{sub:def quad}

Fix a symmetric Cartan matrix $C$ as in \eqref{eqn:cartan}, and set
\begin{equation}
\label{eqn:def zeta diff}
\zeta_{ij}(x) = \frac {x-q^{-d_{ij}}}{x-1}
\end{equation}
for all $i,j \in I$. We will now recall the algebra $\tUUp$ of \eqref{eqn:def quad intro}. 

\medskip

\begin{definition}
\label{def:quad}

Let $\tUUp = \BK\langle e_{i,k} \rangle_{i\in I, k\in \BZ}$, modulo the relation
\begin{equation}
\label{eqn:rel quad 1}
e_i(z) e_j(w) \zeta_{ji} \left(\frac wz\right) = e_j(w) e_i(z) \zeta_{ij} \left(\frac zw\right)
\end{equation}
where
$$
e_i(z) = \sum_{k \in \BZ} \frac {e_{i,k}}{z^k}
$$
The meaning of \eqref{eqn:rel quad 1} is that one clears out the denominators of the $\zeta$ functions, and then equates the coefficients of any $z^k w^l$ in the two sides of the relation.

\end{definition}

\medskip

\noindent Similarly, we define $\tUUm = \BK\langle f_{i,k} \rangle_{i\in I, k\in \BZ}$ modulo the relation
\begin{equation}
\label{eqn:rel quad 2}
f_i(z) f_j(w) \zeta_{ij} \left(\frac zw\right) = f_j(w) f_i(z) \zeta_{ji} \left(\frac wz\right)
\end{equation}
where
$$
f_i(z) = \sum_{k \in \BZ} \frac {f_{i,k}}{z^k}
$$
It is easy to see that $e_{i,k} \mapsto f_{i,-k}$ induces an isomorphism 
$$
\tUUp \xrightarrow{\sim} \tUUm
$$
Note that the algebra $\tUUpm$ is graded by $\pm \nn \times \BZ$, via
\begin{equation}
\label{eqn:grading u}
\deg e_{i,k} = (\bs^i, k) \quad \text{and} \quad \deg f_{i,-k} =(-\bs^i, -k)
\end{equation}
for all $i \in I$, $k \in \BZ$. Above and henceforth, $\BN$ is considered to contain 0, and $\bs^i \in \nn$ denotes the $I$-tuple of integers with a 1 on position $i$ and 0 everywhere else. Let
\begin{equation}
\label{eqn:total degree}
|\bn| = \sum_{i \in I} n_i
\end{equation}
for any $\bn = (n_i)_{i \in I} \in \zz$.

\medskip

\subsection{}
\label{sub:def shuf}

Fix a total order $<$ on the finite set $I$. Consider an infinite collection of variables $z_{i1},z_{i2},\dots$ for all $i \in I$. For any $\bn = (n_i)_{i \in I} \in \BN^I$, let $\bn! = \prod_{i \in I} n_i!$.

\medskip

\begin{definition}
\label{def:big shuf}

(\cite{E,FO}) The \textbf{big shuffle algebra} associated to $C$ is
\begin{equation}
\label{eqn:big shuf}
\CV^+ = \bigoplus_{\bn \in \BN^I} \frac {\BK[z_{i1}^{\pm 1},\dots,z_{in_i}^{\pm 1}]_{i \in I}^{\esym}}{\prod_{i<j \in I} \prod_{a \leq n_i, b \leq n_j} (z_{ia} - z_{jb})}
\end{equation}
endowed with the multiplication
\begin{equation}
\label{eqn:shuf prod}
R(\dots,z_{i1},\dots,z_{in_i},\dots) * R'(\dots,z_{i1},\dots,z_{in_i'},\dots ) = \frac 1{\bn! \bn'!} \cdot
\end{equation}
$$
\eSym \left[R(\dots,z_{i1},\dots,z_{in_i},\dots) R'(\dots,z_{i,n_i+1},\dots,z_{i,n_i+n_i'},\dots ) \mathop{\prod^{i,j \in I}_{1\leq a\leq n_i}}_{n_j < b \leq n_j+n_j'} \zeta_{ij} \left(\frac {z_{ia}}{z_{jb}} \right) \right]
$$
Above and henceforth, ``\emph{sym}" (resp. ``\emph{Sym}") denotes symmetric functions (resp. symmetrization) with respect to the variables $z_{i1},z_{i2},\dots$ for each $i \in I$ separately. 

\end{definition}

\medskip

\noindent Define $\CV^-$ to be the same vector space as \eqref{eqn:big shuf}, but endowed with the multiplication given by the analogue of formula \eqref{eqn:shuf prod} with
$$
\zeta_{ij} \left(\frac {z_{ia}}{z_{jb}} \right) \quad \text{replaced by} \quad \zeta_{ji} \left(\frac {z_{jb}}{z_{ia}} \right)
$$
The algebra $\CV^\pm$ is graded by $\pm \nn \times \BZ$, via the following assignment for any $R^\pm \in \CV^\pm$:
$$
\deg R^\pm(\dots, z_{i1}, \dots, z_{in_i}, \dots) = (\pm \bn, \text{hom deg }R^\pm)
$$
where ``$\text{hom deg}$" denotes total homogeneous degree. Denote the graded pieces by
$$
\CV^\pm = \bigoplus_{\bn \in \nn} \CV_{\pm \bn} = \bigoplus_{(\bn,k) \in \nn \times \BZ} \CV_{\pm \bn,k}
$$

\medskip

\subsection{} 

The following are straightforward exercises, which we leave to the reader.

\medskip

\begin{proposition}
\label{prop:easy 1}

There is an algebra isomorphism
$$
\CV^+ \xrightarrow{\sim} \CV^-
$$
given by $R(\dots,z_{ia},\dots) \mapsto R(\dots,z^{-1}_{ia},\dots)$.

\end{proposition}

\medskip

\begin{proposition}
\label{prop:easy 2}

There exist homomorphisms of $\pm \nn \times \BZ$ graded algebras
\begin{align}
&\tUUp \stackrel{\widetilde{\Upsilon}^+}\longrightarrow \CV^+, \qquad e_{i,k} \mapsto z_{i1}^k \label{eqn:upsilon+} \\
&\tUUm \stackrel{\widetilde{\Upsilon}^-}\longrightarrow \CV^-, \qquad f_{i,k} \mapsto z_{i1}^k \label{eqn:upsilon-}
\end{align}

\end{proposition}

\medskip

\noindent The maps $\widetilde{\Upsilon}^\pm$ are neither injective nor surjective, and one of the main purposes of the present paper is to describe
\begin{equation}
\label{eqn:image}
\oCS^\pm = \text{Im }\widetilde{\Upsilon}^\pm
\end{equation}
(i.e. $\oCS^\pm$ is the $\BK$-subalgebra of $\CV^\pm$ generated by $\{z_{i1}^k\}_{i \in I, k\in \BZ}$) and
\begin{equation}
\label{eqn:kernel}
K^\pm = \text{Ker }\widetilde{\Upsilon}^\pm
\end{equation}
as a two-sided ideal of $\tUUpm$. 

\medskip

\subsection{}
\label{sub:ext}

Let us now enlarge the algebras $\tUUpm$ and $\CV^\pm$, with the goal of making them into bialgebras. We follow the procedure in Subsection \ref{sub:quantum loop group intro}, items \emph{(2)--(3)}. 

\medskip

\begin{definition}
\label{def:ext}

Consider the algebras
\begin{align}
&\tUUg = \frac {\tUUp [h_{i,0}^{\pm 1}, h_{i,1}, h_{i,2},\dots]_{i \in I}}{\left( h_i(z) e_j(w) = e_j(w) h_i(z) \frac {\zeta_{ij} \left(\frac zw\right)}{\zeta_{ji} \left(\frac wz\right)} \right)} \label{eqn:extended rel 1} \\
&\tUUl = \frac {\tUUm [{h'}^{\pm 1}_{i,0}, h'_{i,-1}, h'_{i,-2},\dots]_{i \in I}}{\left( h'_i(z) f_j(w) = f_j(w) h_i'(z) \frac {\zeta_{ji} \left(\frac wz\right)}{\zeta_{ij} \left(\frac zw\right)} \right)}\label{eqn:extended rel 2}
\end{align}
where $h_i(z) = \sum_{k=0}^{\infty} h_{i,k} z^{-k}$ and $h_i'(z) = \sum_{k=0}^{\infty} h'_{i,-k} z^{k}$. Similarly, define
\begin{align}
&\CV^{\geq} = \frac {\CV^+ [h_{i,0}^{\pm 1}, h_{i,1}, h_{i,2},\dots]_{i \in I}}{\left( h_i(z) R(\dots,z_{jb}, \dots) = R(\dots,z_{jb}, \dots) h_i(z) \prod_{(j,b)} \frac {\zeta_{ij}\left(\frac z{z_{jb}} \right)}{\zeta_{ji}\left(\frac {z_{jb}}z \right)} \right)} \label{eqn:extended rel 3} \\
&\CV^{\leq} = \frac {\CV^- [{h'}^{\pm 1}_{i,0}, h'_{i,-1}, h'_{i,-2},\dots]_{i \in I}}{\left( h'_i(z) R(\dots,z_{jb}, \dots) = R(\dots,z_{jb}, \dots) h_i'(z) \prod_{(j,b)} \frac {\zeta_{ji}\left(\frac {z_{jb}}z \right)}{\zeta_{ij}\left(\frac z{z_{jb}} \right)} \right)} \label{eqn:extended rel 4}
\end{align}
One makes sense of the denominators of \eqref{eqn:extended rel 1}, \eqref{eqn:extended rel 2}, \eqref{eqn:extended rel 3}, \eqref{eqn:extended rel 4} by expanding in negative, positive, negative, positive powers of the variable $z$, respectively.

\end{definition}

\medskip

\noindent The following Proposition is straightforward, and proved just like \cite[Proposition 4.1]{Shuf}, so we leave the details as an exercise to the reader. 

\medskip

\begin{proposition}
\label{prop:bialgebra} 

There are coproducts on $\tUUg$, $\tUUl$, $\CV^{\geq}$, $\CV^{\leq}$, given by the following formulas
$$
\Delta(h_i(z)) = h_i(z) \otimes h_i(z)
$$
$$
\Delta(h'_i(z)) = h'_i(z) \otimes h'_i(z)
$$
$$
\Delta(e_i(z)) = h_i(z) \otimes e_i(z) + e_i(z) \otimes 1
$$
$$
\Delta(f_i(z)) = 1 \otimes f_i(z) + f_i(z) \otimes h_i'(z)
$$
for all $i \in I$, while for all $R^\pm \in \CV^\pm$ we set
\begin{multline*}
\Delta (R^+(\dots,z_{i1},\dots,z_{in_i},\dots)) = \\
= \sum_{\{k_i \in \{0,\dots,n_i\}\}_{i \in I}} \frac {\left[ \prod^{j \in I}_{k_j < b \leq n_j} h_j(z_{jb}) \right] \cdot R^+(\dots, z_{i1},\dots , z_{ik_i} \otimes z_{i,k_i+1}, \dots, z_{in_i},\dots)}{\prod^{i \in I}_{1\leq a \leq k_i} \prod^{j \in I}_{k_j < b \leq n_j} \zeta_{ji} \left( \frac {z_{jb}}{z_{ia}} \right)}
\end{multline*}
\begin{multline*}
\Delta (R^-(\dots,z_{i1},\dots,z_{in_i},\dots)) = \\
= \sum_{\{k_i \in \{0,\dots,n_i\}\}_{i \in I}} \frac {R^-(\dots, z_{i1},\dots , z_{ik_i} \otimes z_{i,k_i+1}, \dots, z_{in_i},\dots) \cdot \left[ \prod^{i \in I}_{1 \leq a \leq k_i} h'_i(z_{ia}) \right]}{\prod^{i \in I}_{1\leq a \leq k_i} \prod^{j \in I}_{k_j < b \leq n_j} \zeta_{ij} \left( \frac {z_{ia}}{z_{jb}} \right)}
\end{multline*}
In the latter two formulas above, we expand the denominator as a power series for $|z_{ia}| \ll |z_{jb}|$, and place all the powers of $z_{ia}$ to the left of the $\otimes$ sign and all the powers of $z_{jb}$ to the right of the $\otimes$ sign (for all $i,j \in I$, $1 \leq a \leq k_i$, $k_j < b \leq n_j$). 

\end{proposition}

\medskip

\noindent It is easy to see that \eqref{eqn:upsilon+}--\eqref{eqn:upsilon-} can be extended to bialgebra homomorphisms: 
$$
\tUUg \xrightarrow{\widetilde{\Upsilon}^{\geq}} \CV^{\geq} \qquad \text{and} \qquad \tUUl \xrightarrow{\widetilde{\Upsilon}^{\leq}} \CV^{\leq}
$$
by sending $h_i(z)$ to $h_i(z)$ and $h_i'(z)$ to $h_i'(z)$.

\medskip

\subsection{} We will write $Dz = \frac {dz}{2\pi i z}$ for any variable $z$. For any homogeneous rational function $F(z,w)$ with coefficients in any field $\BK$ of characteristic 0, the notation
\begin{equation}
\label{eqn:contour integral 2 vars}
\int_{|z| \gg |w|} F(z,w) Dz Dw \qquad \left(\text{respectively} \int_{|z| \ll |w|} F(z,w) Dz Dw \right)
\end{equation}
will refer to the constant term in the expansion of $F$ as a power series in $\frac wz$ (respectively $\frac zw$). The notation is motivated by the situation when $\BK = \BC$, and the constant terms above can be calculated as contour integrals, with the variables $z$ and $w$ running over concentric circles that are very far away from each other. Let us further assume that the rational function $F(z,w)$ only has poles of the form $z - wq^k$ for various $k \in \BZ \backslash 0$. Then we may write for all such $k$
$$
F(z,w) = F_{k}(z,w) + \sum_{d=1}^{\infty} \frac {F_{k,d}(w)}{\left(1 - \frac {wq^k}z\right)^d}
$$
where $F_{k}(z,w)$ are rational functions that are regular at $z=wq^k$, and $F_{k,d}(w)$ are Laurent polynomials, only finitely many of which are non-zero. Then we may formally define the residues of the rational function $F(z,w)$ as
$$
\underset{z = wq^k}{\text{Res}} F(z,w) = F_{k,1}(w)
$$
Since everything above is done in the formal setting (wherein $\BK$ is an arbitrary field), the following analogue of Cauchy's residue theorem still holds
$$
\int_{|z| \gg |w|} F(z,w) Dz Dw - \int_{|z| \ll |w|} F(z,w) Dz Dw = \sum_{k \in \BZ \backslash 0} \int \left[ \underset{z = wq^k}{\text{Res}} F(z,w) \right] Dw
$$
Thus, we may define
\begin{equation}
\label{eqn:residue theorem}
\begin{split}
\qquad \quad \int_{|z| = |w|} F(z,w) Dz Dw  &= \int_{|z| \gg |w|} F(z,w) Dz Dw - \sum_{k > 0} \int \left[ \underset{z = wq^k}{\text{Res}} F(z,w) \right] Dw \\
&= \int_{|z| \ll |w|} F(z,w) Dz Dw + \sum_{k < 0} \int \left[ \underset{z = wq^k}{\text{Res}} F(z,w) \right] Dw
\end{split}
\end{equation}
The definition above is motivated by the situation $\BK = \BC$ and $|q|>1$, in which the left-hand side of \eqref{eqn:residue theorem} is the contour integral of $F(z,w)$ where both variables run over the same circle centered at the origin. Even though all our results are purely algebraic and thus hold for any generic $q$, changing the inequality on $|q|$ would require us to change the contours of integration in the upcoming proofs, hence it would be beneficial to keep the assumption $|q|>1$ in mind throughout the paper. Similar formulas apply for rational functions in several variables, for example
\begin{equation}
\label{eqn:contour integral}
\int_{|z_1| \gg \dots \gg |z_n|} F(z_1,\dots,z_n) \prod_{a=1}^n Dz_a
\end{equation}
will denote the constant term in the expansion of $F(z_1,\dots,z_n)$ as a power series in 
$$
\frac {z_2}{z_1}, \dots, \frac {z_n}{z_{n-1}}.
$$
One can also define $\int_{|z_1|=\dots=|z_n|}$ by analogy with \eqref{eqn:residue theorem}, and we recommend keeping the assumption $\BK=\BC$ and $|q|>1$ in mind to make sense of this notation.

\medskip

\subsection{}

Using the notation above, we now define two very important pairings.

\medskip

\begin{definition}
\label{def:pair}

There exist bialgebra pairings
\begin{align}
&\tUUg \otimes \CV^{\leq} \xrightarrow{\langle \cdot, \cdot \rangle_{\tU\CV}} \BK \label{eqn:pair 1} \\
&\CV^{\geq} \otimes \tUUl \xrightarrow{\langle \cdot, \cdot \rangle_{\CV\tU}} \BK\label{eqn:pair 2}
\end{align}
given by
$$
\Big \langle h_i(z), h_j'(w) \Big \rangle_{\tU\CV} = \Big \langle h_i(z), h_j'(w) \Big \rangle_{\CV\tU} = \frac {\zeta_{ij} \left(\frac zw \right)}{\zeta_{ji} \left(\frac wz \right)} 
$$
while for all $R^\pm \in \CV_{\pm \bn}$ and all $i_1,\dots,i_n \in I$, $k_1,\dots,k_n \in \BZ$, we set
\begin{multline}
\label{eqn:pair formula 1}
\Big \langle e_{i_1,k_1} \dots e_{i_n,k_n}, R^- \Big \rangle_{\tU\CV} = \\ =  \int_{|z_1| \gg \dots \gg |z_n|} \frac {z_1^{k_1}\dots z_n^{k_n} R^-(z_1,\dots,z_n)}{\prod_{1\leq a < b \leq n} \zeta_{i_bi_a} \left(\frac {z_b}{z_a} \right)} \prod_{a=1}^n Dz_a  
\end{multline}
\begin{multline}
\label{eqn:pair formula 2}
\Big \langle R^+, f_{i_1,-k_1} \dots f_{i_n,-k_n} \Big \rangle_{\CV\tU} = \\ =  \int_{|z_1| \ll \dots \ll |z_n|} \frac {z_1^{-k_1}\dots z_n^{-k_n} R^+(z_1,\dots,z_n)}{\prod_{1\leq a < b \leq n} \zeta_{i_ai_b} \left(\frac {z_a}{z_b} \right)} \prod_{a=1}^n Dz_a 
\end{multline}
if $\bs^{i_1}+\dots +\bs^{i_n} = \bn$, and 0 otherwise. Implicit in the notation \eqref{eqn:pair formula 1}--\eqref{eqn:pair formula 2} is that we plug the variable $z_a$ into one of the variables $z_{i_a \bullet_a}$ of $R^\pm$, for every $a \in \{1,\dots,n\}$. The choice of such $\bullet_a$ is immaterial, due to the symmetry of $R^\pm$. 

\end{definition}

\medskip

\noindent The proof follows that of \cite[Proposition 4.2]{Shuf} very closely, so we leave the details as an exercise to the reader. The following result will be proved in Section \ref{sec:words}.

\medskip

\begin{proposition}
\label{prop:non-deg}

The restriction of the pairings \eqref{eqn:pair 1}--\eqref{eqn:pair 2} to
\begin{equation}
\label{eqn:pair restricted}
\tUUp \otimes \CV^{-} \xrightarrow{\langle \cdot, \cdot \rangle_{\tU\CV}} \BK \qquad \text{and} \qquad \CV^{+} \otimes \tUUm \xrightarrow{\langle \cdot, \cdot \rangle_{\CV\tU}} \BK
\end{equation}
are non-degenerate in both arguments. 

\end{proposition}

\medskip

\subsection{}  
\label{sub:zig-zag}

We will now define certain subalgebras of $\CV^\pm$ which pair trivially with the kernels \eqref{eqn:kernel}. For any $i \in I$ and any integers $s \leq t$ congruent modulo 2, we write
$$
\{s,\dots,t\}_i = \Big\{ s, s+2,\dots,t-2, t \Big\}
$$
(the index $i \in I$ will play an important role shortly). If $t = s-2$, the progression above is defined to be empty. Given $i\neq j$ in $I$, a pair
$$
Z = \Big( \{s,\dots,t\}_i , \{s',\dots,t'\}_j \Big)
$$
will be called a \textbf{zig-zag} (the choice of $i$ and $j$ will be part of the datum of the zig-zag). To such a zig-zag, we will associate the number
\begin{multline}
\label{eqn:number}
m_Z = \min \left( \Big| \Big\{ (a,b) \in \{s,\dots,t\}_i \times \{s',\dots,t'\}_j \text{ s.t. } a = b + d_{ij} \Big\} \Big|, \right. \\ \left. \Big| \Big\{(a,b) \in \{s,\dots,t\}_i \times \{s',\dots,t'\}_j \text{ s.t. } a = b - d_{ij} \Big\}\Big| \right) \geq 0
\end{multline}
It can be represented graphically as the minimum of the number of northwest-pointing arrows and the number of southwest-pointing arrows in Figure 4.

\begin{figure}[h]
\centering
\includegraphics[scale=0.25]{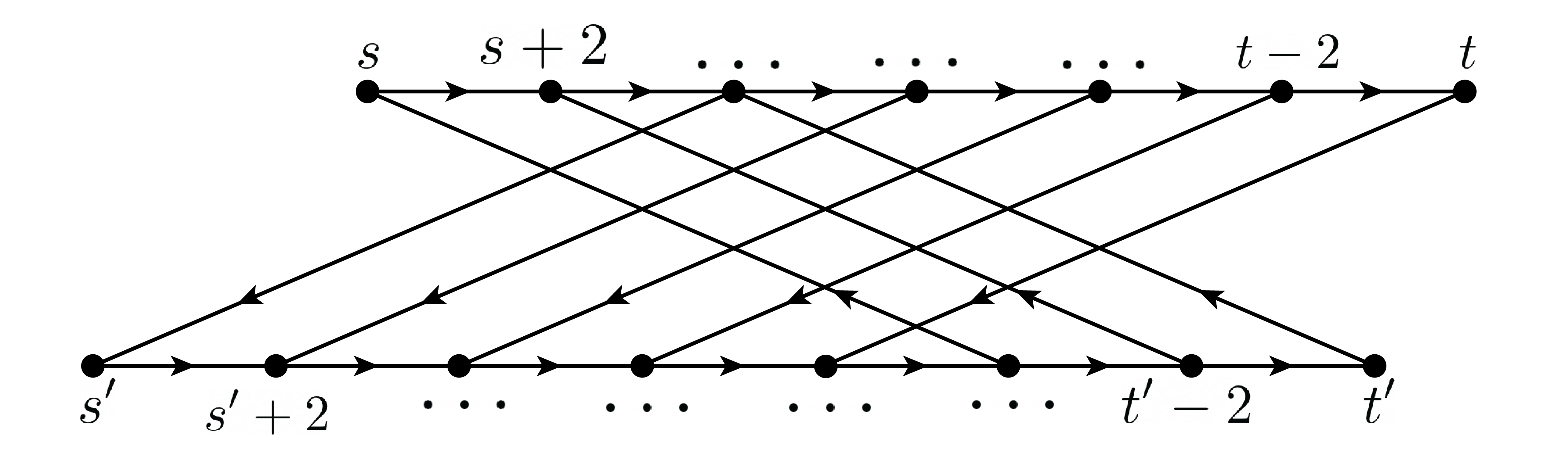} 
\caption{A zig-zag $Z$ for $d_{ij} = -7$ ($m_Z = 3$)}
\end{figure}

\medskip

\noindent Note that $m_Z = 0$ is allowed, and it means the absence of diagonal arrows. Figure 4 above explains the terminology ``zig-zag". We will call a zig-zag \textbf{minimal} if it takes the form in Figure 5 below.

\begin{figure}[h]
\centering
\includegraphics[scale=0.25]{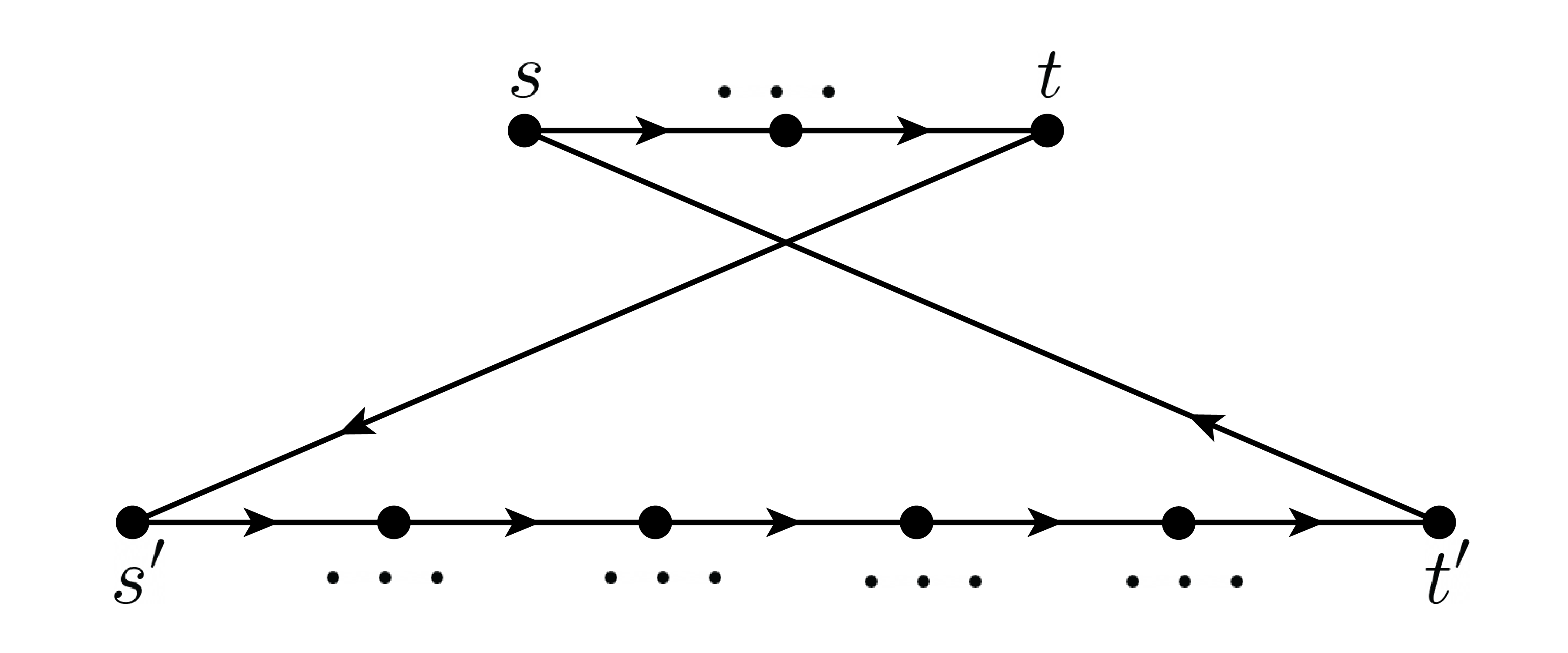} 
\caption{A minimal zig-zag $Z$, for $d_{ij} = -7$, $k=2$, $l=5$}
\end{figure}

\noindent for some non-negative integers $k, l \geq 0$ such that $k+l=-d_{ij}$. In other words, minimal zig-zags have the property that 
$$
t = s+2k, \quad t'=s'+2l, \quad t'=s-d_{ij}, \quad t = s'-d_{ij}
$$
It is easy to see that for a minimal zig-zag, the right-hand side of \eqref{eqn:number} is $\min(1,1) = 1$. We will call a zig-zag \textbf{distinguished} if it can be obtained by ``repeating" a minimal zig-zag $m$ times, for some natural number $m \geq 1$

\begin{figure}[h]
\centering
\includegraphics[scale=0.25]{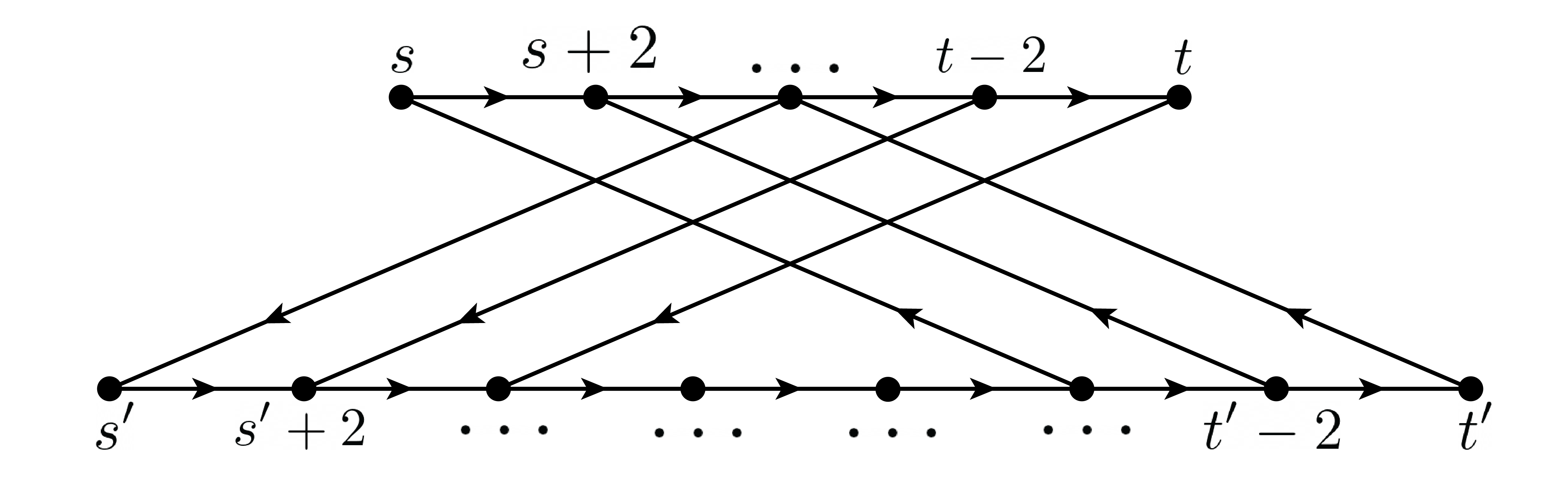} 
\caption{A distinguished zig-zag $Z$, corresponding to ``repeating" the minimal zig-zag in Figure 5 a number of $m=3$ times}
\end{figure}

\noindent In other words, distinguished zig-zags have the property that 
\begin{align}
&t = s+2(k+m-1), \qquad t'=s'+2(l+m-1) \label{eqn:distinguished 1} \\
&t'=s-d_{ij}+2(m-1),  \quad t = s'-d_{ij} +2(m-1) \label{eqn:distinguished 2}
\end{align}
It is easy to see that for such a distinguished zig-zag $Z$, the right-hand side of \eqref{eqn:number} is $\min(m,m) = m$. We will call the number $m$ the \textbf{multiplicity} of $Z$.

\medskip

\subsection{}
\label{sub:wheel}

Using the combinatorics developed in the preceding Subsection, we will now define certain important subalgebras of the shuffle algebra. 

\medskip

\begin{definition}
\label{def:wheel}

Consider the vector space $\CS^\pm \subset \CV^\pm$ consisting of symmetric rational functions (recall that $r$ in formula \eqref{eqn:symmetric rational function} is a Laurent polynomial)
\begin{equation}
\label{eqn:symmetric rational function}
R(\dots, z_{ia}, \dots) = \frac {r(\dots, z_{ia}, \dots)}{\prod_{i < j} \prod_{a \leq n_i, b \leq n_j} (z_{ia} - z_{jb})}
\end{equation}
where for any distinguished zig-zag of multiplicity $m$ (as in Figure 6), we have
\begin{equation}
\label{eqn:wheel}
(x-y)^m \quad \text{divides} \quad \mathop{r\Big|_{z_{i1} = xq^{s}, z_{i2} = xq^{s+2}, \dots, z_{i,\frac {t-s}2+1} = xq^{t}}}_{\text{ }\qquad z_{j1} = yq^{s'}, z_{j2} = yq^{s'+2}, \dots, z_{j,\frac {t'-s'}2+1} = yq^{t'}}
\end{equation}

\end{definition}

\medskip

\noindent If $d_{ij} = 0$, then condition \eqref{eqn:wheel} simply says that $r$ is divisible by $z_{ia}-z_{jb}$ for all $a,b$; thus, we could simply cancel these factors against the corresponding factors in the denominator of \eqref{eqn:symmetric rational function}. If $d_{ij} = -1$, condition \eqref{eqn:wheel} for $m=1$ implies conditions \eqref{eqn:wheel} for all $m \geq 2$: this is based on the elementary observation that if a polynomial in variables $z_1,\dots,z_n$ vanishes at $z_a=z_b=z_c$ for various ``triangles" $\{a,b,c\} \subset \{1,\dots,n\}$, then it vanishes at $z_1=\dots=z_n$ to order equal to the maximal number of triangles that pairwise intersect in 0 or 1 elements. However, as soon as $d_{ij} \leq -2$, conditions \eqref{eqn:wheel} for various natural numbers $m$ are independent. For instance, the symmetric polynomial 
\begin{multline*} 
(1+q^2+q^4)z_{i1}z_{i2}z_{i3}-q^2(z_{i1}z_{i2}+z_{i2}z_{i3}+z_{i3}z_{i1})(z_{j1}+z_{j2}+z_{j3})+\\ q^2(z_{i1}+z_{i2}+z_{i3})(z_{j1}z_{j2}+z_{j2}z_{j3}+z_{j3}z_{j1}) - (1+q^2+q^4) z_{j1}z_{j2}z_{j3}
\end{multline*} 
vanishes at $(z_{i1},z_{i2},z_{i3},z_{j1}) = (x,xq^2,xq^4,xq^2)$, $(z_{i1},z_{i2},z_{j1},z_{j2}) = (x,xq^2,x,xq^2)$, $(z_{i1},z_{j1},z_{j2},z_{j3}) = (xq^2,x,xq^2,xq^4)$ but is not divisible by $(x-y)^2$ when specialized at $(z_{i1},z_{i2},z_{i3},z_{j1},z_{j2},z_{j3}) = (x,xq^2,xq^4,y,yq^2,yq^4)$. 

\medskip 

\begin{remark} When $m = 1$ (i.e. the distinguished zig-zag is minimal) and $k = -d_{ij}$, $l = 0$, then property \eqref{eqn:wheel} is precisely the ``wheel condition" of \cite{E, FO}. 
\end{remark}

\medskip

\begin{proposition}
\label{prop:stronger}

\noindent For any $R \in \CS^\pm$, and any zig-zag $Z$ as in Figure 4, we have
\begin{equation}
\label{eqn:wheel strong}
(x-y)^{m_Z} \quad \text{divides} \quad \mathop{r\Big|_{z_{i1} = xq^{s}, z_{i2} = xq^{s+2}, \dots, z_{i,\frac {t-s}2+1} = xq^{t}}}_{\text{ } \qquad z_{j1} = yq^{s'}, z_{j2} = yq^{s'+2}, \dots, z_{j,\frac {t'-s'}2+1} = yq^{t'}}
\end{equation}
where $r$ is the Laurent polynomial in the numerator of \eqref{eqn:symmetric rational function}.

\end{proposition}

\medskip

\begin{proof} We need to show that if property \eqref{eqn:wheel} holds for all distinguished zig-zags, then it holds for all zig-zags. This is an immediate consequence of the following claim: any zig-zag $Z$ contains a distinguished zig-zag of multiplicity $m_Z$. It remains to prove this claim. Let us assume without loss of generality that $m_Z$ is the number of northwest-pointing arrows in Figure 4. It is easy to see that the left-most northwest-pointing arrow $A_{\nwarrow}$ must intersect the left-most southwest-pointing arrow $A_{\swarrow}$. Thus, the two arrows in question determine a minimal zig-zag as in Figure 5. Any other arrows in the zig-zag $Z$ are obtained by translating either $A_{\nwarrow}$ and $A_{\swarrow}$ by an even integer to the right. By assumption, the translates by $2,4,\dots,2(m_Z-1)$ to the right of $A_{\nwarrow}$ and $A_{\swarrow}$ are still contained in the zig-zag $Z$. Therefore, $Z$ contains a distinguished zig-zag of multiplicity $m_Z$.

\end{proof}

\medskip

\begin{proposition}

$\CS^\pm$ is a subalgebra of $\CV^\pm$.

\end{proposition}

\medskip

\begin{proof} Let us prove the case $\pm = +$, as the case $\pm = -$ is analogous. We will actually show that if $R,R' \in \CS^+$, then $R * R'$ satisfies \eqref{eqn:wheel strong} for any zig-zag $Z$ of the form in Figure 4. In fact, we will prove the stronger fact that every term in the symmetrization \eqref{eqn:shuf prod} satisfies the divisibility property \eqref{eqn:wheel strong}. In other words, consider the expression on the second line of \eqref{eqn:shuf prod} and specialize the variables 
\begin{align*}
&z_{i1}, z_{i2}, \dots, z_{i,\frac {t-s}2+1} \ \ \qquad \text{to} \qquad xq^s, xq^{s+2}, \dots, xq^t \\
&z_{j1}, z_{j2}, \dots, z_{j,\frac {t'-s'}2+1} \qquad \text{to} \qquad yq^{s'}, yq^{s'+2}, \dots, yq^{t'} 
\end{align*}
By the nature of the symmetrization, this entails partitioning the variables $z_{i1},z_{i2},\dots$, $z_{j1}, z_{j2},\dots$ into two disjoint sets: those which are arguments of $R$, and those which are arguments of $R'$. But because of the presence of the functions
$$
\zeta_{ii} \left(\frac {z_{ia} \text{ a variable of }R} {z_{ib} \text{ a variable of }R'}\right) \quad \text{and} \quad \zeta_{jj} \left(\frac {z_{ja} \text{ a variable of }R} {z_{jb} \text{ a variable of }R'}\right)
$$
on the second line of \eqref{eqn:shuf prod}, the specialization in question is non-vanishing only if
\begin{align}
&z_{i1},\dots, z_{i\alpha}, z_{j1},\dots,z_{j\beta} & &\text{are variables of }R' \label{eqn:variables r 1} \\
&z_{i,\alpha+1},\dots, z_{i,\frac {t-s}2+1}, z_{j,\beta+1}, \dots, z_{j,\frac {t'-s'}2+1} & &\text{are variables of }R \label{eqn:variables r 2}
\end{align}
for some $\alpha \in \{0,\dots,\frac {t-s}2 + 1 \}$ and $\beta \in \{0,\dots,\frac {t'-s'}2 + 1 \}$. In other words, the zig-zag $Z$ is partitioned into the two ``consecutive" zig-zags
\begin{align}
&Z_1 = \Big( \{s,\dots,u \}_i, \{s',\dots,u' \}_j \Big) \label{eqn:zig-zag 1} \\
&Z_2 = \Big( \{u+2,\dots,t \}_i, \{u'+ 2,\dots, t' \}_j \Big) \label{eqn:zig-zag 2} 
\end{align}
where $u = s+2(\alpha-1)$ and $u' = s'+2(\beta-1)$. However, the presence of 
$$
\zeta_{ij} \left(\frac {z_{ia} \text{ a variable of }R} {z_{jb} \text{ a variable of }R'}\right) \quad \text{and} \quad \zeta_{ji} \left(\frac {z_{ja} \text{ a variable of }R} {z_{ib} \text{ a variable of }R'}\right)
$$
implies that the (specialization of) the second line of \eqref{eqn:shuf prod} is divisible by a factor of $x-y$ for every arrow in the zig-zag that points from one of the variables in \eqref{eqn:variables r 2} to one of the variables in \eqref{eqn:variables r 1} \footnote{It is a little imprecise of us to talk about the specialization of the rational function $R$, when the condition we need to check pertains to the specialization of the Laurent polynomial $r$ that appears in \eqref{eqn:symmetric rational function}; the difference between these lies precisely in the denominators of the functions $\zeta_{ij}$ for $i\neq j$, which is why we ignore factors $(x-y)$ that might arise from these denominators.}. Thus, the required property \eqref{eqn:wheel strong} is an immediate consequence of the fact that
\begin{equation}
\label{eqn:inequality}
m_Z \leq m_{Z_1} + m_{Z_2} + \Big| \Big\{ \text{arrows pointing from }Z_2 \text{ to } Z_1 \Big\} \Big|
\end{equation}
As shown in the proof of Proposition \ref{prop:stronger}, the zig-zag $Z$ contains a pair of intersecting arrows $(A_{\nwarrow}, A_{\swarrow})$ together with their translates by $2,4,\dots,2(m_Z-1)$ units to the right. Any one of these $m_Z$ pairs of arrows either lies completely in $Z_1$, or it lies completely in $Z_2$, or at least one of the arrows in the pair points from $Z_2$ to $Z_1$. This establishes the inequality \eqref{eqn:inequality}.

\end{proof}

\subsection{} It is easy to note that:
$$
\oCS^\pm \subseteq \CS^\pm \subset \CV^\pm
$$
on account of the fact that $\CS^\pm$ are subalgebras, and they contain the generators $\{z_{i1}^k\}_{i \in I, k \in \BZ}$ of $\oCS^\pm$. Our choice of the conditions \eqref{eqn:wheel} that define $\CS^\pm$ was motivated by the following result. 

\medskip

\begin{proposition}
\label{prop:descends}

The pairings \eqref{eqn:pair restricted} trivially pair anything in the kernels \eqref{eqn:kernel} with anything in the subalgebras $\CS^\pm$, i.e.
\begin{equation}
\label{eqn:descend}
\Big \langle K^+, \CS^- \Big \rangle_{\tU\CS} = 0 = \Big \langle \CS^+, K^- \Big \rangle_{\CS\tU}
\end{equation}
Therefore, the pairings \eqref{eqn:pair restricted} descend to
\begin{align}
&\oCS^{+} \otimes \CS^{-} \xrightarrow{\langle \cdot, \cdot \rangle_{\oCS \CS}} \BK \label{eqn:descend pair 1} \\
&\CS^{+} \otimes \oCS^{-} \xrightarrow{\langle \cdot, \cdot \rangle_{\CS \oCS}} \BK \label{eqn:descend pair 2}
\end{align}
which obviously have the same restriction to
\begin{equation}
\label{eqn:descend pair 3}
\oCS^{+} \otimes \oCS^{-} \xrightarrow{\langle \cdot, \cdot \rangle_{\oCS \oCS}} \BK
\end{equation}

\end{proposition}

\medskip

\noindent Akin to Definition \ref{def:ext}, let us define the extended bialgebras
\begin{equation}
\label{eqn:extended bialgebras}
\oCS^\geq, \oCS^\leq \quad \text{and} \quad \CS^\geq, \CS^\leq
\end{equation}
by adding generators $h_{i,k}, h'_{i,k}$ to $\oCS^\pm$ and $\CS^\pm$, as in Definition \ref{def:ext}. Then \eqref{eqn:descend pair 1}--\eqref{eqn:descend pair 2} extend to bialgebra pairings between the various bialgebras in \eqref{eqn:extended bialgebras}. This also explains the word ``obviously" before \eqref{eqn:descend pair 3}: to show that two bialgebra pairings are identical, one need only check that they match on the generators, in which case the check is trivial.

\medskip

\begin{proof} We will prove the first equality in \eqref{eqn:descend}, as the second one is analogous. The strategy is to consider any linear combination for fixed $\bn$
$$
\sum_{\bs^{i_1}+\dots +\bs^{i_n} = \bn} \sum_{k_1,\dots,k_n \in \BZ} \text{coefficient} \cdot e_{i_1,k_1} \dots e_{i_n,k_n} \in K^+ = \text{Ker }\tUpsilon^+
$$
and to show that it pairs trivially with any $R^- \in \CS_{-\bn}$. To this end, let us start from the right-hand side of \eqref{eqn:pair formula 1} and study how the integral changes as we move the contours of integration toward $|z_1| = \dots = |z_n|$. The following explanation is phrased for $\BK = \BC$ and $q \in \BC^*$ satisfying $|q| > 1$; however, this is just a linguistic tool to keep track of various residues involved as we move contours around (since these residues are all rational functions in $q$, they make sense over any field $\BK$). For any $m \in \{1,\dots,n\}$, consider the expression
\begin{equation}
\label{eqn:xm}
X_m = \sum^{\text{fair partition}}_{\{m,\dots,n\}=A_1 \sqcup \dots \sqcup A_t} \int_{|z_1| \gg \dots \gg |z_{m-1}| \gg |w_1| = \dots = |w_t|} 
\end{equation}
$$
\left[ \underset{(z_{a_s^{(1)}}, \dots, z_{a_s^{(n_s)}}) = (w_s q^{n_s-1}, \dots, w_sq^{1-n_s})}{\text{Res}}  \frac {z_1^{k_1}\dots z_n^{k_n} R^-(z_1,\dots,z_n)}{\prod_{1\leq a < b \leq n} \zeta_{i_bi_a} \left(\frac {z_b}{z_a} \right)} \right]_{\forall s \in \{1,\dots,t\}} \prod_{s=1}^t Dw_s \prod_{a=1}^{m-1} Dz_a
$$
In the notation above, a partition is called \textbf{fair} if all of its constituent sets
\begin{equation}
\label{eqn:set of fair partition}
A_s = \{a_s^{(1)} < \dots < a_s^{(n_s)}\}
\end{equation}
have the property that $\iota(A_s) := i_{ a_s^{(1)}}  = \dots = i_{a_s^{(n_s)}}$ for all $s \in \{1,\dots,t\}$, and the notation on the second line of \eqref{eqn:xm} is shorthand for the following iterated residue
$$
\underset{(z_{a_s^{(1)}}, \dots, z_{a_s^{(n_s)}}) = (w_s q^{n_s-1}, \dots, w_sq^{1-n_s})}{\text{Res}} = \underset{z_{a_s^{(1)}} = z_{a_s^{(2)}} q^2}{\text{Res}} \left[ \underset{z_{a_s^{(2)}} = z_{a_s^{(3)}} q^2}{\text{Res}}  \dots \left[ \underset{z_{a_s^{(n_s-1)}} = z_{a_s^{(n_s)}} q^2}{\text{Res}} \right] \dots \right]
$$
followed by relabeling the variable $z_{a_s^{(n_s)}}$ to $w_sq^{1-n_s}$.

\medskip

\begin{remark}
\label{rem:midpoint}

It is a key part of the argument (and substantially different from the analogous setup studied in \cite[Proposition 3.3]{Quiver 1}) that the new variable $w_s$ whose contour is considered in \eqref{eqn:xm} is the geometric mean of the variables $z_{a_s^{(1)}},\dots,z_{a_s^{(n_s)}}$. This has the effect that our geometric progressions are midpoint aligned, as opposed from those of \loccit which were endpoint aligned.

\end{remark}

\medskip

\begin{claim}
\label{claim:induction}

We have $X_m = X_{m-1}$ for all $m \in \{2,\dots,n\}$.

\end{claim}

\medskip

\noindent Let us first show how Claim \ref{claim:induction} implies \eqref{eqn:descend}. By iterating Claim \ref{claim:induction} a number of $n-1$ times, we conclude that $X_n = X_1$, or more explicitly
\begin{equation}
\label{eqn:end induction 1}
\Big \langle e_{i_1,k_1} \dots e_{i_n,k_n}, R^- \Big \rangle_{\tU\CS} = \sum^{\text{fair partition}}_{\{1,\dots,n\}=A_1 \sqcup \dots \sqcup A_t} \int_{|w_1| = \dots = |w_t|} 
\end{equation}
$$
\left[ \underset{(z_{a_s^{(1)}}, \dots, z_{a_s^{(n_s)}}) = (w_s q^{n_s-1}, \dots, w_sq^{1-n_s})}{\text{Res}} \frac {z_1^{k_1}\dots z_n^{k_n} R^-(z_1,\dots,z_n)}{\prod_{1\leq a < b \leq n} \zeta_{i_bi_a} \left(\frac {z_b}{z_a} \right)} \right]_{\forall s \in \{1,\dots,t\}} \prod_{s=1}^t Dw_s 
$$
However, for any fixed fair partition $\{1,\dots,n\} = \bar{A}_1 \sqcup \dots \sqcup \bar{A}_t$, we claim that 
\begin{equation}
\label{eqn:ofofof}
\left[ \mathop{\underset{(z_{\bar{a}_s^{(1)}}, \dots, z_{\bar{a}_s^{(n_s)}}) = }{\text{Res}}}_{(w_s q^{n_s-1}, \dots, w_sq^{1-n_s})} \text{Sym}\left( \frac {z_1^{k_1}\dots z_n^{k_n} R^-(z_1,\dots,z_n)}{\prod_{1\leq a < b \leq n} \zeta_{i_bi_a} \left(\frac {z_b}{z_a} \right)} \right) \right]_{\forall s \in \{1,\dots,t\}} = 
\end{equation}
$$
= \mathop{\sum^{\text{fair partition}}_{\{1,\dots,n\}=A_1 \sqcup \dots \sqcup A_t}}_{|A_s| = |\bar{A}_s|, \iota(A_s) = \iota(\bar{A}_s), \forall s} \left[ \mathop{\underset{(z_{a_s^{(1)}}, \dots, z_{a_s^{(n_s)}})=}{\text{Res}}}_{(w_s q^{n_s-1}, \dots, w_sq^{1-n_s})} \frac {z_1^{k_1}\dots z_n^{k_n} R^-(z_1,\dots,z_n)}{\prod_{1\leq a < b \leq n} \zeta_{i_bi_a} \left(\frac {z_b}{z_a} \right)} \right]_{\forall s \in \{1,\dots,t\}}
$$
In \eqref{eqn:ofofof}, we write $\Sym$ for symmetrization with respect to all pairs of variables $z_a$ and $z_b$ such that $i_a = i_b$, and we denote the elements of $\bar{A}_s$ \underline{in any order} by
\begin{equation}
\label{eqn:choose bar}
\bar{A}_s = \{\bar{a}_s^{(1)} , \dots , \bar{a}_s^{(n_s)}\}
\end{equation}
for all $s \in \{1,\dots,t\}$. Formula \eqref{eqn:ofofof} is an immediate consequence of the fact that the only poles of the integrand involving variables $z_a$ and $z_b$ with $a<b$ and $i_a = i_b$ are $z_a - z_b q^2$. With \eqref{eqn:ofofof} in mind, relation \eqref{eqn:end induction 1} yields
\begin{equation}
\label{eqn:end induction 2}
\Big \langle e_{i_1,k_1} \dots e_{i_n,k_n}, R^- \Big \rangle_{\tU\CS} = \sum^{\text{fixed fair partition}}_{\{1,\dots,n\}=\bar{A}_1 \sqcup \dots \sqcup \bar{A}_t} \int_{|w_1| = \dots = |w_t|} 
\end{equation}
$$
\left[ \underset{(z_{\bar{a}_s^{(1)}}, \dots, z_{\bar{a}_s^{(n_s)}}) = (w_s q^{n_s-1}, \dots, w_sq^{1-n_s})}{\text{Res}} \Sym \left( \frac {z_1^{k_1}\dots z_n^{k_n} R^-(z_1,\dots,z_n)}{\prod_{1\leq a < b \leq n} \zeta_{i_bi_a} \left(\frac {z_b}{z_a} \right)} \right) \right]_{\forall s \in \{1,\dots,t\}} \prod_{s=1}^t Dw_s 
$$
where in the right-hand side of \eqref{eqn:end induction 2}, we choose a fixed fair partition $\bar{A}_1 \sqcup \dots \sqcup \bar{A}_t$ for any unordered sum
$$
\bn = \sum_{s=1}^{t} n_s \cdot \bs^{\iota(\bar{A}_s)}
$$
(the implication being that $\bar{A}_s$ of \eqref{eqn:choose bar} is chosen to be an arbitrary set of $n_s$ variables among $z_{i1},\dots,z_{in_i}$, where $i = \iota(\bar{A}_s)$). However, it is clear that the right-hand side of \eqref{eqn:end induction 2} only depends on
$$
\Sym \left( \frac {z_1^{k_1}\dots z_n^{k_n}}{\prod_{1\leq a < b \leq n} \zeta_{i_bi_a} \left(\frac {z_b}{z_a} \right)} \right) = \frac {\tUpsilon^+(e_{i_1,k_1} \dots e_{i_n,k_n})}{\prod_{a \neq b} \zeta_{i_bi_a} \left(\frac {z_b}{z_a} \right)}
$$
so if a certain linear combination of $e_{i_1,k_1} \dots e_{i_n,k_n}$ lies in $K^+ = \text{Ker }\tUpsilon^+$, then that linear combination will pair trivially with any $R^- \in \CS^-$. This establishes \eqref{eqn:descend}.

\medskip

\noindent It remains to prove Claim \ref{claim:induction}. To this end, consider the residue theorem
$$
\int_{|z| \gg |w|} f(z,w) Dz Dw = \int_{|z| = |w|} f(z,w) Dz Dw + \sum_{|c| > 1} \int \left[ \underset{z = wc}{\text{Res}} f(z,w) \right] Dw
$$
for any homogeneous rational function $f$, all of whose poles are of the form $z - wc$ (see  \eqref{eqn:residue theorem} for the notation). Consider formula \eqref{eqn:xm}, and let us zoom in on the summand corresponding to a given fair partition $\{m,\dots,n\} = A_1 \sqcup \dots \sqcup A_t$. As we move the (larger) contour of the variable $z_{m-1}$ toward the (smaller) contours of the variables $w_1,\dots,w_t$, one of two things can happen. The first thing is that the larger contour reaches the smaller ones, which leads to the fair partition
$$
\{m-1,\dots,n\} = A_1\sqcup \dots \sqcup A_t \sqcup \{m-1\}
$$
in formula \eqref{eqn:xm} for $m$ replaced by $m-1$. The second thing is that the variable $z_{m-1}$ gets ``caught" in a pole of the form $z_{m-1} = w_s c$ for some $s \in \{1,\dots,t\}$ and some $|c|>1$. However, the apparent poles of the rational function on the second line of \eqref{eqn:xm} that involve both $z_{m-1}$ and some $w_s$ for $s \in \{1,\dots,t\}$ are of the form
\begin{equation}
\label{eqn:two cases}
 \begin{cases} \frac 1{z_{m-1}-w_s q^{n_s+1}} &\text{if } i_{m-1} = \iota(A_s) \\ \prod^{\bullet \in \{n_s-1,n_s-3,\dots,3-n_s,1-n_s\}}_{\bullet+d \geq 0} \frac 1{z_{m-1} - w_s q^{\bullet+d}} &\text{if }i_{m-1} \neq \iota(A_s)\end{cases}
\end{equation}
where $d = d_{i_{m-1}\iota(A_s)}$; the reason we only consider $\bullet+d \geq 0$ in the second option is that we only choose to move the contour of the variable $z_{m-1}$ from infinity to $|w_s|$, i.e. the midpoint of the geometric progression $w_sq^{1-n_s},\dots,w_sq^{n_s-1}$. Let us start with the second option above. The apparent simple pole at
$$
z_{m-1} = w_s q^{\bullet+d} 
$$
for some $\bullet \in \{n_s-1,\dots,1-n_s\}$ such that $\bullet+d \geq 0$ is precisely canceled by the fact that $R^-$ vanishes at the specialization
$$
z_{\iota(A_s)1} = w_s q^{\bullet+2d}, z_{\iota(A_s)2} = w_s q^{\bullet+2d+2}, \dots , z_{\iota(A_s), 1-d} = w_s q^{\bullet}
$$
$$
z_{i_{m-1}1} = w_s q^{\bullet+d} 
$$
(all the powers of $q$ on the first line lie in the arithmetic progression $\{n_s-1,\dots,1-n_s\}$) due to the condition \eqref{eqn:wheel} for the distinguished triangle in Figure 3. 

\medskip

\noindent As for the first option in \eqref{eqn:two cases}, it leads to the fair partition:
$$
\{m-1,\dots,n\} = A_1 \sqcup \dots \sqcup A_{s-1} \sqcup \Big( A_s \sqcup \{m-1\} \Big) \sqcup A_{s+1} \sqcup \dots \sqcup A_t
$$
in formula \eqref{eqn:xm} for $m$ replaced by $m-1$. However, there is a catch: in this new fair partition, the variables that correspond to the $s$-th part are specialized to
$$
w_s q^{n_s+1}, w_s q^{n_s-1},\dots, w_s q^{1-n_s}
$$
In order to match this with formula \eqref{eqn:xm} for $m$ replaced by $m-1$, we need to move the contour of the variable $w_s$ from
\begin{equation}
\label{eqn:move}
|w_s| = |w_r|, \ \forall r \neq s \qquad \text{to} \qquad |w_sq| = |w_r|, \ \forall r \neq s
\end{equation}
It remains to show that no new poles involving $w_s$ and $w_r$ (for an arbitrary $r \neq s$) are produced in the rational function
\begin{equation}
\label{eqn:restricted rational function}
\mathop{\underset{(z_{a_r^{(1)}}, \dots, z_{a_r^{(n_r)}}) = (w_r q^{n_r-1}, \dots, w_r q^{1-n_r})}{\text{Res}}}_{(z_{m-1}, z_{a_s^{(1)}}, \dots, z_{a_s^{(n_s)}}) = (w_s q^{n_s+1}, w_s q^{n_s-1}, \dots, w_sq^{1-n_s})} \left[ \frac {z_1^{k_1}\dots z_n^{k_n} R^-(z_1,\dots,z_n)}{\prod_{1\leq a < b \leq n} \zeta_{i_bi_a} \left(\frac {z_b}{z_a} \right)} \right]
\end{equation}
as we move the contours according to \eqref{eqn:move}. Because all the denominators of \eqref{eqn:restricted rational function} in the variables $w_r,w_s$ are of the form $w_s - w_r q^{\text{integer}}$, it suffices to show that the residue \eqref{eqn:restricted rational function} is regular at both $w_sq - w_r$ and $w_s - w_r$. We will deal with these two possible cases below, but note that for any given $r$ and $s$, only one of these cases can occur according to the parity of $d_{\iota(A_r)\iota(A_s)}$.

\medskip

\noindent \emph{Case 1:} regularity at $w_s q - w_r$. Consider the distinguished zig-zag $Z$ given in the following picture

\begin{figure}[h]
\centering
\includegraphics[scale=0.2]{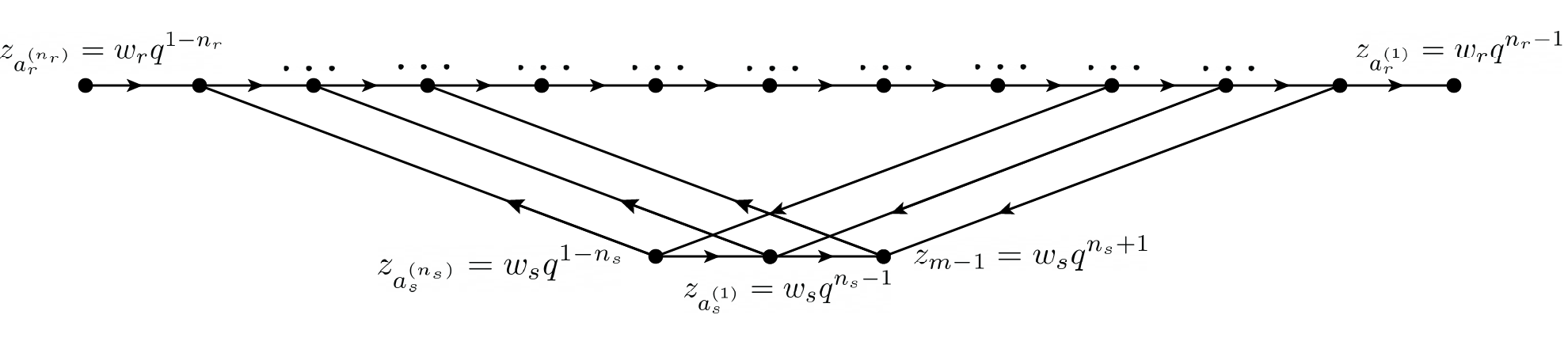} 
\end{figure}

\noindent The denominator of \eqref{eqn:restricted rational function} will contain a factor of $w_sq-w_r$ for every diagonal arrow above whose source $z_b$ and target $z_a$ satisfy $b > a$. However, the order \eqref{eqn:set of fair partition} implies that from every pair of intersecting arrows $A_{\swarrow}$ and $A_{\nwarrow}$, at most one of them could have the property that its source $z_b$ has a larger subscript than its target $z_a$. Since a distinguished zig-zag consists of $m_Z$ disjoint pairs of intersecting arrows $(A_{\swarrow}, A_{\nwarrow})$, we conclude that there can be at most $m_Z$ factors of $(w_sq-w_r)$ in the denominator of \eqref{eqn:restricted rational function}. However, property \eqref{eqn:wheel} implies that the numerator of \eqref{eqn:restricted rational function} is divisible by $(w_sq-w_r)^{m_Z}$, thus canceling the apparent pole.

\medskip

\noindent \emph{Case 2:} regularity at $w_s - w_r$. Consider the zig-zag $Z$ given in the following picture

\begin{figure}[h]
\centering
\includegraphics[scale=0.2]{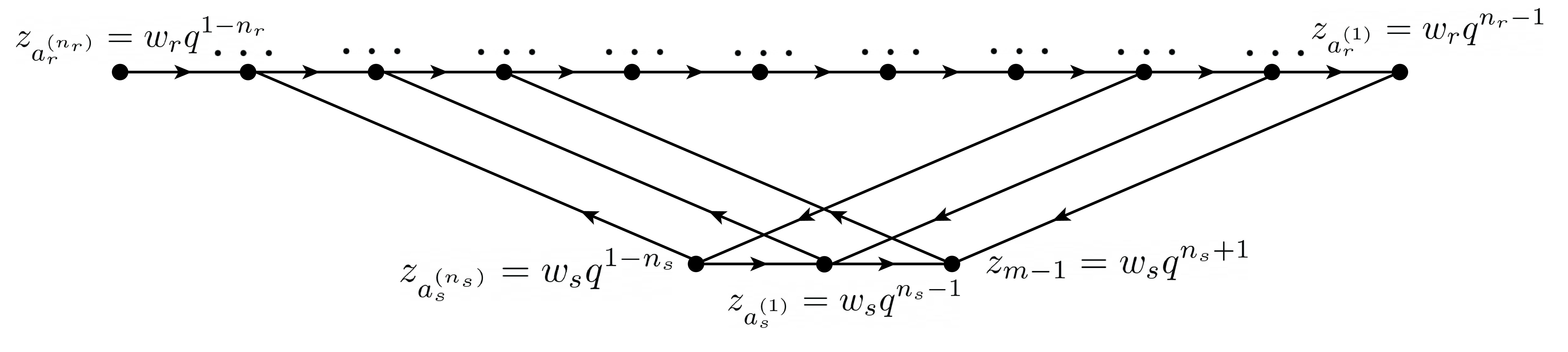} 
\end{figure}

\noindent If we removed the bottom right vertex, the resulting zig-zag $Z'$ would be distinguished. It is easy to see from the picture above that
\begin{equation}
\label{eqn:mzz}
m_{Z} = \begin{cases} m_{Z'}+1 &\text{if } n_s+1-d_{\iota(A_s)\iota(A_r)} \in \{n_r-1,\dots,1-n_r\} \\ m_{Z'} &\text{otherwise} \end{cases}
\end{equation}
The denominator of \eqref{eqn:restricted rational function} will contain at most $m_{Z'}$ factors $w_s-w_r$ which do not involve the bottom right vertex in the Figure above, by the same argument as in \emph{Case 1}. Since only those arrows whose source $z_b$ and target $z_a$ satisfy $b > a$ produce factors $w_s-w_r$ in the denominator, and since the variable $z_{m-1}$ associated to the bottom right vertex has lower subscript than all of the other ones, we see that the bottom right vertex can give rise to at most one factor $w_s-w_r$ and this only happens if $n_s+1-d_{\iota(A_s)\iota(A_r)} \in \{n_r-1,\dots,1-n_r\}$. Thus, the number of factors $w_s-w_r$ in the denominator is at most equal to the right-hand side of \eqref{eqn:mzz}. However, property \eqref{eqn:wheel strong} implies that the numerator of \eqref{eqn:restricted rational function} is divisible by $(w_s-w_r)^{m_Z}$, thus canceling the apparent pole due to \eqref{eqn:mzz}. 

\end{proof}

\begin{remark} In follow-up to Remark \ref{rem:midpoint}, let us note that it was crucial to the argument above that we only moved the contours until all geometric progressions were midpoint aligned. If we tried instead to move the contours so that the geometric progressions were endpoint aligned (as in \cite{Quiver 1}), certain poles of the integrand would hinder us.	
\end{remark}

\medskip

\noindent In Section \ref{sec:words}, we will use the combinatorics of words to prove the following result. 

\medskip

\begin{proposition}
\label{prop:coincide}

We have the identity
\begin{equation}
\label{eqn:coincide}
\oCS^\pm = \CS^\pm
\end{equation}
Moreover, the pairings \eqref{eqn:descend pair 1}--\eqref{eqn:descend pair 2} coincide, thus yielding a pairing
\begin{equation}
\label{eqn:pairing shuffle}
\CS^+ \otimes \CS^- \xrightarrow{\langle \cdot, \cdot \rangle_{\CS\CS}} \BK
\end{equation}
which is non-degenerate in both arguments. 

\end{proposition}

\medskip

\section{Quantum loop groups}
\label{sec:quantum}

\medskip

\noindent The main goal of the present Section is to construct a quotient that fits into the bottom left corner of the diagram below, which will be our desired quantum loop group
\begin{equation}
\label{eqn:desired diagram}
\begin{tikzcd}[column sep=0.05em]
\tUUp \arrow[d,twoheadrightarrow] & \otimes & \CV^- \arrow{rrrrrrrrrrrrrrrrrrrrrrrrrrrrrrrrrrrrrrrrrrrrrrrr}{\langle \cdot, \cdot \rangle_{\widetilde{U}\CV}} & & & & & & & & & & & & & & & & & & & & & & & & & & & & & & & & & & & & & & & & & & & & & & & & \BK\\
\UUp  & \otimes & \CS^- \arrow[u, hook] \arrow{rrrrrrrrrrrrrrrrrrrrrrrrrrrrrrrrrrrrrrrrrrrrrrrru}[swap]{\langle \cdot, \cdot \rangle_{U\CS}}
\end{tikzcd}
\end{equation}
The gist of the construction is that the two vertical arrows above should be dual to each other, i.e. the $\twoheadrightarrow$ map precisely annihilates those elements $\phi$ such that the linear functionals $\langle \phi,-\rangle_{\widetilde{U}\CV}$ cut out the inclusion $\hookrightarrow$. The point will be to describe such $\phi$'s (or more precisely a set of generators for the two-sided ideal of such $\phi$'s) and to also prove that the homomorphism \eqref{eqn:upsilon+} descends to an isomorphism $\Upsilon^+$
$$
\begin{tikzcd}
\tUUp \arrow[d,twoheadrightarrow] \arrow{rd}{\widetilde{\Upsilon}^+}  \\
\UUp \arrow{r}{\sim}[swap]{\Upsilon^+} & \CS^+
\end{tikzcd}
$$
The analogous results also hold with $+$ replaced by $-$.

\medskip

\subsection{} Let us consider a distinguished zig-zag $Z$, as in Subsection \ref{sub:zig-zag}.

\begin{figure}[H]
\centering
\includegraphics[scale=0.25]{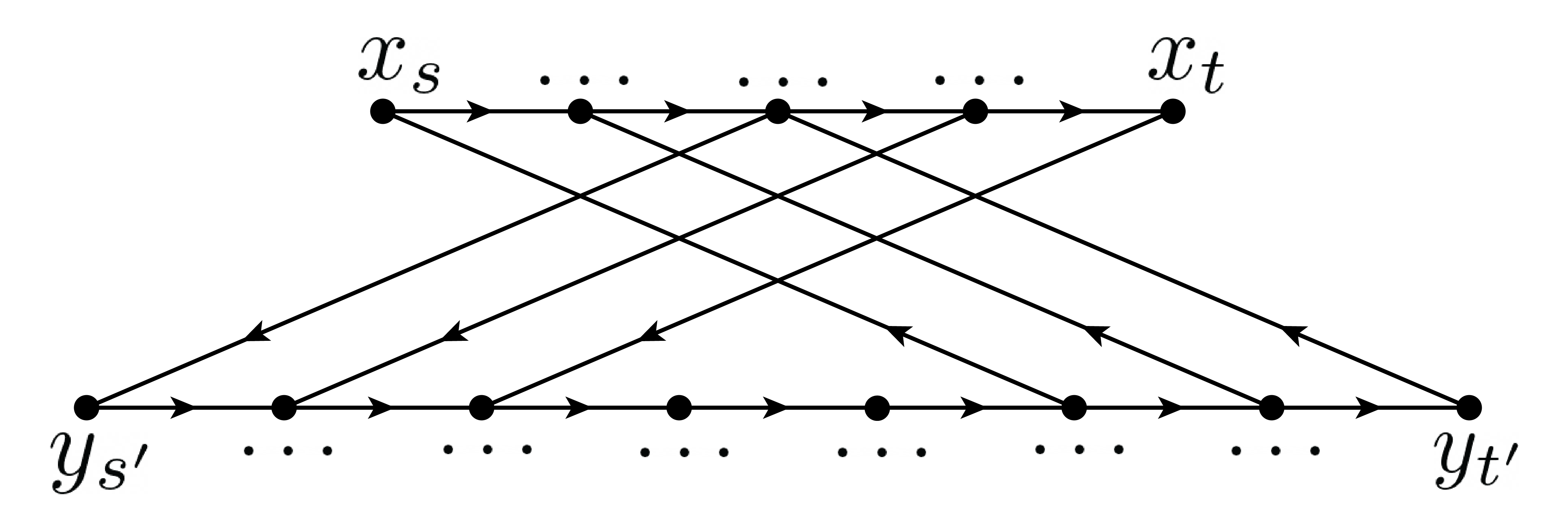} 
\caption{A distinguished zig-zag $Z$, decorated with variables}
\end{figure}

\noindent As an oriented graph, its vertex set will be denoted by $Z_v = \{s,\dots,t\}_i \sqcup \{s',\dots,t'\}_j$, and its edge set will be denoted by $Z_e$. We decorate the vertices of the zig-zag with the variables $x_s, \dots, x_t$ and $y_{s'}, \dots, y_{t'}$, as depicted in Figure 7. Consider the notation
\begin{equation}
\label{eqn:notation 1}
z_c = \begin{cases} x_c &\text{if } c \in \{s,\dots,t\}_i \\ y_c &\text{if } c\in \{s',\dots,t'\}_j \end{cases}
\end{equation}
and 
\begin{equation}
\label{eqn:notation 2}
\iota(c) = \begin{cases} i &\text{if } c \in \{s,\dots,t\}_i \\ j &\text{if } c\in \{s',\dots,t'\}_j \end{cases}
\end{equation}
for any $c \in Z_v = \{s,\dots,t\}_i \sqcup \{s',\dots,t'\}_j$. Then we define
\begin{equation}
\label{eqn:notation 3}
e_{\iota(c)}(z_c) = \begin{cases} e_i(x_c) &\text{if } c \in \{s,\dots,t\}_i \\ e_j(y_c) &\text{if } c\in \{s',\dots,t'\}_j \end{cases}
\end{equation}
Finally, we will write
\begin{equation}
\label{eqn:notation 4}
\barx_c = x_c q^c, \qquad \bary_c = y_c q^c, \qquad \barz_c = z_cq^c
\end{equation}
for all applicable $c \in Z_v$ (the exponent of $q^c$ refers to the underlying integer of $c$).

\medskip

\begin{definition}
\label{def:cycle}

Consider any collection of edges $H \subset Z_e$ which contains no oriented cycle. Choose a total order of the set of vertices $Z_v = \{s,\dots,t\}_i \sqcup \{s',\dots,t'\}_j$ that is compatible with $H$, i.e. there is an arrow only from a larger vertex to a smaller vertex. Then we consider the formal series
$$
e_H\left(x_s, \dots, x_t, y_{s'}, \dots, y_{t'} \right) \in \tUUp [\![ x^{\pm 1}_s,\dots,x^{\pm 1}_t, y^{\pm 1}_{s'}, \dots, y^{\pm 1}_{t'} ]\!]
$$
given by the formula
\begin{equation}
\label{eqn:def eh}
\frac {\prod_{c < c' \in Z_v} \left[ (-1)^{\delta_{\iota(c)i} \delta_{\iota(c')j} + \delta_{\iota(c)\iota(c')} \nu_{c'-c}} (z_c - z_{c'} q^{-d_{\iota(c)\iota(c')}}) \right]}{\prod_{\{c' \rightarrow c\} \in H} \left(\barz_{c'} - \barz_c \right)}  \mathop{\prod_{c \in Z_v \text{ in}}}_{\text{descending order}} e_{\iota(c)}(z_c) 
\end{equation}
where $\delta_{kl}$ is the Kronecker delta symbol for any $k,l \in I$, and $\nu_{c'-c}$ is 1 (respectively 0) if $c'$ is greater (respectively smaller) than $c$ with respect to the usual ordering of the integers, which may differ from the total order $<$. Note that the ratio in \eqref{eqn:def eh} is a Laurent polynomial in the variables $x_s,\dots,x_t,y_{s'},\dots,y_{t'}$, because every linear factor in the denominator is canceled by an appropriate factor in the numerator.

\end{definition}

\medskip

\begin{proposition}
\label{prop:doesn't matter}

The right-hand side of \eqref{eqn:def eh} does not depend on the chosen total order compatible with $H$.

\end{proposition}

\medskip

\begin{proof} Any two total orders compatible with $H$ can be related by a sequence of transpositions compatible with $H$, i.e. we switch the order of $c$ and $c'$ which are next to each other in the total order and not connected by an edge in $H$. The fact that the right-hand side of \eqref{eqn:def eh} is unaffected by this transposition is a consequence of
$$
\left(z_c - z_{c'} q^{-d_{\iota(c)\iota(c')}} \right) e_{\iota(c')}(z_{c'}) e_{\iota(c)}(z_c) = - \left(z_{c'} - z_{c} q^{-d_{\iota(c')\iota(c)}} \right) e_{\iota(c)}(z_{c}) e_{\iota(c')}(z_{c'})
$$
which is precisely \eqref{eqn:rel quad 1}. 

\end{proof}

\medskip

\subsection{} For $H \subset Z_e$ as above, let us compute
$$
E_H(x_s,\dots,x_t, y_{s'}, \dots, y_{t'}) = \tUpsilon^+(e_H(x_s,\dots,x_t, y_{s'}, \dots, y_{t'}))
$$
as elements of $\CS^+[\![ x^{\pm 1}_s,\dots,x^{\pm 1}_t, y^{\pm 1}_{s'}, \dots, y^{\pm 1}_{t'} ]\!]$. We will consider the formal series
$$
\delta(u) = \sum_{k \in \BZ} u^k
$$
known as the \textbf{delta function}, and note that it has the property that
\begin{equation}
\label{eqn:delta identity}
\delta\left(\frac uv\right)P(u) = \delta\left(\frac uv\right)P(v)
\end{equation}
for any Laurent polynomial $P(u)$. Because of this, we might say that the delta function ``matches" the variables $u$ and $v$. 

\medskip

\begin{proposition}
\label{prop:compute eh}

For any $H$ as above, we have
\begin{equation}
\label{eqn:form eh}
E_H(x_s,\dots,x_t,y_{s'},\dots,y_{t'}) = \eSym \left[ \frac {\prod_{a=1}^{\frac {t-s}2+1} \delta \left(\frac {x_{s+2a-2}}{z_{ia}} \right) \prod_{b=1}^{\frac {t'-s'}2+1} \delta \left(\frac {y_{s'+2b-2}}{z_{jb}} \right)}{\prod_{\{c' \rightarrow c\} \in H} (\barz_{c'} - \barz_c)} \right.
\end{equation}
$$
\left. \prod_{a, b} \frac {(z_{ia} - z_{jb} q^{-d_{ij}})(z_{jb} - z_{ia} q^{-d_{ij}})}{z_{ia}-z_{jb}} \cdot \frac {\prod_{a \neq a'} (z_{ia} - z_{ia'}q^{-d_{ii}}) }{\prod_{a < a'} (z_{ia} - z_{ia'})} \cdot  \frac {\prod_{b \neq b'} (z_{jb} - z_{jb'}q^{-d_{jj}})}{\prod_{b < b'} (z_{jb} - z_{jb'})} \right] 
$$
On the second line, the indices $a,a'$ (resp. $b,b'$) go from 1 to $\frac {t-s}2$ (resp. $\frac {t'-s'}2$). 

\medskip

\end{proposition}

\noindent In order for the RHS of formula \eqref{eqn:form eh} to be an appropriately defined formal series with values in $\CS^+$, we would need to replace in the denominator of the first line
\begin{equation}
\label{eqn:square}
\barz_c \quad \text{by} \quad z_{\iota(c)\square(c)} q^{c}
\end{equation}
for every $c \in \{s,\dots,t\}_i \sqcup \{s',\dots,t'\}_j$, where the natural number $\square(c)$ is chosen such that the delta functions on the first line of \eqref{eqn:form eh} ``match" $z_{\iota(c)\square(c)}$ with $z_c$. This also ensures that the linear factors $\barz_{c'} - \barz_c$ on the first line are canceled by various linear factors in the numerator of the second line, as expected from \eqref{eqn:def eh}. Having made this clarification, we henceforth keep the notation $\barz_c$ for better legibility.

\medskip

\begin{proof} Since $\tUpsilon^+(e_i(x)) = \delta\left(\frac x{z_{i1}}\right)$, it is an immediate consequence of \eqref{eqn:shuf prod} that
\begin{multline*}
\tUpsilon^+ \left( \mathop{\prod_{c \in Z_v \text{ in}}}_{\text{descending order}} e_{\iota(c)}(z_c) \right) = \\ = \Sym \left[ \prod_{c < c' \in Z_v} \zeta_{\iota(c')\iota(c)} \left( \frac {z_{\iota(c')\square(c')}}{z_{\iota(c)\square(c)}} \right) \prod_{a=1}^{\frac {t-s}2+1} \delta \left(\frac {x_{s+2a-2}}{z_{ia}} \right) \prod_{b=1}^{\frac {t'-s'}2+1} \delta \left(\frac {y_{s'+2b-2}}{z_{jb}} \right) \right]
\end{multline*}
where $\square(c) \in \BN$ is defined in \eqref{eqn:square}. By the same token, we have
$$
\tUpsilon^+(e_H) = \Sym \left[ \frac {\prod_{c < c' \in Z_v} \left[ (-1)^{\delta_{\iota(c)i} \delta_{\iota(c')j}+ \delta_{\iota(c)\iota(c')} \nu_{c'-c}} (z_{\iota(c)\square(c)} - z_{\iota(c')\square(c')} q^{-d_{\iota(c)\iota(c')}}) \right]}{\prod_{\{c' \rightarrow c\} \in H} \left(z_{\iota(c')\square(c')} q^{c'} - z_{\iota(c)\square(c)} q^c \right)} \right.
$$
$$
\left. \prod_{c < c' \in Z_v} \frac {z_{\iota(c')\square(c')} - z_{\iota(c)\square(c)} q^{-d_{\iota(c)\iota(c')}}}{z_{\iota(c')\square(c')} - z_{\iota(c)\square(c)}} \prod_{a=1}^{\frac {t-s}2+1} \delta \left(\frac {x_{s+2a-2}}{z_{ia}} \right) \prod_{b=1}^{\frac {t'-s'}2+1} \delta \left(\frac {y_{s'+2b-2}}{z_{jb}} \right) \right]
$$
where we used \eqref{eqn:delta identity} to match a Laurent polynomial in the variables $z_c$ to the same Laurent polynomial in the variables $z_{\iota(c)\square(c)}$. The formula above matches \eqref{eqn:form eh}, according to the paragraph immediately under the statement of Proposition \ref{prop:compute eh}. 

\end{proof}

\medskip

\subsection{} 
\label{sub:selection}

We will now define a specific linear combination of the elements $E_H$ (as $H$ runs over the subsets of $Z_e$ without oriented cycles) which is 0, and thus the corresponding linear combination of elements $e_H$ will lie in the kernel of $\tUpsilon^+$. Recall that the distinguished zig-zag $Z$ consists of $m$ southwest pointing arrows (below, we indicate the source and target of the arrows in terms of the associated variables)
$$
\Big\{ x_{t-2\alpha} \rightarrow y_{t-2\alpha +d_{ij}} \Big\}_{0 \leq \alpha < m}
$$
and $m$ northwest pointing arrows
$$
\Big\{ y_{t'-2\alpha} \rightarrow x_{t'-2\alpha +d_{ij}} \Big\}_{0 \leq \alpha < m}
$$
It is thus natural to cover the zig-zag $Z$ by the $m$ overlapping trapezoids
\begin{equation}
\label{eqn:trapezoid}
\begin{tikzcd}
& x_{t - 2\alpha -2k} \arrow[r] & \cdots \arrow[r] & x_{t-2\alpha} \arrow[dlll] \\
y_{t' - 2\alpha -2l} \arrow[r] & \cdots \arrow[r] & \cdots \arrow[r] & \cdots \arrow[r] & y_{t' - 2\alpha} \arrow[ulll]
\end{tikzcd}
\end{equation}
as $\alpha$ goes from 0 to $m-1$ (recall that $t'-t = d_{ij}+2l = -d_{ij}-2k$ as a consequence of \eqref{eqn:distinguished 1}--\eqref{eqn:distinguished 2}, so the diagonal arrows above shift indices by $d_{ij}$). The four edges of a trapezoid will be marked by the symbols $\swarrow$, $\downarrow$, $\nwarrow$, $\uparrow$, where $\uparrow$ (respectively $\downarrow$) refers to the composition of all the top (respectively bottom) arrows in the trapezoid. A \textbf{selection} is a choice
$$
S = \{\e_0,\dots,\e_{m-1}\}
$$
of an edge $\e_\alpha \in \{\swarrow, \downarrow, \nwarrow, \uparrow\}$ of the $\alpha$-th trapezoid, for all $\alpha \in \{0,\dots,m-1\}$, subject to the following conditions for all $\alpha \in \{0,\dots,m-2\}$:

\medskip

\begin{itemize}[leftmargin=*]

\item If $\e_{\alpha}$ is $\swarrow$ or $\uparrow$, then $\e_{\alpha+1}$ can only be $\swarrow$ or $\downarrow$ 

\medskip

\item If $\e_{\alpha}$ is $\nwarrow$ or $\downarrow$, then $\e_{\alpha+1}$ can only be $\nwarrow$ or $\uparrow$ 

\medskip

\end{itemize}

\noindent The \textbf{sign} of a selection $S$ is $(-1)^{\sigma(S)}$, where
\begin{equation}
\label{eqn:sign}
\sigma(S) = \Big| \Big\{ 1 \leq \alpha < m \text{ s.t. } \e_{\alpha} \text{ is } \nwarrow \text{ or } \uparrow \Big\} \Big|
\end{equation}
and for the chosen edge $\e_{\alpha}$ in the $\alpha$-th trapezoid, we write
\begin{equation}
\label{eqn:def delta}
\delta(\e_\alpha) = \begin{cases} \barx_{t-2\alpha} - \bary_{t'-2\alpha - 2l} &\text{if } \e_{\alpha} \text{ is } \swarrow  \\ \bary_{t'-2\alpha - 2l}  - \bary_{t'-2\alpha}  &\text{if } \e_{\alpha} \text{ is } \downarrow  \\ \bary_{t'-2\alpha} - \barx_{t - 2\alpha - 2k} &\text{if } \e_{\alpha} \text{ is } \nwarrow  \\ \barx_{t - 2\alpha - 2k} - \barx_{t-2\alpha} &\text{if } \e_{\alpha} \text{ is } \uparrow  
\end{cases}
\end{equation}
(the variables in question are as in \eqref{eqn:trapezoid}, see also \eqref{eqn:notation 4} for the bar notation). 

\medskip

\begin{proposition}
\label{prop:selection combi}

With the terminology above, we have the following combinatorial identity
\begin{equation}
\label{eqn:selection combi}
\sum_{S \text{ selection}} (-1)^{\sigma(S)} \prod_{\alpha = 0}^{m-1} \delta(\e_{\alpha}) = 0
\end{equation}
for any distinguished zig-zag $Z$.

\end{proposition}

\medskip

\begin{proof} We will prove the required statement by induction on $m$ (the base case $m=1$ is trivial, as all the signs are 1 and all products in \eqref{eqn:selection combi} consist of a single linear factor). For any selection 
$$
S = \{\e_0, \dots, \e_{m-2}, \e_{m-1}\}
$$
let $S' = \{\e_0, \dots, \e_{m-2}\}$. Once we have chosen $S'$, then the two possible choices for $\e_{m-1}$ yield the following contribution to $\delta(\e_{m-1})$ in \eqref{eqn:selection combi}:

\medskip

\begin{itemize}[leftmargin=*]

\item $\barx_{t-2(m-1)} - \bary_{t'-2(l+m-1)} + \bary_{t'-2(l+m-1)} - \bary_{t'-2(m-1)}$ if $\e_{m-2}$ is $\swarrow$ or $\uparrow$

\medskip

\item $\barx_{t - 2(k+m-1)} - \barx_{t- 2(m-1)} + \bary_{t'-2(m-1)} - \barx_{t - 2(k+m-1)}$ if $\e_{m-2}$ is $\nwarrow$ or $\downarrow$

\end{itemize}

\medskip 

\noindent In both cases, the quantity above is $\pm (\barx_{t-2(m-1)} - \bary_{t'-2(m-1)})$. The sign is precisely accounted for by the discrepancy between $\sigma(S)$ and $\sigma(S')$, so we conclude that 
$$
\Big( \text{LHS of \eqref{eqn:selection combi} for } m \Big) = \Big( \text{LHS of \eqref{eqn:selection combi} for } m-1 \Big) \cdot \Big( \barx_{t-2(m-1)} - \bary_{t'-2(m-1)} \Big)
$$
As the RHS is zero by the induction hypothesis, so is the LHS. This completes the proof of the induction step. 

\end{proof}

\medskip

\subsection{} 
\label{sub:refined selection}

A small issue with trapezoids, as defined above, is that their $\uparrow$ and $\downarrow$ edges typically consist of more that one arrow in the oriented graph $Z$. Thus, let us define a \textbf{refined selection} to be a multiset
$$
\mathscr{S} = \{\e_0,\dots,\e_{m-1}\}
$$
where $\e_{\alpha} \in Z_e$ lies in the $\alpha$-th trapezoid, for all $\alpha \in \{0,\dots,m-1\}$. Intuitively, if the $\alpha$-th trapezoid is of the form \eqref{eqn:trapezoid}, there are $k+l+2$ choices for $\e_{\alpha}$: $k$ choices from the top, $l$ choices from the bottom, and $2$ choices from the diagonals. Rigorously, a refined selection involves the choice of one of the following types of edges from the $\alpha$-th trapezoid, for each $\alpha \in \{0,\dots,m-1\}$
$$
\e_{\alpha} = \begin{cases} x_{t-2\alpha} \rightarrow y_{t'-2\alpha - 2l} &\text{or} \\
y_{t'-2\alpha - 2v } \rightarrow y_{t'-2\alpha -2(v-1)} &\text{or} \\
y_{t'-2\alpha} \rightarrow x_{t-2\alpha-2k} &\text{or} \\
x_{t-2\alpha -2u} \rightarrow x_{t-2\alpha -2(u-1)} \end{cases}
$$
for some $u \in \{1,\dots,k\}$ and $v \in \{1,\dots,l\}$. For each type of edge above, we define $\delta(\e_{\alpha})$ as the difference of the variable on the left minus the variable on the right, with bars on top (see \eqref{eqn:def delta}). Also define $\sigma(\mathscr{S})$ by the same formula \eqref{eqn:sign} as $\sigma(S)$. It is easy to observe that \eqref{eqn:selection combi} implies 
\begin{equation}
\label{eqn:refined selection combi}
\sum_{\mathscr{S} \text{ refined selection}} (-1)^{\sigma(\mathscr{S})} \prod_{\alpha = 0}^{m-1} \delta(\e_{\alpha}) = 0
\end{equation}
simply because $\sum_{u=1}^k (\barx_{t-2\alpha  - 2u} - \barx_{t-2\alpha - 2(u-1)}) = \barx_{t-2\alpha -2k} - \barx_{t-2\alpha}$ (as well as the analogous formula for $\bary$'s instead of $\barx$'s). The following result crucially uses the defining properties of a selection contained in the two bullets in Subsection \ref{sub:selection}.

\medskip

\begin{proposition}
\label{prop:complement refined selection}

Let $H \subset Z_e$ be the complement of a refined selection. Then $H$ does not contain any oriented cycles.

\end{proposition}

\medskip

\begin{proof} Assume the contrary, that the complement $H$ of a refined selection $\mathscr{S}$ contained an oriented cycle. Because of the specific shape of the oriented graph $Z$, this can only happen if there are two intersecting diagonal edges in $H$, together with all the horizontal edges connecting the endpoints of the diagonal ones. In other words, the complement $H$ of $\mathscr{S}$ must contain a cycle of the form
\begin{multline}
x_u \rightarrow y_{u'-2l'} \rightarrow y_{u'-2(l'-1)} \rightarrow \dots y_{u'-2} \rightarrow y_{u'} \rightarrow \\ \rightarrow x_{u-2k'} \rightarrow x_{u-2(k'-1)} \rightarrow \dots \rightarrow x_{u-2} \rightarrow x_u \label{eqn:cycle}
\end{multline}
for some non-negative integers $k', l'$ such that $k'+l' = -d_{ij}$, and some $u = t-2\alpha$, $u' = t'-2\alpha'$ with $\alpha,\alpha' \in \{0,\dots,m-1\}$ satisfying
$$
u'-2l' = u +d_{ij} \quad \Leftrightarrow \quad u-2k' = u' + d_{ij}
$$
The assumptions above, together with $t'-2l = t+d_{ij} \Leftrightarrow t-2k = t'+d_{ij}$, yield
$$
\alpha \geq \alpha' \Leftrightarrow k \geq k' \Leftrightarrow l \leq l'
$$
We have three cases to analyze.

\medskip

\noindent \emph{Case 1:} If $\alpha = \alpha'$, then $\e_{\alpha}$ cannot be $\swarrow$ or $\nwarrow$, due to the presence of these edges in the cycle \eqref{eqn:cycle}. However, if $\e_{\alpha}$ were $\uparrow$ or $\downarrow$, the fact that $k = k'$ and $l=l'$ implies that one of the edges
$$
x_{t-2\alpha -2k} \rightarrow x_{t-2\alpha -2(k-1)} \rightarrow \dots \rightarrow x_{t-2\alpha -2} \rightarrow x_{t-2\alpha }
$$
$$
y_{t'-2\alpha -2l} \rightarrow y_{t'-2\alpha -2(l-1)} \rightarrow \dots \rightarrow y_{t'-2\alpha-2} \rightarrow y_{t'-2\alpha}
$$
lies in $\mathscr{S}$. Then the edge in question does not lie in $H$, which contradicts the presence of the cycle \eqref{eqn:cycle} in $H$.

\medskip

\noindent \emph{Case 2:} If $\alpha > \alpha'$, then $l < l'$ and for the same reason as in Case 1
\begin{equation}
\label{eqn:prohibition 1}
\e_{\alpha'}, \dots, \e_{\alpha} \text{ cannot be } \downarrow
\end{equation}
But $\e_{\alpha'}$ also cannot be $\nwarrow$, which leaves only two choices: $\e_{\alpha'}$ must be $\swarrow$ or $\uparrow$. In either of these two cases, the rules of selections imply that
$$
(\e_{\alpha'+1},\dots, \e_{\alpha}) = (\swarrow, \dots, \swarrow)
$$
The fact that $\e_{\alpha}$ must be $\swarrow$ contradicts the existence of the cycle \eqref{eqn:cycle} in the complement of the selection.

\medskip

\noindent \emph{Case 3:} If $\alpha < \alpha'$, then $k < k'$ and for the same reason as in Case 1
\begin{equation}
\label{eqn:prohibition 2}
\e_{\alpha}, \dots, \e_{\alpha'} \text{ cannot be } \uparrow
\end{equation}
But $\e_{\alpha}$ also cannot be $\swarrow$, which leaves only two choices: $\e_{\alpha}$ must be $\nwarrow$ or $\downarrow$. In either of these two cases, the rules of selections imply that
$$
(\e_{\alpha+1},\dots, \e_{\alpha'}) = (\nwarrow, \dots, \nwarrow)
$$
The fact that $\e_{\alpha'}$ must be $\nwarrow$ contradicts the existence of the cycle \eqref{eqn:cycle} in the complement of the selection. 

\end{proof}

\subsection{}

If we divide \eqref{eqn:refined selection combi} by
$$
\prod_{\{c' \rightarrow c\} \in Z_e} (\barz_{c'} - \barz_c)
$$
then together with \eqref{eqn:form eh} (which we may invoke due to Proposition \ref{prop:complement refined selection}) we have
\begin{equation}
\label{eqn:equality refined selection 1}
\sum_{\mathscr{S} \text{ refined selection}} (-1)^{\sigma(\mathscr{S})} E_{Z_e \backslash \mathscr{S}} = 0
\end{equation}
However, we owe the reader a clarification about the notation $E_{Z_e \backslash \mathscr{S}}$. If $\mathscr{S}$ were simply a set of edges of $Z_e$, then the complement $H = Z_e \backslash \mathscr{S}$ has no oriented cycles by Proposition \ref{prop:complement refined selection}, and we may use formula \eqref{eqn:form eh}. However, in general $\mathscr{S}$ is a multiset of edges of $Z_e$, because a given horizontal edge in Figure 7 may appear more than once in $\mathscr{S}$. If we denote the multiplicity of an edge in the multiset $\mathscr{S}$ by $\mu_{\mathscr{S}}(c'\rightarrow c) \in \{0,1,\dots\}$, then we define the following analogue of \eqref{eqn:form eh}
$$
E_{Z_e \backslash \mathscr{S}}(x_s,\dots,x_t,y_{s'},\dots,y_{t'}) = \Sym \left[ \frac {\prod_{a=1}^{\frac {t-s}2+1} \delta \left(\frac {x_{s+2a-2}}{z_{ia}} \right) \prod_{b=1}^{\frac {t'-s'}2+1} \delta \left(\frac {y_{s'+2b-2}}{z_{jb}} \right)}{\prod_{\{c' \rightarrow c\} \in Z_e} (\barz_{c'} - \barz_c)^{1-\mu_{\mathscr{S}}(c' \rightarrow c)}} \right.
$$
$$
\left. \prod_{a, b} \frac {(z_{ia} - z_{jb} q^{-d_{ij}})(z_{jb} - z_{ia} q^{-d_{ij}})}{z_{ia}-z_{jb}} \cdot \frac {\prod_{a \neq a'} (z_{ia} - z_{ia'}q^{-d_{ii}}) }{\prod_{a < a'} (z_{ia} - z_{ia'})} \cdot  \frac {\prod_{b \neq b'} (z_{jb} - z_{jb'}q^{-d_{jj}})}{\prod_{b < b'} (z_{jb} - z_{jb'})} \right] 
$$
as an element of $\CS^+[\![ x^{\pm 1}_s,\dots,x^{\pm 1}_t, y^{\pm 1}_{s'}, \dots, y^{\pm 1}_{t'} ]\!]$. Similarly, we write
\begin{equation}
\label{eqn:def eh refined}
e_{Z_e\backslash \mathscr{S}}\left(x_s, \dots, x_t, y_{s'}, \dots, y_{t'} \right) = \frac 1{\prod_{\{c' \rightarrow c\} \in Z_e} \left(\barz_{c'} - \barz_c \right)^{1-\mu_{\mathscr{S}(c'\rightarrow c)}}}
\end{equation}
$$
\prod_{c < c' \in Z_v} \Big[ (-1)^{\delta_{\iota(c)i} \delta_{\iota(c')j}+ \delta_{\iota(c)\iota(c')} \nu_{c'-c}} (z_c - z_{c'} q^{-d_{\iota(c)\iota(c')}}) \Big]  \mathop{\prod_{c \in Z_v \text{ in}}}_{\text{descending order}} e_{\iota(c)}(z_c) 
$$
for any refined selection $\mathscr{S}$ (above, $c<c'$ denotes any total order on $Z_v$ compatible with $Z_e\backslash \mathscr{S}$, defined just like the analogous notion in \eqref{eqn:def eh}).

\medskip

\begin{proposition}
\label{prop:equality refined selection}

For any distinguished zig-zag $Z$, the formal series
\begin{equation}
\label{eqn:equality refined selection 2}
\rho_Z := \sum_{\mathscr{S} \text{ refined selection}} (-1)^{\sigma(\mathscr{S})} e_{Z_e \backslash \mathscr{S}}
\end{equation}
lies in $K^+[\![ x^{\pm 1}_s,\dots,x^{\pm 1}_t, y^{\pm 1}_{s'}, \dots, y^{\pm 1}_{t'} ]\!]$, where $K^+ = \emph{Ker }\tUpsilon^+$ from \eqref{eqn:kernel}.

\end{proposition}

\medskip

\noindent The Proposition is an immediate consequence of \eqref{eqn:equality refined selection 1}, since by a straightforward analogue (in fact, a slight generalization) of Proposition \ref{prop:compute eh}, we have
$$
\tUpsilon^+(e_{Z_e\backslash \mathscr{S}}) = E_{Z_e\backslash \mathscr{S}}
$$
For any Laurent polynomial $\tau(x_s,\dots,x_t,y_{s'},\dots,y_{t'})$, we define
\begin{align*}
&e_{Z_e\backslash \mathscr{S},\tau} = \Big[\tau(x_s,\dots,x_t,y_{s'},\dots,y_{t'}) e_{Z_e\backslash \mathscr{S}}(x_s,\dots,x_t,y_{s'},\dots,y_{t'}) \Big]_{\text{constant term}} \\
&\rho_{Z,\tau} = \Big[\tau(x_s,\dots,x_t,y_{s'},\dots,y_{t'}) \rho_Z(x_s,\dots,x_t,y_{s'},\dots,y_{t'}) \Big]_{\text{constant term}}
\end{align*}
and Proposition \ref{prop:equality refined selection} implies that $\rho_{Z,\tau} \in K^+$.

\medskip

\subsection{} 

Our interest in the series $\rho_Z$ is motivated by the following result.

\medskip

\begin{proposition}
\label{prop:realize}

An element $R \in \CV^-$ lies in $\CS^-$ if and only if
\begin{equation}
\label{eqn:realize pairing}
\left \langle \tUUp \cdot \rho_{Z,\tau_\ell}(x_s,\dots,x_t,y_{s'},\dots,y_{t'}) \cdot \tUUp, R \right \rangle_{\tU\CV} = 0
\end{equation}
for all distinguished zig-zags $Z$, where $\{\tau_\ell\}_{\ell \in \BZ}$ denotes an arbitrary collection of Laurent polynomials of homogeneous degree $\ell$ such that
\begin{equation}
\label{eqn:tau neq}
\tau_\ell \left(q^{\frac {t-s}2}, \dots , q^{\frac {s-t}2}, q^{\frac {t'-s'}2}, \dots , q^{\frac {s'-t'}2} \right) \neq 0
\end{equation}

\end{proposition}

\medskip

\begin{proof} We will use the following statement as a stepping stone.

\medskip

\begin{claim}
\label{claim:stepping stone} 

Consider a distinguished zig-zag $Z$ as in Figure 7, and an element
$$
R \in \CV_{-\bs^i\left( \frac {t-s}2 + 1 \right)-\bs^j\left( \frac {t'-s'}2 + 1 \right)}
$$
which satisfies \eqref{eqn:wheel} for all distinguished zig-zags $Z' \subsetneq Z$. Then $R$ satisfies \eqref{eqn:wheel} for the zig-zag $Z$ if and only if
\begin{equation}
\label{eqn:stepping stone}
\Big \langle \rho_{Z,\tau_\ell}(x_s,\dots,x_t,y_{s'},\dots,y_{t'}), R \Big \rangle_{\tU\CV} = 0
\end{equation}
where $\{\tau_{\ell}\}_{\ell \in \BZ}$ denotes an arbitrary collection of Laurent polynomials of homogeneous degree $\ell$ that satisfy property \eqref{eqn:tau neq}.
 
\end{claim}

\medskip

\noindent Let us show how Claim \ref{claim:stepping stone} allows us to prove Proposition \ref{prop:realize}, which we will do by induction on the number of vertices of the zig-zag $Z$ (the base case is an empty zig-zag). The ``only if" statement of \eqref{eqn:realize pairing} is an immediate consequence of \eqref{eqn:descend} and Proposition \ref{prop:equality refined selection}. For the ``if" statement, consider a distinguished zig-zag $Z$ as in Figure 7, and let us fix a partition of the variables of $R$ into two groups

\medskip

\begin{itemize}[leftmargin=*]

\item The big variables $u_1, \dots, u_p$ of some colors $i_1,\dots, i_p$

\medskip

\item The small variables $x_s,\dots,x_t,y_{s'},\dots,y_{t'}$ of colors $i,\dots,i,j,\dots,j$

\end{itemize}

\medskip

\noindent With respect to this choice of variables, we may expand
\begin{equation}
\label{eqn:partial expand}
R = \sum_{k_1,\dots,k_p \in \BZ} \frac {u_1^{-k_1} \dots u_p^{-k_p}  R_{k_1,\dots,k_p}(x_s,\dots,x_t,y_{s'},\dots,y_{t'})}{\prod_{a < b} \left(1-\frac {u_b}{u_a} \right) \prod_{b,c} \left(1-\frac {z_c}{u_b}\right)}
\end{equation}
(the products in the denominator run over various pairs of indices: $a$ and $b$ over $\{1,\dots,p\}$, and $c$ over $\{s,\dots,t\}_i \sqcup \{s', \dots, t'\}_j$; $z_c$ are defined as in \eqref{eqn:notation 1}). It is very important that only finitely many of the $R_{k_1,\dots,k_p}$ are non-zero, which is a consequence of the fact that any $R \in \CV^-$ is a Laurent polynomial divided by a specific collection of linear factors. Then, formula \eqref{eqn:pair formula 1} takes the following form
\begin{multline}
\Big \langle e_{i_1,k'_1} \dots e_{i_p,k'_p} \cdot \rho_{Z,\tau}, R  \Big \rangle_{\tU\CV} = \int_{|u_1|\gg \dots \gg |u_p| \gg |x_a|,|y_b|}  \\ \sum_{k_1, \dots, k_p \in \BZ} \frac {q^\bullet u_1^{k'_1-k_1} \dots u_p^{k_p'-k_p}}{\prod_{a < b} \left(1 - \frac {u_{b} q^\bullet}{u_a} \right)} \cdot \left \langle \rho_{Z,\tau}, \frac {R_{k_1,\dots, k_p}}{\prod_{b,c} \left(1 - \frac {z_cq^{\bullet}}{u_b} \right)}  \right \rangle_{\tU\CV} \label{eqn:really big}
\end{multline}
where the denominator goes over all allowable indices (in particular, $c$ runs over $\{s,\dots,t\}_i \sqcup \{s', \dots, t'\}_j$ and $z_c$ is given by \eqref{eqn:notation 1}) and the various $\bullet$'s are placeholders for various integers that will not be important in the subsequent argument. To show that $R$ satisfies \eqref{eqn:wheel} for the given zig-zag $Z$, we must show that
\begin{equation}
\label{eqn:divide particular}
(x-y)^{m_Z} \ \Big| \ R_{k_1,\dots,k_p}(x q^s,\dots,x q^t,y q^{s'},\dots,yq^{t'})
\end{equation}
for all $k_1,\dots,k_p$. We will do so by induction on the finite set 
$$
\Big\{ (k_1,\dots,k_p) \text{ s.t. } R_{k_1,\dots,k_p} \neq 0 \Big\}
$$
in lexicographic order (i.e. first compare $k_1$, then compare $k_2$ to break ties etc). Condition \eqref{eqn:realize pairing} implies that the left-hand side of \eqref{eqn:really big} is zero for all $k_1',\dots,k_p'$. Therefore, the same is true for the right-hand side, and we obtain 
\begin{equation}
\label{eqn:really small}
0 = \Big \langle \rho_{Z,\tau}, R_{k'_1,\dots,k_p'} \Big \rangle_{\tU\CV} + \underline{\mathop{\sum_{(k_1,\dots,k_p) < (k_1',\dots,k_p')}}_{\text{lexicographically}} \Big \langle \rho_{Z,\tau}, \flat \cdot R_{k_1,\dots,k_p} \Big \rangle_{\tU\CV}}
\end{equation}
where $\flat$ stands for symmetric Laurent polynomials in $x_s,\dots,x_t,y_{s'},\dots,y_{t'}$ that arise from the power series expansion of the factors $1 - \frac {z_cq^\bullet}{u_b}$ in \eqref{eqn:really big} (multiplication by $\flat$ preserves the divisibility conditions \eqref{eqn:wheel}). The induction hypothesis implies that \eqref{eqn:divide particular} holds for all terms in the underlined sum in the right-hand side of \eqref{eqn:really small}; thus Claim \ref{claim:stepping stone} implies that the underlined sum is 0. Therefore, so is the first pairing in \eqref{eqn:really small}, and invoking Claim \ref{claim:stepping stone} again implies that 
$$
(x-y)^{m_Z} \ \Big| \ R_{k_1',\dots,k_p'} (x q^s,\dots,x q^t,y q^{s'},\dots,yq^{t'})
$$
The induction step is thus complete. 

\begin{remark}
\label{rem:ideal}

The argument above shows that \eqref{eqn:realize pairing} would still remain valid if one replaced the two-sided ideal generated by the $\rho_{Z,\tau_\ell}$ by the corresponding left ideal (and by an analogous proof, one could use instead the corresponding right ideal).

\end{remark}

\medskip

\noindent It remains to prove Claim \ref{claim:stepping stone}. The ``only if" statement follows from \eqref{eqn:descend} and Proposition \ref{prop:equality refined selection}. For the ``if" statement, note that \eqref{eqn:pair formula 1} and \eqref{eqn:def eh refined} imply the following formula for any refined selection $\mathscr{S}$ and any Laurent polynomial $\tau$
\begin{equation}
\label{eqn:pair e with r}
\Big \langle e_{Z_e \backslash \mathscr{S}, \tau}, R \Big \rangle_{\tU\CV} =\mathop{\int_{|z_{c'}| \gg |z_c|}}_{\text{if }\{c'\rightarrow c\} \notin \mathscr{S}} \frac {\omega_{\tau}(x_s,\dots,x_t,y_{s'},\dots,y_{t'})}{\prod_{\{c' \rightarrow c\} \in Z_e} (\barz_{c'} - \barz_c)^{1-\mu_{\mathscr{S}}(c' \rightarrow c)}}  \prod_c Dz_c
\end{equation}
where
$$
\omega_{\tau} = (\tau \cdot r)(x_s,\dots,x_t,y_{s'},\dots,y_{t'})  \prod_{a < a' \in \{s,\dots,t\}_i} (x_{a'}-x_{a})\prod_{b < b' \in \{s',\dots,t'\}_j} (y_{b'}-y_{b})
$$
is a Laurent polynomial. \footnote{The connection between the rational function $R$ and the Laurent polynomial $r$ is 
\begin{equation}
\label{eqn:special r}
R(x_s,\dots,x_t,y_{s'},\dots,y_{t'}) = \frac {r(x_s,\dots,x_t,y_{s'},\dots,y_{t'})}{\prod_{a \in \{s,\dots,t\}_i, b \in \{s',\dots,t'\}_j} (x_a - y_b)} 
\end{equation}
see \eqref{eqn:symmetric rational function}. In all formulas involving $R$ and $r$, we plug $x_s,\dots,x_t$ (respectively $y_{s'}, \dots, y_{t'}$) into the variables of color $i$ (respectively $j$) of $R$ and $r$.} Running the argument in the proof of Proposition \ref{prop:descends} (with the terminology therein), we deduce the following analogue of \eqref{eqn:end induction 2}
\begin{equation}
\label{eqn:slight imprecision}
\Big \langle e_{Z_e \backslash \mathscr{S}, \tau}, R \Big \rangle_{\tU\CV} ``=" \sum^{\text{fixed fair partition}}_{\{1,\dots,n\}=\bar{A}_1 \sqcup \dots \sqcup \bar{A}_t} \int_{|w_1| = \dots = |w_t|}   \prod_{s=1}^t Dw_s 
\end{equation}
$$
\left[ \underset{(z_{c_s^{(1)}},\dots,z_{c_s^{(n_s)}}) = (w_s q^{n_s-1}, \dots, w_sq^{1-n_s})}{\text{Res}} \Sym \left(\frac {\omega_{\tau}(x_s,\dots,x_t,y_{s'},\dots,y_{t'})}{\prod_{\{c' \rightarrow c\} \in Z_e} (\barz_{c'} - \barz_c)^{1-\mu_{\mathscr{S}}(c' \rightarrow c)}} \right) \right]_{\forall s \in \{1,\dots,t\}}
$$
where each set $\bar{A}_s = \{z_{c_s^{(1)}}, \dots, z_{c_s^{(n_s)}}\}$ consists either of all $x$'s or all $y$'s. 

\medskip

\noindent The ``=" symbol in \eqref{eqn:slight imprecision} is due to the following slight imprecision, that we will now correct. In proving \eqref{eqn:end induction 2}, we used Claim \ref{claim:induction}, which in turn used the fact that $R$ satisfies \eqref{eqn:wheel} for the distinguished zig-zag $Z$ and all of its sub zig-zags. However, in the present situation, our assumption is that $R$ satisfies \eqref{eqn:wheel} for all proper sub zig-zags of $Z$, but we do not have this property for $Z$ itself. Therefore, the discrepancy between the LHS and RHS of \eqref{eqn:slight imprecision} is due to the unique residue where property \eqref{eqn:wheel} is invoked for $Z$ itself. This correesponds to $t = 2$ and
\begin{equation}
\label{eqn:troublesome residue}
\bar{A}_1 = \{x_s,\dots,x_t\} \quad \text{and} \quad  \bar{A}_2 = \{y_{s'},\dots,y_{t'}\}
\end{equation}
However, due to the denominators $(\barz_{c'}-\barz_c)^{1-\mu_{\mathscr{S}}(c' \rightarrow c)}$ in \eqref{eqn:slight imprecision}, the appearance of the residue \eqref{eqn:troublesome residue} requires the refined selection $\mathscr{S}$ to not contain any horizontal edge in Figure 7. By the defining criteria of refined selections in Subsection \ref{sub:refined selection}, the only two refined selections which can produce the troublesome residues \eqref{eqn:troublesome residue} are
$$
\mathscr{S}_{\swarrow} = \Big\{ \swarrow, \dots, \swarrow \Big\} \qquad \text{and} \qquad \mathscr{S}_{\nwarrow} = \Big\{ \nwarrow, \dots, \nwarrow \Big\}
$$
It is easy to compute the corresponding residues, and we obtain for $\mathscr{S}_{\swarrow}$
$$
\mathop{\underset{(x_s,\dots,x_t) = (x q^{\frac {t-s}2}, \dots, xq^{\frac {s-t}2})}{\text{Res}}}_{(y_{s'},\dots,y_{t'}) = (y q^{\frac {t'-s'}2}, \dots, yq^{\frac {s'-t'}2})} \Sym \left(\frac {\omega_{\tau}(x_s,\dots,x_t,y_{s'},\dots,y_{t'})}{\prod_{\{c' \rightarrow c\} \in Z_e} (\barz_{c'} - \barz_c)^{1-\mu_{\mathscr{S}_{\swarrow}}(c' \rightarrow c)}} \right) = 
$$
\begin{equation}
\label{eqn:formula sw arrow}
= \frac {\omega_{\tau}(x q^{\frac {t-s}2},\dots,x q^{\frac {s-t}2},y q^{\frac {t'-s'}2},\dots,y q^{\frac {s'-t'}2})}{(y-x)^{m_Z}x^{\frac {t-s}2-1}y^{\frac {t'-s'}2-1} q^\#}
\end{equation}
(recall that $s+t=s'+t'$ for a distinguished zig-zag) for some integer $\#$ that will not make a difference in the subsequent argument. The case of $\nwarrow$ is given by formula \eqref{eqn:formula sw arrow} with $x,y,s,t,s',t'$ replaced by $y,x,s',t',s,t$. Therefore, in order to correct \eqref{eqn:slight imprecision}, we need to add the following term to its right-hand side that accounts for the two troublesome refined selections $\mathscr{S}_{\swarrow}$ and $\mathscr{S}_{\nwarrow}$
$$
\Big( \delta_{\mathscr{S} \mathscr{S}_{\swarrow}} + (-1)^{m_Z-1} \delta_{\mathscr{S} \mathscr{S}_{\nwarrow}} \Big) \left[ \underset{y = x}{\text{Res}} \  \frac {\omega_{\tau}(x q^{\frac {t-s}2},\dots,x q^{\frac {s-t}2},y q^{\frac {t'-s'}2},\dots,y q^{\frac {s'-t'}2})}{(y-x)^{m_Z}x^{\frac {t-s}2-1}y^{\frac {t'-s'}2-1} q^{\#}} \right]_{\text{constant term in }x}
$$
\footnote{In the formula above, the sign $(-1)^{m_Z-1}$ comes from two contributions: first we have $(-1)^{m_Z}$ because of the discrepancy between $(x-y)^{m_Z}$ and $(y-z)^{m_Z}$; secondly, we have $(-1)$ because of the discrepancy between taking the residue at $x=y$ and the residue at $y=x$.}; specifically, we encounter the residue from Case 1 in the proof of Proposition \ref{prop:descends}, where the two geometric progressions are midpoint aligned, and the shape they trace out matches the zig-zag $Z$. Taking the sum of \eqref{eqn:slight imprecision} thus corrected over all refined selections $\mathscr{S}$, weighted by the sign $(-1)^{\sigma(\mathscr{S})}$, yields 
\begin{equation}
\label{eqn:rho}
\Big \langle \rho_{Z,\tau}, R \Big \rangle_{\tU\CV} = \sum^{\text{fixed fair partition}}_{\{1,\dots,n\}=\bar{A}_1 \sqcup \dots \sqcup \bar{A}_t} \int_{|w_1| = \dots = |w_t|}  \prod_{s=1}^t Dw_s 
\end{equation}
$$
\left[ \underset{(z_{c_s^{(1)}},\dots,z_{c_s^{(n_s)}}) = (w_s q^{n_s-1}, \dots, w_sq^{1-n_s})}{\text{Res}} \Sym \left(\sum_{\mathscr{S}} (-1)^{\sigma(\mathscr{S})} \frac {\omega_{\tau}(x_s,\dots,x_t,y_{s'},\dots,y_{t'})}{\prod_{\{c' \rightarrow c\} \in Z_e} (\barz_{c'} - \barz_c)^{1-\mu_{\mathscr{S}}(c' \rightarrow c)}} \right) \right]_{\forall s \in \{1,\dots,t\}} 
$$
$$
+ 2 \left[ \underset{y = x}{\text{Res}} \ \frac {\omega_{\tau}(x q^{\frac {t-s}2},\dots,x q^{\frac {s-t}2},y q^{\frac {t'-s'}2},\dots,y q^{\frac {s'-t'}2})}{(y-x)^{m_Z}x^{\frac {t-s}2-1}y^{\frac {t'-s'}2-1} q^{\#}}\right]_{\text{constant term in }x}
$$
where the third line is due to the contribution of the troublesome refined selections $\mathscr{S}_{\swarrow}$ and $\mathscr{S}_{\nwarrow}$, which account for the difference between the LHS and RHS of \eqref{eqn:slight imprecision}. Since the second line of the expression above vanishes (this is because it is equal to the RHS of \eqref{eqn:end induction 2} with the LHS replaced by
$$
\Big \langle \text{LHS of \eqref{eqn:equality refined selection 1}}, R \Big \rangle_{\oCS\CS}
$$
and the LHS in the equation above is 0) we conclude that the left-hand side of \eqref{eqn:rho} is 0 if and only if the third line of \eqref{eqn:rho} is 0. However, because $R$ satisfies \eqref{eqn:wheel} for all proper sub zig-zags of $Z$, we already have
$$
(y-x)^{m_Z-1} \Big | r \left(x q^{\frac {t-s}2},\dots,x q^{\frac {s-t}2},y q^{\frac {t'-s'}2},\dots,y q^{\frac {s'-t'}2}\right) 
$$
Therefore, we conclude that $\langle \rho_{Z,\tau},R \rangle_{\tU\CV} = 0$ if and only if
$$
\frac {\omega_{\tau} \left(x q^{\frac {t-s}2},\dots,x q^{\frac {s-t}2},y q^{\frac {t'-s'}2},\dots,y q^{\frac {s'-t'}2}\right)}{(y-x)^{m_Z}x^{\frac {t-s}2-1}y^{\frac {t'-s'}2-1} q^{\#}}
$$
is a Laurent polynomial in $x$ and $y$ of homogeneous degree 0 (if the homogeneous degree were non-zero, then the constant term would vanish on general grounds). This implies that
$$
(y-x)^{m_Z} \Big | r \left(x q^{\frac {t-s}2},\dots,x q^{\frac {s-t}2},y q^{\frac {t'-s'}2},\dots,y q^{\frac {s'-t'}2}\right) 
$$
as long as $\tau$ is chosen homogeneous of degree $\frac {t-s}2-1+\frac {t'-s'}2-1+m_Z-\deg r$ and satisfies \eqref{eqn:tau neq} (if $r$ is not homogeneous, then we perform the argument above for its homogeneous components).  

\end{proof}

\subsection{}

Consider the two-sided ideal
\begin{equation}
\label{eqn:j}
J^+ = (\text{coefficients of } \rho_{Z} )_{Z \text { distinguished zig-zag}} \subset \tUUp
\end{equation}
and define
$$
\UUp = \tUUp \Big/ J^+
$$
By Proposition \ref{prop:realize}, we have $(J^+)^\perp = \CS^-$ with respect to the pairing $\langle \cdot, \cdot \rangle_{\tU \CV}$, and thus we obtain the pairing $\langle \cdot, \cdot \rangle_{U\CS}$ in diagram \eqref{eqn:desired diagram}. If all the vector spaces in said diagram were finite-dimensional (or at least finite-dimensional in every $\pm \nn \times \BZ$ degree), then the non-degeneracy of $\langle \cdot, \cdot \rangle_{\tU \CV}$ would imply that $(\CS^-)^\perp = J^+$ and yield the non-degeneracy of $\langle \cdot, \cdot \rangle_{U\CS}$. However, because all these vector spaces are in general infinite-dimensional, care must be taken when proving the following result.

\medskip

\begin{proposition}
\label{prop:final non-deg}

The pairing on the bottom of diagram \eqref{eqn:desired diagram}, namely
\begin{equation}
\label{eqn:pairing final}
\UUp \otimes \CS^- \xrightarrow {\langle \cdot, \cdot \rangle_{U\CS}} \BK
\end{equation}
is non-degenerate in both arguments. 

\end{proposition} 

\medskip

\noindent The natural analogue of Proposition \ref{prop:final non-deg} where we switch the roles of $+$ and $-$ also holds. We will prove Proposition \ref{prop:final non-deg} in the next Section using the combinatorics of words. But for now, let us deduce from it the proofs of our main Theorems.

\medskip

\begin{proof} \emph{of Theorem \ref{thm:connect}:} As explained in the paragraph preceding the statement of Proposition \ref{prop:final non-deg}, the pairings on top of either diagram \eqref{eqn:desired diagram intro} descend to pairings on the bottom. The non-degeneracy of the latter is established in Proposition \ref{prop:final non-deg}.  

\end{proof}

\begin{proof} \emph{of Theorems \ref{thm:main} and \ref{thm:iso}:} We showed in \eqref{eqn:equality refined selection 2} that
$$
J^+ \subseteq K^+
$$
where $K^+ = \text{Ker }\tUpsilon^+$. Thus, we obtain a surjective algebra homomorphism
$$
\UUp = \tUUp/J^+ \xrightarrow{\Upsilon^+} \tUUp/K^+ \stackrel{\text{Prop.  \ref{prop:coincide}}}= \CS^+
$$
which intertwines the pairings \eqref{eqn:pairing final} and \eqref{eqn:pairing shuffle}, as explained in the paragraph after the statement of Proposition \ref{prop:descends}. As such, any element in the kernel of $\Upsilon^+$ would also have to be in the radical of the pairing \eqref{eqn:pairing final}. As the latter is non-degenerate due to Proposition \ref{prop:final non-deg}, we conclude that $\Upsilon^+$ is an isomorphism, i.e.
\begin{equation}
\label{eqn:jk}
J^+ = K^+
\end{equation}
This implies Theorem \ref{thm:iso}. The compatibility of the pairings in diagram \eqref{eqn:desired diagram} and the non-degeneracy of \eqref{eqn:pairing final} imply that
\begin{equation}
\label{eqn:ij}
I^+ = J^+
\end{equation}
We conclude that $I^+ = J^+ = K^+$, thus implying Theorem \ref{thm:main}. 

\end{proof}

\begin{remark} 
\label{rem:replace}

In fact, Proposition \ref{prop:realize} shows that we could have replaced the ideal $J^+$ of \eqref{eqn:j} in the preceding argument by
\begin{equation}
\label{eqn:jj}
(\text{coefficients of } \rho_{Z,\tau_\ell} )_{Z \text { distinguished zig-zag}} \subset \tUUp
\end{equation}
where $\{\tau_\ell\}_{\ell \in \BZ}$ denotes any collection of Laurent polynomials of homogeneous degree $\ell$ that satisfy \eqref{eqn:tau neq}. Therefore, we conclude that $J^+$ is equal to the ideal \eqref{eqn:jj} for any such collection of Laurent polynomials $\tau_\ell$.

\end{remark}

\medskip

\begin{proposition}
\label{prop:equivalent}

When $Z$ is the particular minimal zig-zag in Figure 3, the relations $\rho_{Z,\tau_\ell} = 0$ (for any collection $\{\tau_\ell\}_{\ell \in \BZ}$ of Laurent polynomials of homogeneous degree $\ell$ satisfying \eqref{eqn:tau neq})  are equivalent to relation \eqref{eqn:loop serre} in the algebra $\tUUp$.

\end{proposition}

\medskip

\begin{proof} Let $\CL = \text{LHS of \eqref{eqn:loop serre}}$. Since $Z$ is minimal, Claim \ref{claim:stepping stone} shows that
$$
\Big \langle \rho_{Z,\tau}(x_s,\dots,x_t,y_{s' = t'}), R \Big \rangle_{\tU\CV} = \text{non-zero multiple of } R\left(q^s,\dots, q^t, q^{s'=t'}\right)
$$
for any homogeneous rational function $R \in \CV_{-(1-d_{ij})\bs^i-\bs^j}$, where $\tau$ is chosen so that its homogeneous degree balances out that of $R$. Above, $s' = t' = \frac {s+t}2$ and $s = t + 2d_{ij}$. By the non-degeneracy of the pairing $\langle \cdot , \cdot\rangle_{\tU\CV}$, the fact that
\begin{equation}
\label{eqn:pair serre}
\Big \langle \CL, R \Big \rangle_{\tU\CV} = \text{a delta function times } R\left(q^s,\dots, q^t, q^{s'=t'}\right)
\end{equation}
would imply that $\CL$ is proportional to $\rho_{Z,\tau}$, which is precisely what we need to prove. In order to prove \eqref{eqn:pair serre}, let us invoke \eqref{eqn:pair formula 1} to obtain
$$
\Big \langle e_i(z_1) \dots e_i(z_k) e_j(w) e_i(z_{k+1}) \dots e_i(z_{1-d_{ij}}), R \Big \rangle = $$
$$
\left[ \frac {R(z_1,\dots,z_{1-d_{ij}},w)}{\prod_{1\leq a < b \leq 1-d_{ij}} \zeta_{ii} \left(\frac {z_b}{z_a} \right) \prod_{a=1}^k \zeta_{ji}\left(\frac w{z_a} \right)\prod_{a=k+1}^{1-d_{ij}} \zeta_{ij}\left(\frac {z_a}w \right)} \right]^{\text{ expand as}}_{|z_1| \gg \dots \gg |z_k| \gg |w| \gg |z_{k+1}| \gg \dots \gg |z_{1-d_{ij}}|}
$$
Taking $\Sym \sum_{k = 0}^{1-d_{ij}} (-1)^k {1-d_{ij} \choose k}_{q}$ of the expression above yields $\langle \CL, R \rangle_{\tU\CV}$ in the left-hand side, while in the right-hand side we obtain (by invoking formulas (6.48) and (6.49) of \cite{NT}, themselves based on an identity developed in \cite{DJ})
$$
\text{non-zero multiple of }\delta \left(\frac {z_1}{wq^{d_{ij}}} \right) \delta \left(\frac {z_2}{wq^{d_{ij}+2}} \right) \dots  \delta \left(\frac {z_{1-d_{ij}}}{wq^{-d_{ij}}} \right) R(z_1,\dots,z_{1-d_{ij}},w)
$$
By \eqref{eqn:delta identity}, the quantity above vanishes if and only if $R(q^s,\dots, q^t, q^{s'=t'}) = 0$, precisely as we needed to prove in \eqref{eqn:pair serre}. 

\end{proof}

\noindent The proof of Proposition \ref{prop:equivalent} shows that just like the power series expansions of different linear combinations of rational functions can produce the same product of delta functions, different relations in the algebra $\tUUp$ can be equivalent modulo \eqref{eqn:rel quad 1}.

\medskip

\section{Words}
\label{sec:words}

\medskip

\noindent In the present Section, we will use techniques stemming from the combinatorics of words in order to prove Proposition \ref{prop:coincide}, which will also allow us to construct bases of the shuffle algebras $\CS^\pm \subset \CV^\pm$. We will also prove Proposition \ref{prop:final non-deg} by ``covering" our infinite-dimensional quantum loop groups and shuffle algebras by finite-dimensional vector subspaces, also indexed by words.

\medskip

\subsection{}

The treatment of the present Section is inspired (in chronological order) by \cite{LR, L, NT} and most closely \cite{Quiver 1}. Let us consider the set of letters
$$
\left\{ i^{(k)} \right\}_{i \in I, k \in \BZ}
$$
A word is simply a sequence of letters
\begin{equation}
\label{eqn:word}
w = \left[ i_1^{(k_1)} \dots i_n^{(k_n)} \right]
\end{equation}
We will call $|w| := n$ the length of a word as above, and 
\begin{equation}
\label{eqn:sequence of exponents}
\ow = (k_1,\dots,k_n)
\end{equation}
the sequence of exponents of $w$. By analogy with the constructions of Subsection \ref{sub:def quad}, the degree of $w$ is defined as
$$
\deg w = (\bs^{i_1} + \dots + \bs^{i_n}, k_1 + \dots + k_n) \in \nn \times \BZ
$$ 
Denote the set of all words by $\CW$. To $w \in \CW$ as above, we associate the element
\begin{equation}
\label{eqn:associated word}
e_w = e_{i_1,k_1} \dots e_{i_n, k_n} \in \tUUp
\end{equation}
In the opposite algebra, we will use the notation
\begin{equation}
\label{eqn:associated word}
f_w = f_{i_1,-k_1} \dots f_{i_n, -k_n} \in \tUUm
\end{equation}

\subsection{} Let us fix a total order $<$ on $I$, and associate to it the following total order on the set of letters:
$$
i^{(k)} < j^{(l)} \quad \text{if} \quad \begin{cases} k > l \\ \ \ \text{or} \\ k = l \text{ and } i<j \end{cases}
$$
Then we have the corresponding total lexicographic order on the set of words:
\begin{equation}
\label{eqn:lex order}
\left[i_1^{(k_1)} \dots i_n^{(k_n)} \right] < \left[j_1^{(l_1)} \dots j_m^{(l_m)} \right]
\end{equation}
if $i_1^{(k_1)} = j_1^{(l_1)}, \dots, i_x^{(k_x)} = j_x^{(l_x)}$ and either $i_{x+1}^{(k_{x+1})} < j_{x+1}^{(l_{x+1})}$ or $x = n < m$.

\medskip

\begin{definition}
\label{def:non-increasing}

A word \eqref{eqn:word} is called \textbf{non-increasing} if
\begin{equation}
\label{eqn:non-increasing}
\begin{cases} k_a < k_b + \Big| \Big\{ s \in \{a,\dots,b-1\} \text{ s.t. } i_s \neq i_b \Big\}\Big| \\ \ \ \text{or} \\ k_a = k_b + \Big| \Big\{ s \in \{a,\dots,b-1\} \text{ s.t. } i_s \neq i_b \Big\} \Big| \text{ and } i_a \geq i_b \end{cases}
\end{equation}
for all $1 \leq a < b \leq n$. Let $\CW_{\emph{non-inc}}$ denote the set of non-increasing words.

\end{definition}

\medskip

\begin{lemma}
\label{lem:finite}

There are finitely many non-increasing words of given degree, which are bounded above by any given word $v$. 

\end{lemma} 

\medskip

\begin{proof} Let us assume we are counting non-increasing words $[i_1^{(k_1)} \dots i_n^{(k_n)}]$ with $k_1+\dots+k_n = k$ for fixed $n$ and $k$. The fact that such words are bounded above implies that $k_1$ is bounded below. But then the inequality \eqref{eqn:non-increasing} implies that $k_2,\dots,k_n$ are also bounded below. The fact that $k_1+\dots+k_n$ is fixed implies that there can only be finitely many choices for the exponents $k_1,\dots,k_n$. Since there are also finitely many choices for $i_1,\dots,i_n \in I$, this concludes the proof.  

\end{proof}

\medskip

\begin{proposition} 
\label{prop:basis}

The set $\{e_w\}_{w \in \CW_{\emph{non-inc}}}$ is a linear basis of $\tUUp$.

\end{proposition}

\medskip

\begin{proof} By running the proof of \cite[Proposition 3.11]{Quiver 1} for the quiver $Q$ with vertex set $I$ and a single arrow from $i$ to $j$ if and only if $i < j$, we can prove that
\begin{equation}
\label{eqn:basis}
\{e_w\}_{w \in \CW_{\text{non-inc}}}
\end{equation}
linearly spans $\tUUp$. Indeed, the only thing this fact requires in \loccit is relation \cite[(3.20)]{Quiver 1}, which takes the same form as our \eqref{eqn:quad intro} (albeit with different parameters instead of $q^{d_{ij}}$, though these do not affect the validity of the argument). In order to show that the elements \eqref{eqn:basis} are linearly independent in $\tUUp$, we invoke \cite[Proposition 2.16]{Quiver 3}, for the same quiver $Q$ as above. We will explain the details (without proof) in Subsection \ref{sub:lead}, see the discussion  between formula \eqref{eqn:leading word pairing} and the proof of Proposition \ref{prop:non-deg}. 

\end{proof}

\begin{remark}
\label{rem:more}

Note that the proof above actually establishes the stronger statement (cf. \cite[Proposition 3.11]{Quiver 1}): for any $i_1,\dots,i_n \in I$ and $k_1,\dots,k_n \in \BZ$ we have
\begin{equation}
	\label{eqn:straighten}
	e_{i_1,k_1} \dots e_{i_n,k_n} \in \mathop{\sum_{[j_1^{(l_1)} \dots j_n^{(l_n)}] \text{ non-increasing and } \geq \ [i_1^{(k_1)} \dots i_n^{(k_n)}]}}_{\min(k_a) - \beta(n) \leq \min(l_a) \leq \max(l_a) \leq \max(k_a) + \beta(n)} \BK \cdot e_{j_1,l_1} \dots e_{j_n,l_n}
\end{equation}
for some constant $\beta(n)$ that only depends on $n$. The stronger statement above says that, as we use \eqref{eqn:quad intro} to reorder the $e_{i,k}$'s so that the right-hand side of the expression above consists only of non-increasing words, we only increase/decrease the $k_a$'s by a bounded amount during the whole procedure. 

\end{remark}

\medskip

\subsection{}
\label{sub:finitely many}

Let us give an application of property \eqref{eqn:straighten}. As we saw in Remark \ref{rem:replace}, one may replace $J^+$ by the ideal generated by $\rho_{Z,\tau_k}$, as $Z$ runs over all distinguished zig-zags as in Figure 7, and $\tau_k(x_s,\dots,x_t,y_{s'},\dots,y_{t'})$ runs over generic Laurent polynomials of any homogeneous degree $k \in \BZ$. For given $Z$, let us write
$$
\bn = \bs^i  \cdot \left( \frac {t-s}2 + 1 \right) + \bs^j  \cdot \left( \frac {t'-s'}2 + 1 \right) \qquad \Rightarrow \qquad n = |\bn| = \frac {t-s+t'-s'}2+2
$$
Then let us choose $\tau_k$ such that
$$
\tau_{k+nd}(x_s,\dots,x_t,y_{s'},\dots,y_{t'}) = \tau_k (x_s,\dots,x_t,y_{s'},\dots,y_{t'}) \cdot (x_s\dots x_t y_{s'} \dots y_{t'})^d
$$
for all $k,d \in \BZ$. The effect that this choice has is that
\begin{align}
&\text{if } \qquad \rho_{Z,\tau_k} = \mathop{\sum_{i_1,\dots,i_n \in I}}_{k_1,\dots,k_n \in \BZ} \text{coefficient} \cdot e_{i_1,k_1} \dots e_{i_n,k_n} \label{eqn:first} \\
&\text{then } \rho_{Z,\tau_{k+nd}} = \mathop{\sum_{i_1,\dots,i_n \in I}}_{k_1,\dots,k_n \in \BZ} \text{coefficient} \cdot e_{i_1,k_1+d} \dots e_{i_n,k_n+d} \label{eqn:second}
\end{align}
for all $k,d \in \BZ$. Thus, the finitely many coefficients in \eqref{eqn:first} for $k \in \{1,\dots,n\}$ determine the coefficients for all $k$. We conclude that
\begin{equation}
\label{eqn:finitely many}
\rho_{Z,\tau_k} = \mathop{\sum_{i_1,\dots,i_n \in I}}_{\frac kn - \gamma(n) \leq k_1,\dots,k_n \leq \frac kn + \gamma(n)} \text{coefficient} \cdot e_{i_1,k_1} \dots e_{i_n,k_n}
\end{equation}
for all $k \in \BZ$, for some large enough $\gamma(n)$ that only depends on $n$. 

\medskip

\subsection{}
\label{sub:lead}

Given a monomial in $\{z_{i\bullet}\}_{i \in I}^{\bullet \in \BN}$, we will consider all ways to order its variables
\begin{equation}
\label{eqn:monomial}
\mu = z_{i_1\bullet_1}^{-l_1} \dots z_{i_n \bullet_n}^{-l_n}
\end{equation}
such that $(i_a,\bullet_a) \neq (i_b,\bullet_b)$ for all $a \neq b$ (in other words, we do not split the powers in monomials). The \textbf{leading word} of such a monomial is defined as the (lexicographically) largest word
\begin{equation}
\label{eqn:leading word}
w_\mu = \left [i_1^{(k_1)} \dots i_n^{(k_n)} \right] 
\end{equation}
where:
\begin{equation}
\label{eqn:formula bijection}
k_a = l_a - \Big| \Big\{ s< a \text{ s.t. } i_s > i_a \Big\} \Big|  + \Big| \Big\{ t > a \text{ s.t. } i_t < i_a \Big\} \Big|
\end{equation}
among all ways to order the variables in \eqref{eqn:monomial}. Following the proof of \cite[Lemma 4.8]{Quiver 1} for the quiver $Q$ defined in the previous Subsection (vertex set $I$ and a single arrow from $i$ to $j$ iff $i < j$), the leading word of any monomial is non-increasing. It is easy to see that the leading word only depends on $\text{Sym } \mu$, namely the symmetrization with respect to $z_{i1}, z_{i2}, \dots$ for each $i \in I$ separately.  

\medskip

\noindent Conversely, the \textbf{associated polynomial} of a non-increasing word \eqref{eqn:leading word} is defined by $\Sym \ \eqref{eqn:monomial}$, with the $l_a$'s determined by \eqref{eqn:formula bijection}. The functions ``leading word" and ``associated polynomial" yield mutually inverse bijections between the sets
\begin{equation}
\label{eqn:bijection}
\begin{tikzcd}
\text{Sym(monomial)}\ar[r, bend right] &
\text{non-increasing words}\arrow[l, bend right]
\end{tikzcd}
\end{equation}
More generally, the leading word of any $R \in \CV_{-\bn}$, denoted by $\text{lead}(R)$, will be the lexicographically largest of the leading words \eqref{eqn:leading word} for all the monomials which appear with non-zero coefficient in the Laurent polynomial
\begin{equation}
\label{eqn:all monomials}
R(\dots,z_{ia},\dots) \prod^{i < j \in I}_{1\leq a \leq n_i, 1 \leq b \leq n_j} \left(1 - \frac {z_{jb}}{z_{ia}} \right)
\end{equation}
(the reason we place $z_{jb}$ in the numerator and $z_{ia}$ in the denominator will be made clear in \eqref{eqn:leading word pairing}). Conversely, any non-increasing word $w$ appears as the leading word
\begin{equation}
\label{eqn:w as lead}
w = \text{lead}\left[ \frac {\text{Sym } \mu}{\prod^{i < j \in I}_{1\leq a \leq n_i, 1 \leq b \leq n_j} \left(1 - \frac {z_{jb}}{z_{ia}} \right)} \right]
\end{equation}
with $\mu$ and $w$ connected by \eqref{eqn:monomial}--\eqref{eqn:formula bijection}. Note that the denominator of the formula above is already symmetric, so it can be freely moved in and out of the symmetrization.

\medskip

\noindent We have the following analogue of \cite[Proposition 3.24]{Quiver 1}, which is actually the initial motivation for our introduction of non-increasing words.
\begin{equation}
\label{eqn:leading word pairing}
\Big \langle e_w, R \Big \rangle_{\tU\CV} \text{ is } \begin{cases} \neq 0 &\text{if }w = \text{lead}(R) \\ = 0 &\text{if }w > \text{lead}(R) \end{cases}
\end{equation}
To prove this, let us note that formula \eqref{eqn:pair formula 1} immediately implies that
\begin{equation}
	\label{eqn:hart}
\qquad \Big \langle e_{\left [i_1^{(k_1)} \dots i_n^{(k_n)} \right]}, R \Big  \rangle_{\tU\CV} \text{ is } \begin{cases} \neq 0 &\text{if }(k_1,\dots,k_n) = (\ell_1,\dots,\ell_n) \\ = 0 &\text{if }(k_1,\dots,k_n) < (\ell_1,\dots,\ell_n) \text{ lexicographically }\end{cases}
\end{equation}
where $R \sim z_1^{-\ell_1} \dots z_n^{-\ell_n}$ plus terms of lower order as $|z_1|\gg \dots \gg |z_n|$ (above and hereafter, the notation $\sim$ refers to proportionality up to a non-zero constant). If
$$
R(\dots,z_{ia},\dots) = \frac {\text{Sym } \mu}{\prod^{i < j \in I}_{1\leq a \leq n_i, 1 \leq b \leq n_j} \left(1 - \frac {z_{jb}}{z_{ia}} \right)}
$$
we have
$$
\ell_a = l_a -\Big| \Big\{s<a \text{ s.t. } i_s > i_a \Big\} \Big| + \Big| \Big\{t>a \text{ s.t. } i_t < i_a \Big\} \Big|
$$
where $\text{Sym }\mu \sim z_1^{-l_1} \dots z_n^{-l_n}$ plus terms of lower order as $|z_1|\gg \dots \gg |z_n|$. The word $\text{lead}(R)$ is simply the lexicographically largest of the words
$$
\left[ i_1^{(\ell_1)} \dots i_n^{(\ell_n)} \right]
$$
that can arise for various choices of $i_1,\dots,i_n \in I$. Therefore, \eqref{eqn:hart} implies \eqref{eqn:leading word pairing}.

\medskip 

\noindent Formula \eqref{eqn:leading word pairing} immediately implies the linear independence of the elements $e_w$, as $w$ runs over non-increasing words. Indeed, if one were able to write such an element $e_w$ as a linear combination of elements $e_v$ for various $v>w$, then we would contradict \eqref{eqn:leading word pairing} for $R$ being the rational function in the right-hand side of \eqref{eqn:w as lead}. 

\medskip

\begin{proof} \emph{of Proposition \ref{prop:non-deg}:} With formulas \eqref{eqn:pair formula 1}--\eqref{eqn:pair formula 2} in mind, the non-degeneracy of the pairing in the $\CV^\pm$ argument is simply restating the well-known fact that if a rational function $F$ in variables $z_1,\dots,z_n$ vanishes when expanded as a power series (in some relative order of its variables) then $F = 0$.

\medskip

\noindent Let us now prove non-degeneracy in the $\tUUpm$ argument (we will focus on $\tUUp$ without loss of generality). Consider any non-zero element
$$
\phi = \sum_{v \in \CW_{\text{non-inc}}} \underbrace{\gamma_v}_{\in \BK} \cdot e_v
$$
and let $w$ be the smallest word such that $\gamma_w \neq 0$. Since we showed in \eqref{eqn:w as lead} that there exists an element $R \in \CV^-$ such that $\text{lead}(R) = w$, then \eqref{eqn:leading word pairing} implies that
$$
\Big \langle \phi, R \Big \rangle_{\tU\CV} = \sum_{v \in \CW_{\text{non-inc}}} \gamma_v \Big \langle e_v, R \Big \rangle_{\tU\CV} = \gamma_w \Big \langle e_w,R \Big \rangle_{\tU\CV} \neq 0
$$
This implies the non-degeneracy of the pairing in the $\tUUp$ argument.  

\end{proof}

\medskip

\subsection{} We are now ready to prove Proposition \ref{prop:coincide}. Recall the surjective homomorphisms
$$
\widetilde{\Upsilon}^\pm : \tUUpm \twoheadrightarrow \oCS^\pm
$$
and let $E_w = \widetilde{\Upsilon}^+(e_w)$, $F_w = \widetilde{\Upsilon}^-(f_w)$ for any word $w$. 

\medskip

\begin{definition} A word $w$ is called \textbf{standard} if 
$$
E_w \notin \sum_{v > w} \BK \cdot E_v
$$
Let $\CW_{\emph{stan}}$ denote the set of standard words. 
 
\medskip

\end{definition}

\noindent Formula \eqref{eqn:straighten} implies that any standard word is non-increasing. The converse is definitely not true, since the definitions readily imply the fact that
$$
w \in \CW_{\text{non-inc}} \backslash \CW_{\text{stan}} \quad \Leftrightarrow \quad e_w \in \sum_{v > w} \BK \cdot e_v + \text{Ker }\widetilde{\Upsilon}^+
$$
Thus, the discrepancy between standard and non-increasing words is a measure of $K^+$.

\medskip

\begin{proof} \emph{of Proposition \ref{prop:coincide}:} We will prove \eqref{eqn:coincide} for $\pm = -$, as the case $\pm = +$ is analogous. All words in the present proof will have fixed length $n$. For any words
$$
v = \left[ i_1^{(k_1)} \dots i_n^{(k_n)} \right] \qquad \text{and} \qquad w  = \left[ j_1^{(l_1)} \dots j_n^{(l_n)} \right]
$$
the following formula is easy to deduce from \eqref{eqn:pair formula 1} (see \cite[Remark 3.16]{Quiver 1} for a proof in an almost identical setup)
$$
\Big \langle E_v, F_w \Big \rangle_{\oCS\oCS} = \int_{|z_1| \gg \dots \gg |z_n|} \mathop{\sum_{\sigma \in S(n)}}_{i_a = j_{\sigma(a)}, \forall a} z_1^{k_1 - l_{\sigma(1)}} \dots z_n^{k_n - l_{\sigma(n)}} \prod^{a<b}_{\sigma(a)>\sigma(b)} \frac {\zeta_{i_ai_b} \left(\frac {z_a}{z_b} \right)}{\zeta_{i_bi_a} \left(\frac {z_b}{z_a} \right)} \prod_{a=1}^n Dz_a  
$$
Expanding the integrand as $|z_1| \gg \dots \gg |z_n|$ reveals that $\langle E_v, F_w \rangle_{\oCS\oCS} \neq 0$ only if there exists $\sigma \in S(n)$ and integers
\begin{equation}
\label{eqn:integers}
\Big\{ c_{a,b} \leq 0 \Big\}_{a < b, \sigma(a) > \sigma(b)}
\end{equation}
such that
\begin{equation}
\label{eqn:equality}
k_a + \sum^{t > a}_{\sigma(t) < \sigma(a)} c_{at} - \sum^{s<a}_{\sigma(s) > \sigma(a)} c_{sa} = l_{\sigma(a)}
\end{equation}
Consider the infinite directed graph $G_{n}$ whose vertices are collections of integers
$$
(k_1,\dots,k_n) \text{ such that } k_a \leq k_{a+1} + 1
$$
(the inequality stems from the definition of non-increasing words) and whose edges are
\begin{equation}
\label{eqn:edge directed graph}
(k_1,\dots,k_n) \rightarrow ( l_1, \dots, l_n)
\end{equation}
if and only if there exists $\sigma \in S(n)$ and integers \eqref{eqn:integers} such that \eqref{eqn:equality} holds (although the graph is defined to be directed, it is easy to see that an edge \eqref{eqn:edge directed graph} exists if and only if the opposite edge exists). The preceding discussion implies that
\begin{equation}
\label{eqn:only if}
\Big \langle E_v, F_w \Big \rangle_{\oCS\oCS} \neq 0 \quad \Rightarrow \quad \text{exists edge } \overline{v} \rightarrow \overline{w}
\end{equation}
in $G_{n}$, for any non-increasing words $v,w$ (above, $\overline{v}$ denotes the sequence of exponents of the word $v$, see \eqref{eqn:sequence of exponents}). The following result was proved in \cite[Lemma 3.18]{Quiver 1}. 

\medskip

\begin{lemma} All connected components of $G_{n}$ are finite. \end{lemma}

\medskip

\noindent With the Lemma above in mind, let us define the following \underline{finite-dimensional} vector spaces, for any connected component $H \subset G_{n}$:
\begin{align*}
&\oCS^+_H = \sum_{w \in \CW_{\text{non-inc}}, \ \ow \in H} \BK \cdot E_w \subset \oCS^+ \\
&\oCS_H^- = \sum_{w \in \CW_{\text{non-inc}}, \ \ow \in H} \BK \cdot F_w \subset \oCS^-
\end{align*}
Because of \eqref{eqn:only if}, the descended pairing \eqref{eqn:descend pair 3} satisfies
$$
\Big \langle \oCS^+_H, \oCS_{H'}^-  \Big \rangle_{\oCS\oCS} = 0
$$
if $H \neq H'$. However, since the descended pairing is non-degenerate in both arguments (a property which it inherits from the non-degeneracy of $\langle\cdot,\cdot \rangle_{\tU\CV}$ and $\langle\cdot,\cdot \rangle_{\CV\tU}$ in the $\CV$ argument), then so is its restriction
\begin{equation}
\label{eqn:restricted pairing}
\oCS^+_H \otimes \oCS_{H}^- \xrightarrow{\langle \cdot, \cdot \rangle_{\oCS\oCS}} \BK
\end{equation}
for any connected component $H \subset G_{n}$. The following statements are straightforward consequences of the finite-dimensionality of the vector spaces $\oCS^\pm_H$.

\medskip

\begin{proposition}
\label{prop:direct}

For any $n\in \BN$, we have
\begin{equation}
\label{eqn:direct sum 1}
\bigoplus_{|\bn|=n} \oCS^\pm_{\bn} = \bigoplus_{H \text{ a connected component of } G_n} \oCS^\pm_H
\end{equation}
and
\begin{equation}
\label{eqn:direct sum 2}
\oCS^+_H = \bigoplus_{w \text{ standard}}^{\ow \in H} \BK \cdot E_w \qquad \text{and} \qquad \oCS^-_H = \bigoplus_{w \text{ standard}}^{\ow \in H} \BK \cdot F_w
\end{equation}

\end{proposition}

\medskip

\begin{proof} Let us prove the statements for $\pm = +$. Because the $E_w$'s span $\oCS^+$ as $w$ runs over all non-increasing words, all that we need to do to prove \eqref{eqn:direct sum 1} is to show that there are no linear relations among the various direct summands of the RHS. To this end, assume that we had a relation
$$
\sum_{H \text{ a connected component of } G_n} \alpha_H = 0
$$
for various $\alpha_H \in \oCS_H^+$. Pairing the relation above with a given $\oCS^-_H$ implies that
$$
\left \langle \alpha_H, \oCS^-_H \right \rangle_{\oCS\oCS} = 0
$$
Because the pairing \eqref{eqn:restricted pairing} is non-degenerate, this implies that $\alpha_H = 0$. As for \eqref{eqn:direct sum 2}, it holds because any vector space spanned by vectors $\alpha_1,\dots,\alpha_k$ has a basis consisting of those $\alpha_i$'s which cannot be written as linear combinations of $\{\alpha_j\}_{j > i}$. 

\end{proof}

\medskip 

\noindent We are now poised to complete the proof of Proposition \ref{prop:coincide}. Consider any $R \in \CS_{-\bn}$ with $|\bn| = n$. From \eqref{eqn:pair formula 1}, it is easy to see that
$$
\Big \langle E_{[i_1^{(k_1)} \dots i_n^{(k_n)} ]}, R \Big \rangle_{\oCS\CS} = 0
$$
if $k_1$ is small enough. However, by Lemma \ref{lem:finite}, there are only finitely many non-increasing words $w$ of given degree with $k_1$ bounded below. This implies that
$$
\Big \langle E_w, R \Big \rangle_{\oCS\CS} \neq 0
$$
only for finitely many non-increasing words $w$. Letting $H_1, \dots, H_t \subset G_n$ denote the connected components which contain the sequences of exponents of the aforementioned words, then \eqref{eqn:leading word pairing} and the non-degeneracy of the pairings \eqref{eqn:restricted pairing} imply that there exists an element
$$
R' \in \oCS^-_{H_1} \oplus \dots \oplus \oCS^-_{H_t} \subset \oCS^-
$$
such that $\langle E_w, R \rangle_{\oCS\CS} = \langle E_w, R' \rangle_{\oCS\oCS}$ for all non-increasing words $w$. Then the non-degeneracy of the descended pairing \eqref{eqn:descend pair 1} in the second argument (which follows from Propositions \ref{prop:non-deg} and \ref{prop:descends}) implies that $R = R' \in \oCS^-$, as we needed to prove. The aforementioned non-degeneracy also implies the non-degeneracy of the pairing \eqref{eqn:pairing shuffle}. 

\end{proof}

\medskip

\subsection{}

Throughout the remainder of the present Section, we will only consider words
$$
w = \left[i_1^{(k_1)} \dots i_n^{(k_n)} \right]
$$
of fixed degree $\deg w = (\bn, k) \in \nn \times \BZ$. For any \underline{finite} set $T \subset \CW_{\text{non-inc}}$ of words, we define
\begin{equation}
\label{eqn:u finite}
\tUUpT = \bigoplus_{w \in T} \BK \cdot e_w
\end{equation}
(recall that as $w$ runs over $\CW_{\text{non-inc}}$, the $e_w$'s yield a basis of $\tUUp$). Also let
\begin{equation}
\label{eqn:v finite}
\CV^{-,T} \subset \CV_{-\bn,-k}
\end{equation}
denote the subspace of rational functions $R$, such that all monomials appearing in the Laurent polynomial \eqref{eqn:all monomials} have leading word in $T$. Thus, the restriction
\begin{equation}
\label{eqn:finite pairing}
\tUUpT \otimes \CV^{-,T} \xrightarrow{\langle \cdot, \cdot \rangle_{\tU\CV}^T} \BK
\end{equation}
is a non-degenerate pairing of finite-dimensional vector spaces (indeed, because the two vector spaces have the same dimension, it suffices to show non-degeneracy in the second factor, which follows immediately from \eqref{eqn:leading word pairing}). 

\medskip

\begin{proof} \emph{of Proposition \ref{prop:final non-deg}:} Let us consider a finite set of non-increasing words $T$, of some fixed degree $(\bn,k)$, to be chosen in \eqref{eqn:choice of t}. The only thing we postulate for the time being is that the set $T$ can be chosen ``arbitrarily large", i.e. to contain any given finite set of words. Let $\CS^{-,T} = \CS^- \cap \CV^{-,T}$. Proposition \ref{prop:realize} implies that
$$
\CS^{-,T} \subseteq \left(J^+ \cap \tUUpT \right)^\perp
$$
where the orthogonal complement is defined with respect to the pairing \eqref{eqn:finite pairing}. Our goal will be to prove the opposite inclusion, namely
\begin{equation}
\label{eqn:opposite inclusion}
\left(J^+ \cap \tUUpT\right)^\perp \subseteq \CS^{-,T}
\end{equation}
Let us first use \eqref{eqn:opposite inclusion} to conclude the proof of Proposition \ref{prop:final non-deg}. Because \eqref{eqn:finite pairing} is a non-degenerate pairing of finite-dimensional vector spaces, we would have the following implication
$$
\CS^{-,T} = \left(J^+ \cap \tUUpT\right)^\perp  \qquad \Rightarrow \qquad J^+ \cap \tUUpT = \left(\CS^{-,T}\right)^\perp
$$
Since every element $\phi \in \tUUp$ lies in $\tUUpT$ for some large enough finite set $T$, if $\phi$ pairs trivially with the whole of $\CS^-$, then the formula above implies that $\phi \in J^+$. This establishes the non-degeneracy of the pairing \eqref{eqn:pairing final} in the first argument. Since the non-degeneracy in the second argument is inherited from that of $\langle \cdot , \cdot \rangle_{\tU\CV}$ (see Proposition \ref{prop:non-deg}, proved earlier in this Section), we would be done. \\

\noindent To prove \eqref{eqn:opposite inclusion}, we will revisit the proof of Proposition \ref{prop:realize}. Specifically, after we stated Claim \ref{claim:stepping stone}, we showed that given $R \in \CV^-$,
\begin{equation}
\label{eqn:easy}
\Big \langle \phi, R \Big \rangle_{\tU\CV} = 0, \ \forall \phi \in J^+ \qquad \qquad \qquad \quad \Rightarrow  R \in \CS^-
\end{equation}
For the task at hand, we need to show that given $R \in \CV^{-,T}$,
\begin{equation}
\label{eqn:hard}
\Big \langle \phi, R \Big \rangle_{\tU\CV} = 0, \ \forall \phi \in J^+ \cap \tUUpT \qquad \Rightarrow R \in \CS^{-,T}
\end{equation}
Concretely, we fix $R \in \CV^{-,T}$ and we will show that the only $\phi$'s that Proposition \ref{prop:realize} needs in order to ensure the implication \eqref{eqn:easy} actually lie in $\tUUpT$; this would yield the implication \eqref{eqn:hard} and we would be done. To this end, let us choose
\begin{equation}
\label{eqn:choice of t}
T = \left \{w = \left[i_1^{(k_1)} \dots i_n^{(k_n)} \right] \text{ such that }\deg w = (\bn, k) \text{ and } \right. 
\end{equation}
$$
\left. \sum_{s \in A} k_s \geq -M|A| + m|A|^2 - \Big|\Big\{(s>t) \text { s.t. } s \in A, t \notin A, i_s \neq i_t \Big\} \Big|, \forall A \subseteq \{1,\dots,n\} \right\} 
$$
for natural numbers $m,M$ (since $M$ can be arbitrarily large, this would ensure the fact that any finite set of non-increasing words can be contained in $T$). The reason we subtract by $|(s>t) \text { s.t. } s \in A, t \notin A, i_s \neq i_t|$ in the inequality above is the straightforward fact (which we leave to the interested reader) that the inequality holds for the leading word \eqref{eqn:leading word} of a monomial \eqref{eqn:monomial} if and only if it holds for the analogous word associated to any other order of the variables in the monomial.

\medskip

\noindent As shown in the proof of Proposition \ref{prop:realize} (we will reuse the notation therein), the only $\phi$'s we need to consider in \eqref{eqn:easy} are those of the form
\begin{equation}
\label{eqn:middle}
e_{i_1,k_1} \dots e_{i_p,k_p} \rho_{Z,\tau_{k-k_1-\dots-k_p}}
\end{equation}
for some monomial that appears in the numerator of $R$
\begin{equation}
\label{eqn:look at numerator}
R = \frac {\ldots + u_1^{-k_1} \dots u_p^{-k_p}  R_{k_1,\dots,k_p}(x_s,\dots,x_t,y_{s'},\dots,y_{t'}) + \dots}{\prod_{a < b} \left(1-\frac {u_b}{u_a} \right) \prod_{b,c} \left(1-\frac {z_c}{u_b}\right)}
\end{equation}
(see \eqref{eqn:partial expand}) with non-zero $R_{k_1,\dots,k_p}$. Moreover, we may assume the word
\begin{equation}
\label{eqn:word is non-increasing}
\left[i_1^{(k_1)} \dots i_p^{(k_p)} \right]
\end{equation}
is non-increasing. This is because in the proof of Proposition \ref{prop:realize}, one can reorder the variables $u_1,\dots,u_p$ arbitrarily; choosing the order which leads to the maximal word \eqref{eqn:word is non-increasing} ensures that the aforementioned word is non-increasing, as explained in Subsection \ref{sub:lead}. 

\medskip

\noindent Using \eqref{eqn:finitely many}, we may write \eqref{eqn:middle} is a linear combination of products of the form
\begin{equation}
\label{eqn:bottom}
e_{i_1,k_1} \dots e_{i_p,k_p} e_{i_{p+1},k_{p+1}} \dots e_{i_n,k_n}
\end{equation}
where $k_{p+1},\dots,k_n$ are within $\gamma(n-p)$ from their average value. It remains to show that any product \eqref{eqn:bottom} can be expressed, using \eqref{eqn:straighten}, as a linear combination of products of $e_{i,k}$'s that correspond to non-increasing words in $T$.

\medskip

\noindent Let us consider the constant $b = 2\max(\beta(1),\dots,\beta(n))+2n$, with the $\beta$'s as in \eqref{eqn:straighten}. In the product \eqref{eqn:bottom}, we consider the largest index $x \in \{0,\dots,p\}$ such that
\begin{equation}
\label{eqn:final ineq}
k_x < k_{x+1} - b
\end{equation}
(we make the convention that $k_0 = -\infty$). We may use relation \eqref{eqn:straighten} to write 
$$
e_{i_{x+1},k_{x+1}} \dots e_{i_n,k_n} = \sum_{\left[j_{x+1}^{(l_{x+1})} \dots j_n^{(l_n)} \right] \in \CW_{\text{non-inc}}} \text{coefficient} \cdot e_{j_{x+1},l_{x+1}} \dots e_{j_n,l_n} 
$$
We need to make two observations about the $l$'s that appear in the formula above.

\medskip

\begin{itemize}[leftmargin=*]

\item All the numbers $l_{x+1},\dots,l_n$ are within a global constant away from their average

\medskip

\item The large difference between $k_x$ and $k_{x+1}$ ensures that all concatenated words
\begin{equation}
\label{eqn:concatenated words}
\left[ i_1^{(k_1)} \dots i_x^{(k_x)} j_{x+1}^{(l_{x+1})} \dots j_n^{(l_n)} \right]
\end{equation}
which arise in the procedure above are non-increasing (recall that the word \eqref{eqn:word is non-increasing} was non-increasing to begin with, and thus so are all of its prefixes)

\end{itemize}

\medskip

\noindent Thus, it remains to show that the concatenated words that appear in \eqref{eqn:concatenated words} are in $T$. Because $R \in \CV^{-,T}$, then \eqref{eqn:choice of t} and \eqref{eqn:look at numerator} imply that
\begin{align}
&\sum_{s \in B} k_s \geq -M|A| + m|A|^2 - \Big|\Big\{(s > t) \text { s.t. } s \in A, t \notin A, i_s \neq i_t \Big\}\Big| \label{eqn:last 1} \\
&\sum_{s \in B} k_s + \sum_{s=x+1}^n l_s \geq -M|A| + m|A|^2 - \Big|\Big\{(s > t) \text { s.t. } s \in A, t \notin A, i_s \neq i_t \Big\}\Big| \label{eqn:last 2}
\end{align}
for $A = B$ or $A = B \sqcup \{x+1,\dots,n\}$ respectively, where $B \subseteq \{1,\dots,x\}$ is arbitrary. The latter formula holds because $l_{x+1}+\dots+l_n = k_{x+1} + \dots + k_n$. Assume for the purpose of contradiction that the defining property of $T$ is violated for
$$
A = B \sqcup C
$$
with $B \subseteq \{1,\dots,x\}$ and $C$ a proper subset of $\{x+1,\dots,n\}$, i.e.
\begin{equation}
\label{eqn:last 3}
\sum_{s \in B} k_s + \sum_{s \in C} l_s < -M|A| + m|A|^2 - \Big|\Big\{(s > t) \text { s.t. } s \in A, t \notin A, i_s \neq i_t \Big\} \Big|
\end{equation}
We claim that \eqref{eqn:last 1}--\eqref{eqn:last 2} and \eqref{eqn:last 3} are incompatible (for $m$ chosen large enough compared to the constant mentioned in the first bullet above). Indeed, letting $\mu$ be the average mentioned in the first bullet above, relations \eqref{eqn:last 1}--\eqref{eqn:last 2} imply
\begin{align}
&\sum_{s \in B} k_s \geq -M |B|+m|B|^2 - \dots \label{eqn:last 4} \\
&\sum_{s \in B} k_s + \mu(n-x) \geq -M(|B|+n-x) + m(|B|+n-x)^2 - \dots \label{eqn:last 5}
\end{align}
while \eqref{eqn:last 3} gives us
\begin{equation}
\label{eqn:last 6}
\sum_{s \in B} k_s + \mu  y < -M (|B|+y)+m(|B|+y)^2 + \dots 
\end{equation}
where $y = |C|$ lies in $\{1,\dots,n-x-1\}$, and the ellipses in formulas \eqref{eqn:last 4}--\eqref{eqn:last 6} denote global constants. Subtracting \eqref{eqn:last 4} from \eqref{eqn:last 6} yields
$$
\mu < -M+m(y +2|B|) + \dots
$$
and subtracting \eqref{eqn:last 6} from \eqref{eqn:last 5} yields 
$$
\mu > -M + m(n-x+y + 2|B|) - \dots
$$
The two inequalities above are incompatible if $m$ is chosen large enough compared to the global constants that appear among the ellipses, thus yielding the desired contradiction. 

\end{proof} 

\medskip

\section{$K$-theoretic Hall algebras}
\label{sec:k-ha}

\medskip

\noindent We will now study an important incarnation of quantum loop groups in geometric representation theory, namely $K$-theoretic Hall algebras associated to quivers. This is an idea that many mathematicians and physicists have been developing in recent decades \footnote{See for example the seminal work \cite{G}, where many of the $K$-theoretic constructions considered herein were developed, albeit with respect to a significantly different $\BC^*$ action. This choice of action leads to quite different notions of quantum loop groups (compare our relation \eqref{eqn:quad intro} with \cite[relation (7.6.1i)]{G}). The analogue of the present paper in the case of the quantum loop group considered in \cite{G} falls under the treatment of \cite[Sections 4 and 5]{Quiver 1} and \cite[Section 5]{Quiver 3}.}, but the particular incarnation we will be working with follows the philosophy of Schiffmann-Vasserot. We follow \cite{VV}, where the authors develop this line of thought in the setting of equivariant $K$-theory.

\medskip

\subsection{}
\label{sub:quiver}

Let us consider a quiver $Q$, i.e. an oriented graph, with vertex set $I$ and arrow set $E$. We allow multiple arrows between any two vertices, but \underline{no loops}. Given $\bn = (n_i)_{i \in I} \in \nn$, we will call 
$$
Z_{\bn} = \bigoplus_{\oij \in E} \Hom(\BC^{n_i}, \BC^{n_j})
$$
the vector space of $\bn$-dimensional representations of the quiver $Q$ (if there are multiple arrows from vertex $i$ to vertex $j$, then there are multiple copies of $\Hom(\BC^{n_i}, \BC^{n_j})$ in the direct sum above). Quiver representations will always be taken modulo change of basis, i.e. modulo the action of the algebraic group
$$
G_{\bn} = \prod_{i \in I} GL_{n_i}(\BC)
$$
by conjugation. Thus, the stack of $\bn$-dimensional representations of $Q$ is
$$
\fZ_{\bn} = Z_{\bn} / G_{\bn}
$$
In the present paper, we will mostly be concerned with the cotangent stack to $\fZ_{\bn}$. To describe it, consider the quadratic moment map
$$
T^*Z_{\bn} \xrightarrow{\mu_{\bn}} \bigoplus_{i \in I} \text{End}(\BC^{n_i})
$$
given by
\begin{equation}
\label{eqn:moment map}
\mu_{\bn} \left( \BC^{n_i} \xrightarrow{\phi_e} \BC^{n_j}, \BC^{n_j} \xrightarrow{\phi^*_e} \BC^{n_i} \right)_{\forall e = \oij \in E} = \sum_{e = \oij \in E} \left( \underbrace{\phi_e \phi_e^*}_{\in \text{End}(\BC^{n_j})} - \underbrace{\phi_e^* \phi_e}_{\in \text{End}(\BC^{n_i})} \right)
\end{equation}
Then we have the following presentation of the cotangent stack
$$
T^*\fZ_{\bn} = \mu_{\bn}^{-1}(0)/G_{\bn}
$$
for any $\bn \in \nn$.

\medskip

\subsection{}
\label{sub:torus action}

Fix a total order on $I$, and let us assume that all the arrows in the quiver $Q$ point from $i$ to $j$ only if $i<j$. This is not a significant restriction, as points of $T^*Z_{\bn}$ always come in pairs: for every linear map corresponding to an arrow, there exists a linear map pointing in the opposite direction. In fact, the only place where the orientation of arrows matters is in a certain monomial twist that defines the convolution product on the $K$-theoretic Hall algebra (see \cite[Subsection 2.2.2]{VV} or \cite[Remark 2.5]{Quiver 1}). Let us consider the torus action
\begin{equation}
\label{eqn:action}
\BC^* \curvearrowright T^*\fZ_{\bn}
\end{equation}
(called a ``normal weight function" in \cite{VV}) defined as follows. For every pair of vertices $i<j$ in $I$, fix an indexing
$$
e_1,\dots,e_{-d_{ij}}
$$
of the set of edges from $i$ to $j$ (this defines the numbers $d_{ij} = d_{ji} \in \BZ_{\leq 0}$ for all $i \neq j$). Then we define the action \eqref{eqn:action} by requiring $t \in \BC^*$ to rescale
\begin{align*}
&\phi_{e_c} \mapsto \phi_{e_c} t^{2-2c-d_{ij}} \\
&\phi_{e_c}^* \mapsto \phi_{e_c}^* t^{2c+d_{ij}} 
\end{align*}
for all $c \in \{1,\dots,-d_{ij}\}$. It is easy to see that $\BC^*$ preserves $\mu^{-1}_{\bn}(0)$, since the weights of $\phi_e$ and $\phi_e^*$ always add up to $2$. As the $\BC^*$ action above commutes with the $G_{\bn}$ action by conjugation, this ensures that the action \eqref{eqn:action} is well-defined.

\medskip

\begin{definition}
\label{def:k-ha}

The $K$-theoretic Hall algebra (with respect to the torus \eqref{eqn:action}) is
\begin{equation}
\label{eqn:k-ha int}
K_{\BC^*} = \bigoplus_{\bn \in \nn} K_{\BC^*}(T^*\fZ_{\bn})
\end{equation}
It is a module over $K_{\BC^*}(\emph{pt}) = \BZ[q^{\pm 1}]$, so we may also consider its localization
\begin{equation}
\label{eqn:k-ha loc}
K_{\BC^*, \eloc} = K_{\BC^*} \bigotimes_{\BZ[q^{\pm 1}]} \BQ(q)
\end{equation}
The multiplication on \eqref{eqn:k-ha int} and \eqref{eqn:k-ha loc} is given by a suitable convolution product (see \cite[Subsection 2.2]{VV}, as well as \cite[Subsection 2.3]{Quiver 1} for notation closer to ours).

\end{definition}

\medskip

\subsection{}
\label{sub:nilp}

Following \cite{VV}, we call a point 
$$
\left( \BC^{n_i} \xrightarrow{\phi_e} \BC^{n_j}, \BC^{n_j} \xrightarrow{\phi^*_e} \BC^{n_i} \right)_{\forall e = \oij \in E} \in T^*\fZ_{\bn}
$$
\textbf{nilpotent} if there exist flags of subspaces
$$
\Big\{ 0 = F^i_0 \subseteq F^i_1 \subseteq \dots \subseteq F^i_r = \BC^{n_i} \Big\}_{i \in I} 
$$
such that
\begin{equation}
\label{eqn:filtration}
\phi_e(F^i_p)\subseteq F^j_{p-1} \qquad \text{and} \qquad \phi_e^*(F^j_p) \subseteq F^i_{p-1}
\end{equation}
for all $p \in \{1,\dots,r\}$ and all arrows $e = \oij$. The set
\begin{equation}
\label{eqn:nilp substack}
\Lambda_{\bn} \subset T^*\fZ_{\bn}
\end{equation}
of nilpotent points is a closed substack.

\medskip

\begin{definition}
\label{def:k-ha nilp}

The $K$-theoretic Hall algebra supported on the substack \eqref{eqn:nilp substack} is
\begin{equation}
\label{eqn:k-ha int nilp}
K^{\emph{nilp}}_{\BC^*} = \bigoplus_{\bn \in \nn} K_{\BC^*}(T^*\fZ_{\bn})_{\Lambda_{\bn}}
\end{equation}
Also consider
\begin{equation}
\label{eqn:k-ha loc nilp}
K_{\BC^*, \eloc}^{\emph{nilp}} = K_{\BC^*}^{\emph{nilp}} \bigotimes_{\BZ[q^{\pm 1}]} \BQ(q)
\end{equation}
by analogy with \eqref{eqn:k-ha loc}. The convolution product that makes \eqref{eqn:k-ha int}--\eqref{eqn:k-ha loc} into algebras works equally well to make \eqref{eqn:k-ha int nilp}--\eqref{eqn:k-ha loc nilp} into algebras. The natural maps
$$
K^{\emph{nilp}}_{\BC^*} \rightarrow K_{\BC^*} \qquad \text{and} \qquad K_{\BC^*, \eloc}^{\emph{nilp}} \rightarrow K_{\BC^*, \eloc}
$$
are algebra homomorphisms.

\end{definition}

\medskip

\subsection{}
\label{sub:map to shuffle}

Consider the closed embedding of the origin
$$
0 \hookrightarrow T^*Z_{\bn}
$$
and let us consider the following composition
\begin{equation}
\label{eqn:iota 1}
\iota_{\bn} : K_{\BC^*}(T^*\fZ_{\bn})_{\Lambda_{\bn}} \xrightarrow{\text{natural map}} K_{\BC^* \times G_{\bn}}(T^*Z_{\bn}) \xrightarrow{\text{restriction to }0} K_{\BC^* \times G_{\bn}}(\text{pt})
\end{equation}
where the ``natural map" takes a coherent sheaf on $T^*\fZ_{\bn}$ supported on $\Lambda_{\bn}$ and interprets it as a $G_{\bn}$-equivariant coherent sheaf on the affine space $T^*Z_{\bn}$. One has
$$
K_{\BC^* \times G_{\bn}} (\pt) = \BZ[q^{\pm 1}][z_{i1}^{\pm 1}, \dots, z_{in_i}^{\pm 1}]^{\sym}_{i \in I}
$$
where $z_{ia}$ denotes the $a$-th elementary character of a maximal torus of $GL_{n_i}(\BC) \subset G_{\bn}$. Therefore, putting the maps \eqref{eqn:iota 1} together over all $\bn \in \BN$ yields a map
\begin{equation}
\label{eqn:iota 2}
\iota : K^{\text{nilp}}_{\BC^*} \longrightarrow \CV^{+,\text{geom}}_{\text{int}} := \bigoplus_{\bn \in \nn} \BZ[q^{\pm 1}][z_{i1}^{\pm 1}, \dots, z_{in_i}^{\pm 1}]^{\sym}_{i \in I}
\end{equation}
Tensoring everything with $\BQ(q)$ yields a map
\begin{equation}
\label{eqn:iota 3}
\iota_{\text{loc}} : K^{\text{nilp}}_{\BC^*,\loc} \longrightarrow \CV^{+,\text{geom}} := \bigoplus_{\bn \in \nn} \BQ(q)[z_{i1}^{\pm 1}, \dots, z_{in_i}^{\pm 1}]^{\sym}_{i \in I}
\end{equation}
Note that we could have defined the maps above with $K^{\text{nilp}}$ replaced by $K$ everywhere, simply by removing the subscript $\Lambda_{\bn}$ in \eqref{eqn:iota 1}. However, Theorem \ref{thm:k-ha} would cease to hold in this new setup, essentially because of the failure of Proposition \ref{prop:k-ha}.

\medskip

\subsection{}
\label{sub:zeta geom}

The maps $\iota$ and $\iota_{\text{loc}}$ are algebra homomorphisms, if the codomain of either map is endowed with the shuffle product \eqref{eqn:shuf prod}, but with $\zeta_{ij}(x)$ replaced by \footnote{Note that our rational function $\zeta_{ij}$ differs by an overall monomial from the function denoted by $\zeta_{ij}'$ in \cite{Quiver 1}; this discrepancy is innocuous, and we made it in order to match with the conventions in the present paper. On the $K$-theoretic Hall algebra side, the monomial in question can be implemented by appropriately twisting the convolution product, see \cite[Remark 2.5]{Quiver 1}.}
\begin{equation}
\label{eqn:zeta geometric}
\zeta_{ij}^{\text{geom}}(x) = \begin{cases} \displaystyle \frac {x-q^{-2}}{x-1} &\text{if } i = j \\ q^{-d_{ij}} \displaystyle \prod_{c = 0}^{-d_{ij}-1} \left(1-xq^{2c+d_{ij}} \right)&\text{if } i < j \\ \displaystyle \prod_{c = 1}^{-d_{ij}} \left(1-\frac {q^{2c+d_{ij}}}x \right)&\text{if } i > j
\end{cases}
\end{equation}
Because
\begin{equation}
\label{eqn:zeta to zeta}
\zeta_{ij}^{\text{geom}}(x) = \zeta_{ij}(x) \cdot \begin{cases} 1 &\text{if }i = j \\ \displaystyle (1-x)\prod_{c=1}^{-d_{ij}-1} \left(1 - x q^{2c+d_{ij}} \right) &\text{if } i < j \\ \displaystyle \left(1-\frac 1x\right)\prod_{c=1}^{-d_{ij}-1} \left(1 - \frac {q^{2c+d_{ij}}}x \right) &\text{if } i > j \end{cases}
\end{equation}
for all $i,j \in I$, it is easy to see that the linear map
\begin{equation}
\label{eqn:v geom}
\CV^+ \xrightarrow{\Omega} \CV^{+,\text{geom}}
\end{equation}
(we take $\BK = \BQ(q)$ in the definition of $\CV^+$ in \eqref{eqn:big shuf}) given for any $R \in \CV_{\bn}$ by
\begin{equation}
\label{eqn:big factor}
\Omega(R) = R \prod^{i < j \in I}_{1 \leq a \leq n_i, 1 \leq b \leq n_j} \left[ \left(1 - \frac {z_{ia}}{z_{jb}} \right) \prod_{c=1}^{-d_{ij}-1} \left(1 - \frac {z_{ia}q^{2c+d_{ij}}}{z_{jb}} \right) \right]
\end{equation}
is an algebra homomorphism. With this in mind, we conclude that
\begin{equation}
\label{eqn:s geom}
\CS^{+,\text{geom}} := \Omega(\CS^+)
\end{equation}
is a subalgebra of $\CV^{+,\text{geom}}$. Finally, define the subalgebra
\begin{equation}
\label{eqn:intersection}
\CS^{+,\text{geom}}_{\text{int}} = \CS^{+,\text{geom}} \cap \CV^{+,\text{geom}}_{\text{int}}
\end{equation}
of $\CV^{+,\text{geom}}_{\text{int}}$. 

\medskip

\begin{proposition}
\label{prop:prop}

A Laurent polynomial $r(\dots,z_{ia},\dots) \in \CV^{+,\emph{geom}}$ lies in $\CS^{+,\emph{geom}}$ if and only if
\begin{equation}
\label{eqn:prop}
(x-y)^{M_Z} \quad \text{divides} \quad \mathop{r\Big|_{z_{i1} = xq^{s}, z_{i2} = xq^{s+2}, \dots, z_{i,\frac {t-s}2+1} = xq^{t}}}_{\text{ }\qquad z_{j1} = yq^{s'}, z_{j2} = yq^{s'+2}, \dots, z_{j,\frac {t'-s'}2+1} = yq^{t'}}
\end{equation}
for any zig-zag $Z$ as in Figure 4, where 
\begin{equation}
\label{eqn:big number}
M_Z = m_Z + 
\end{equation}
$$
+ \Big|\Big\{(a,b,c) \in \{s,\dots,t\}_i \times \{s',\dots,t'\}_j \times \{1,\dots,-d_{ij}-1\} \text{ s.t. } a - b + 2c + d_{ij} = 0\Big\}\Big|
$$
and $m_Z$ is given by \eqref{eqn:number}. The analogous result holds for $\CS^{+,\emph{geom}}_{\emph{int}}$, defined by \eqref{eqn:intersection}.

\end{proposition}

\medskip

\begin{proof} We start with the ``if" statement. Letting $Z$ be the zig-zag consisting of one top vertex and one bottom vertex situated at distance $\in \{d_{ij}+2,\dots,-d_{ij}-2\}$ apart, property \eqref{eqn:prop} implies the fact that $r$ is divisible by
$$
1 - \frac {z_{ia}q^{2c+d_{ij}}}{z_{jb}}
$$
for all $i \neq j$, all $a,b$ and all $c \in \{1,\dots,-d_{ij}-1\}$. This implies that 
$$
r = \Omega(R)
$$
for some $R \in \CV^+$. If we further want $R \in \CS^+$ (which is equivalent to $r \in \CS^{+,\text{geom}}$) then we would need $R$ to satisfy conditions \eqref{eqn:wheel strong} for any zig-zag. Because of the linear factors involved in formula \eqref{eqn:big factor}, this is precisely equivalent to \eqref{eqn:prop}. The ``only if" statement is also established by the above argument. 

\end{proof}

\medskip

\subsection{}

Our main technical step in establishing Theorem \ref{thm:k-ha} is the following.

\medskip

\begin{proposition}
\label{prop:k-ha} 

The image of $\iota$ (resp. $\iota_{\emph{loc}}$) lies in $\CS^{+,\emph{geom}}_{\emph{int}}$ (resp. $\CS^{+,\emph{geom}}$).

\end{proposition}

\medskip

\begin{proof} We will prove the required statement for $\iota$, as the one for $\iota_{\text{loc}}$ is analogous. Throughout the present proof, we will fix a basis of $\{\BC^{n_i}\}_{i \in I}$ compatible with the maximal torus $T_{\bn} \subset G_{\bn}$ whose weights are the symbols $z_{i1},\dots,z_{in_i}$. Let us fix a zig-zag as in Figure 4, corresponding to $i < j$ in $I$, and we will only focus on a subset of the basis vectors
\begin{align}
&\BC^{n_i} = \dots \oplus \BC v_s \oplus \BC v_{s+2} \oplus \dots \oplus \BC  v_{t-2} \oplus \BC v_t \oplus \dots \label{eqn:basis 1} \\
&\BC^{n_j} = \dots \oplus \BC w_{s'} \oplus \BC w_{s'+2} \oplus \dots \oplus \BC  w_{t'-2} \oplus \BC w_{t'} \oplus \dots \label{eqn:basis 2}
\end{align}
We will consider the affine subspace $\BA \hookrightarrow T^*Z_{\bn}$ which parameterizes collections of linear maps $(\phi_e,\phi_e^*)_{e \in E}$, whose only non-zero matrix coefficients in the basis \eqref{eqn:basis 1}--\eqref{eqn:basis 2} are
\begin{align*} 
&\phi_{e_c}(v_\bullet) = x_{\bullet}^{\bullet'} \cdot w_{\bullet'}, \qquad \text{if } \bullet' = \bullet + 2c + d_{ij} - 2 \\
&\phi_{e_c}^*(w_{\bullet'}) = y_{\bullet'}^{\bullet} \cdot v_{\bullet}, \qquad \text{if } \bullet = \bullet' - 2c - d_{ij}
\end{align*}
for any $c \in \{1,\dots,-d_{ij}\}$, $\bullet \in \{s,\dots,t\}$ and $\bullet' \in \{s',\dots,t'\}$. Let us consider the following closed and open subsets of $\BA$, respectively 
$$
V = \BA \cap \mu_{\bn}^{-1}(0) \cap \Lambda_{\bn} \qquad \text{and} \qquad U = \BA \backslash V
$$
where $\mu_{\bn}$ is the quadratic map \eqref{eqn:moment map}, and $\Lambda_{\bn}$ is the closed subset of nilpotent points $(\phi_e,\phi_e^*)$ (we note a slight imprecision: here we are interpreting $\Lambda_{\bn}$ as living inside the affine space $T^*Z_{\bn}$, while in \eqref{eqn:nilp substack} it was defined as living inside $T^*\fZ_{\bn}$). 

\medskip

\begin{claim}
\label{claim:closed}

We may write $V = \cup_{k=1}^N W_k$, where each $W_k$ is the closure of the locus of $\BA$ determined by at least $M_Z$ independent equations of the form
$$
x_{\bullet}^{\bullet'} = \text{expression in other }x\text{'s and } y\text{'s}
$$
or
$$
y^{\bullet}_{\bullet'} = \text{expression in other }x\text{'s and } y\text{'s}
$$
Above, $M_Z$ is the number from \eqref{eqn:big number}.

\end{claim}

\medskip

\noindent Let us first show how Claim \ref{claim:closed} implies Proposition \ref{prop:k-ha}. First of all, specializing the variables $z_{ia}$ and $z_{jb}$ as in \eqref{eqn:prop} precisely corresponds to reducing the equivariance from the torus $\BC^* \times T_{\bn}$ to a subtorus $H$, where every matrix coefficient
\begin{align*}
&x_\bullet^{\bullet'} \text{ is scaled by the weight }\chi:H \rightarrow \BC^*,\text{and } \\
&y_{\bullet'}^{\bullet} \text{ is scaled by the inverse weight }\chi^{-1}
\end{align*}
where $\chi = \frac xy$ in the notation of \eqref{eqn:prop}. Consider the filtration of $V$ by closed subsets
$$
\varnothing = V_0 \subset V_1 \subset \dots \subset V_N = V \qquad \text{where} \qquad V_k = W_1 \cup \dots \cup W_k
$$
and let $U_k = \BA \backslash V_k$. For every $k$, we have the excision long exact sequence in equivariant algebraic $K$-theory
$$
K_H(U_{k-1} \cap W_k) \xrightarrow{\sigma_k} K_H(U_{k-1}) \xrightarrow{\tau_k} K_H(U_{k}) \rightarrow 0
$$
The assumption on $W_k$ implies that any element of $\text{Im }\sigma_k$ is a multiple of $(1-\chi)^{M_Z}$, and thus so is any element in $\text{Ker }\tau_k$. Iterating this claim $N$ times implies that any element of
\begin{equation}
\label{eqn:kernel final}
\text{Ker } \left( K_H(\BA) \xrightarrow{\tau_N \circ \dots \circ \tau_1} K_H(U) \right)
\end{equation}
is a multiple of $(1-\chi)^{M_Z}$. With this in mind, take any element $\iota_{\bn}(\alpha)$ and restrict it to the $H$-equivariant $K$-theory of the affine subspace $\BA$. Denote the resulting element $\beta$; the fact that $\alpha$ is supported on $\mu^{-1}_{\bn}(0) \cap \Lambda_{\bn}$ implies that $\beta$ is supported on the closed subset $V \subset \BA$, and thus $\beta$ lies in the kernel \eqref{eqn:kernel final}. As all elements in the kernel \eqref{eqn:kernel final} are multiples of $(1-\chi)^{M_Z}$, we conclude the sought-for divisibility condition in \eqref{eqn:prop}.

\medskip

\noindent It remains to prove Claim \ref{claim:closed}. To this end, let us draw two kinds of arrows between a vertex $a \in \{s,\dots,t\}$ and a vertex $b \in \{s',\dots,t'\}$:

\medskip

\begin{itemize}[leftmargin=*]

\item An arrow from $a$ to $b$ (respectively from $b$ to $a$) if $a = b-d_{ij}$ (respectively if $a = b + d_{ij}$); these will be called ``long" arrows.

\medskip 

\item An arrow from $a$ to $b$ and an arrow from $b$ to $a$ if $a-b \in \{d_{ij}+2,\dots,-d_{ij}-2\}$; these will be called ``short" arrows.

\end{itemize}

\medskip

\noindent Note that $M_Z$ of \eqref{eqn:big number} is equal to $m_Z+\upsilon$, where $\upsilon$ is the number of pairs of short arrows $a \leftrightarrow b$. For every pair of short arrows $a \leftrightarrow b$, note that one of the equations cutting out the closed subset $V \hookrightarrow \BA$ is
\begin{equation}
\label{eqn:xy}
x_a^b y_b^a = 0
\end{equation}
Otherwise, the fact that $x_a^b \neq 0$ and $y_b^a \neq 0$ simultaneously would violate the nilpotency condition (if $v_a \in F_p^i$ for some $p$, then $x_a^b \neq 0$ implies that $w_b \in F_{p-1}^j$, but $y_b^a \neq 0$ implies that $v_a \in F_{p-2}^i$; since this works for all $p$, we would obtain a contradiction). We conclude that $V$ is covered by the closure of the subsets
\begin{equation}
\label{eqn:defining equations}
W^\circ_k = \Big\{ f_a^b = 0 \text{ for every pair of short arrows } a \leftrightarrow b \Big\}
\end{equation}
as $k$ goes over the $2^{\upsilon}$ sets $(f_a^b \in \{x_a^b, y_b^a\})_{\text{pair of short arrows }a \leftrightarrow b}$. This already yields at least $\upsilon$ equations in Claim \ref{claim:closed}, but we need $M_Z = m_Z + \upsilon$ equations. To obtain this slightly bigger number of equations, let us use the long arrows; recall from the proof of Proposition \ref{prop:stronger} that the zig-zag $Z$ contains two intersecting long arrows
\begin{equation}
\label{eqn:two arrows}
_{b}\swarrow^{a} \qquad \text{ and } \qquad ^{a'}\nwarrow_{b'}
\end{equation}
and their $m_Z-1$ successive translates to the right by 2. Let us consider the following equations cutting out $V \hookrightarrow \BA$ (all these equations come from $\mu_{\bn}^{-1}(0)$)
\begin{align}
&x_a^b y_b^{a-2} + \dots = 0 \label{eqn:equation 1} \\
&x_{a-2}^b y_b^{a-4} + \dots = 0 \label{eqn:equation 2}  
\end{align}
$$
\dots
$$
\begin{align}
&x_{a'+2}^b y_b^{a'} + \dots = 0 \label{eqn:equation 3}  \\
&x_{a'}^b y_{b+2}^{a'} + \dots = 0 \label{eqn:equation 4}  
\end{align}
$$
\dots 
$$
\begin{align}
&x_{a'}^{b'-4} y_{b'-2}^{a'} + \dots = 0 \label{eqn:equation 5} \\
&x_{a'}^{b'-2} y_{b'}^{a'} + \dots = 0 \label{eqn:equation 6} 
\end{align}
(the ellipses stand for other sums of products of $x$'s and $y$'s, which will not matter in the subsequent argument). 

\medskip

\begin{claim}
\label{claim:cut out}

Besides the defining equations \eqref{eqn:defining equations}, the closed subset $V \cap W_k^\circ$ is determined by at least one more equation involving the variables in \eqref{eqn:equation 1}--\eqref{eqn:equation 6}.

\end{claim}

\noindent We may apply Claim \ref{claim:cut out} for each of the $m_Z$ pairs of long arrows given by \eqref{eqn:two arrows} and its translates. For each of these pairs, the ``one more equation" prescribed by Claim \ref{claim:cut out} always involves a variable with an index among the left-most endpoints (i.e. $a'$ and $b$) of the arrows \eqref{eqn:two arrows}. Since these endpoints are distinct for all the $m_Z$ pairs of long arrows, we conclude that the closed subset $V \cap W_k^\circ$ is determined by at least $m_Z$ more equations than the number $\upsilon$ of pairs of short arrows. As $M_Z = m_Z + \upsilon$, this concludes the proof of Claim \ref{claim:closed}.

\medskip

\begin{proof} \emph{of Claim \ref{claim:cut out}:} Among all the short pairs of variables 
\begin{equation}
\label{eqn:pairs}
(x_{a-2}^b, y_b^{a-2}), \dots (x_{a'+2}^b, y_b^{a'+2}), (x_{a'}^b, y_b^{a'}), \dots (x_{a'}^{b'-2}, y_{b'-2}^{a'})
\end{equation}
the set $W_k^\circ$ involves at least one of the variables being equal to 0. If from any such pair, the other variable were also equal to 0, then we would have obtained the required ``one more equation". Therefore, we may assume that among all the pairs of variables \eqref{eqn:pairs}, only one is equal to 0. However, this requires either

\medskip

\begin{itemize}[leftmargin=*]

\item $y_b^{a-2} \neq 0$, in which case \eqref{eqn:equation 1} gives us the extra equation $x_a^b = \dots$, or

\medskip

\item $x_{a'}^{b'-2} \neq 0$, in which case \eqref{eqn:equation 6} gives us the extra equation $y_{b'}^{a'} = \dots$, or

\medskip

\item at least one of the equations \eqref{eqn:equation 2}--\eqref{eqn:equation 5} being of the form
\begin{equation}
\label{eqn:equation}
x_{\bullet}^{\bullet'} y_{*}^{*'} + \dots \neq 0
\end{equation}
with both  $x_{\bullet}^{\bullet'} \neq 0$ and $y_{*}^{*'} \neq 0$, in which case \eqref{eqn:equation} yields an extra equation. 

\end{itemize}

\end{proof} \end{proof}

\medskip

\begin{proof} \emph{of Theorem \ref{thm:k-ha}:} By Proposition \ref{prop:k-ha}, we have an algebra homomorphism
\begin{equation}
\label{eqn:final iota}
\iota_{\loc} : K^{\text{nilp}}_{\BC^*, \text{loc}} \longrightarrow \CS^{+,\text{geom}}
\end{equation}
It remains to prove that the map above is surjective. However, $\CS^{+,\text{geom}} \cong \CS^+$ by \eqref{eqn:s geom}. Since Proposition \ref{prop:coincide} showed that $\CS^+$ is generated by $\{e_{i,k}\}_{i\in I, k\in \BZ}$ as a $\BQ(q)$-algebra, then the same is true of $\CS^{+,\text{geom}}$. However, all of the $e_{i,k}$'s are contained inside the image of $\iota_{\text{loc}}$ (when $\bn = \bs^i$, any quiver representation is nilpotent due to the absence of edge loops; the stack of quiver representations in this dimension is $\pt / \BC^*$, and its $K$-theory is precisely $\BZ[q^{\pm 1}][z_{i1}^{\pm 1}]$), hence all the generators of $\CS^{+,\text{geom}}$ are in the image of the map \eqref{eqn:final iota}. This implies that \eqref{eqn:final iota} is surjective. 

\end{proof}

\medskip

\subsection{} Theorem \ref{thm:k-ha} describes the image of the $K$-theoretic Hall algebra inside the shuffle algebra. To say something about the $K$-theoretic Hall algebra itself, we need to modify the definition of our quantum loop groups.

\medskip

\begin{definition}
\label{def:geom}

Define $\tUUp^{\emph{geom}}$ as in Subsection \ref{sub:quantum loop group intro}, but replacing \eqref{eqn:quad intro} by
\begin{equation}
\label{eqn:quad modified}
e_i(z)e_j(w) \prod_{c=0}^{-d_{ij}-1} (z-wq^{2c+d_{ij}}) = e_j(w)e_i(z) \prod_{c=0}^{-d_{ij}-1} (zq^{2c+d_{ij}}-w)
\end{equation}
for any $i \neq j$ (when $i = j$, relation \eqref{eqn:quad intro} is unchanged). 

\end{definition}

\medskip

\noindent Relation \eqref{eqn:quad modified} appeared in the context of $K$-theoretic Hall algebras in \cite{Nak} (it also arose in different contexts, e.g. \cite{J} in connection with vertex representations).

\medskip

\begin{proposition}
\label{prop:modified}

There exists a surjective algebra homomorphism
\begin{equation}
\label{eqn:k-ha modified}
\tUUp^{\emph{geom}} \twoheadrightarrow K^{\emph{nilp}}_{\BC^*, \emph{loc}}
\end{equation}
that sends $e_{i,k}$ to the $k$-th power of the tautological line bundle on $T^* \fZ_{\bs^i}$. 

\end{proposition}

\medskip

\begin{proof} Let us show that relation \eqref{eqn:quad modified} holds in $K^{\text{nilp}}_{\BC^*, \text{loc}}$. To this end, note that 
\begin{equation}
\label{eqn:stack}
T^*\fZ_{\bs^i+\bs^j} = \Big\{\Big( x_1,\dots,x_{-d_{ij}},y_1,\dots,y_{-d_{ij}} \Big) \text{ s.t. } \sum_{i=1}^{-d_{ij}} x_iy_i = 0 \Big\} \Big / \BC^* \times \BC^*
\end{equation}
The nilpotent substack $\Lambda_{\bs^i+\bs^j}$ consists of those tuples $(x_1,\dots,x_{-d_{ij}},y_1,\dots,y_{-d_{ij}})$ where either all the $x$'s are 0, or all the $y$'s are 0. Thus
\begin{equation}
\label{eqn:substack}
\Lambda_{\bs^i+\bs^j} = \Big( \BA^{-d_{ij}} \cup \BA^{-d_{ij}} \Big) \Big / \BC^* \times \BC^*
\end{equation}
with the two copies of $\BA^{-d_{ij}}$ intersecting at the origin. By the Mayer-Vietoris sequence (\cite[Appendix D]{VV}), the $K$-theory of the nilpotent substack is the direct sum of two copies of the $K$-theory of $\BA^{-d_{ij}}$, modulo the identification of the structure sheaf of the origin in the aforementioned two copies. Since the structure sheaf of the origin can be expressed as an equivariant Koszul complex, we conclude that
\begin{equation}
\label{eqn:quotient}
K_{\BC^*}(\Lambda_{\bs^i+\bs^j}) = \frac {\BZ[q^{\pm 1}, \ell_1^{\pm 1}, \ell_2^{\pm 1}] \bigoplus \BZ[q^{\pm 1}, \ell_1^{\pm 1}, \ell_2^{\pm 1}]}{\left( \prod_{c=0}^{-d_{ij}-1} \left(1 - \frac {\ell_1 q^{2c+d_{ij}}}{\ell_2} \right) -  \prod_{c=0}^{-d_{ij}-1} \left(1 - \frac {\ell_2 q^{2c+d_{ij}}}{\ell_1} \right) \right)}
\end{equation}
where $\ell_1$ and $\ell_2$ denote the elementary characters of the two factors of $\BC^* \times \BC^*$ in \eqref{eqn:substack}. The images of the LHS and RHS of \eqref{eqn:quad modified} in the $K$-theoretic Hall algebra are
$$
\delta\left(\frac w{\ell_1} \right) \delta\left(\frac z{\ell_2}\right) \prod_{c=0}^{-d_{ij}-1} \left(1 - \frac {w q^{2c+d_{ij}}}{z} \right)
$$
and
$$
\delta\left(\frac w{\ell_1} \right) \delta\left(\frac z{\ell_2}\right) \prod_{c=0}^{-d_{ij}-1} \left(1 - \frac {z q^{2c+d_{ij}}}{w} \right)
$$
respectively, and they are manifestly equal in the quotient \eqref{eqn:quotient}. While the formulas above seem like they differ from the LHS and RHS of \eqref{eqn:quad modified} by certain monomials, this discrepancy is precisely accounted for by line bundle twists in the formula for the $K$-theoretic Hall product (see \cite[Remark 2.5]{Quiver 1}). We do not recall the particular twist here, as it is precisely jury-rigged so that \eqref{eqn:quad modified} holds as stated.

\medskip

\begin{remark}

The stronger relation \eqref{eqn:quad intro} does not hold in $K_{\BC^*}(\Lambda_{\bs^i+\bs^j})$, as it is mapped not to 0, but to a zero-divisor.

\end{remark}

\medskip

\noindent Relation \eqref{eqn:quad intro} when $i=j$ holds for the same reason as \cite[Theorem A]{VV}, because the absence of edge loops means that the stack $\Lambda_{2\bs^i}$ is the same as in the case of the $\fsl_2$ quiver with one vertex and no arrows. Finally, the surjectivity of the map \eqref{eqn:k-ha modified} was proved in \cite[Theorem A]{VV}, see also the earlier unpublished work \cite[Section 7.10]{G}. In fact, the argument of \loccit proves the surjectivity of the following map
$$
\tUUp^{\text{geom}}_{\text{int}} \twoheadrightarrow K^{\text{nilp}}_{\BC^*}
$$
where the left-hand side is the $\BZ[q^{\pm 1}]$-subalgebra of $\tUUp^{\text{geom}}$ generated by
\begin{equation}
\label{eqn:divided powers}
\frac {e_{i,k}^n}{[n]!_q}
\end{equation}
for all $i \in I$, $k \in \BZ$ and $n \geq 0$ (recall that $[n]!_q = [1]_q \dots [n]_q$ with $[n]_q = \frac {q^n-q^{-n}}{q-q^{-1}}$). We refer the reader to \cite[Lemma 2.4.8]{VV} and \cite[Section 7.11]{G} for details. 

\end{proof}

\end{document}